\documentclass[A4,12pt]{amsart}

\usepackage{pstricks-add}

\usepackage[latin1]{inputenc}
\usepackage{amsmath,amsfonts,pstricks,bezier,mathrsfs}
\usepackage{hyperref,amssymb,latexsym,indentfirst}
\usepackage[T1]{fontenc}
\usepackage{xspace,ulem,cancel}
\usepackage{color}
\usepackage{upgreek} % letra grega nao fica em italico
\usepackage{fancyhdr}
\usepackage{fancybox}
\usepackage{float}   % H insere figura exatamente no local
\usepackage{setspace}
\usepackage{pstricks-add}
\usepackage{comment}

\hypersetup{
 pdftitle={OurPaper},
 pdfauthor={Nois},
 colorlinks=true,
 linkcolor=blue,
 citecolor=blue,
 filecolor=blue,
 urlcolor=blue}

\newtheorem{theorem}{Theorem} [section]

\newtheorem{definition}{Definition} [section]
\newtheorem{lemma}{Lemma}  [section]
\newtheorem{remark}{Remark}  [section]
\newtheorem{corollary}{Corollary} [section]

\newcommand{\R}{{\mathbb R}}

\newcommand{\N}{{\mathbb N}}

\newcommand{\C}{{\mathbb C}}

\def \ds {\displaystyle}

\newcommand{\cqd}{\hspace*{\fill}{\mbox{\framebox [1.5ex]{}}}}

\definecolor{lightgreen2}{rgb}{.4,.6,.4}
\definecolor{lightred2}{rgb}{.6,.6,.4}
\definecolor{lightbrown2}{rgb}{0.6,0.3,0.3}
\definecolor{lightgreen}{rgb}{.4,.6,.4}
\newrgbcolor{xdxdff}{0.49019607843137253 0.49019607843137253 1}
\newrgbcolor{xdxdff}{0.49 0.49 1}
\newrgbcolor{zzttqq}{0.6 0.2 0}
\newrgbcolor{cczzqq}{0.8 0.6 0}

\usepackage{graphicx}
\usepackage[format=plain,labelfont=bf,font=small]{caption}
\usepackage{subfigure}%replaced by subcaption
\usepackage{xcolor}
\usepackage[arrow, matrix, curve]{xy}
\usepackage{tikz}
\usepackage{float}

\usepackage{tikz-cd}
\usepackage{caption}
\usepackage{bm}

\title[Inertial Manifolds, Saddle Point Property and Dichotomy]{A Unified Theory for Inertial Manifolds, Saddle Point Property and Exponential Dichotomy}

\author[A. Carvalho]{Alexandre N. Carvalho}
\thanks{Partially supported by FAPESP grants 2020/14075-6 and CNPq 306213/2019-2}

\author[ P. Lappicy]{ Phillipo Lappicy}
\thanks{Partially supported by grants FAPESP 2017/07882-0 and 2018/18703-1}

\author[E. Moreira]{Estefani M. Moreira}
\thanks{Supported by CAPES grant PROEX 7547361/D and by FAPESP grants 2018/00065-9 and 2020/00104-4}

\author[A. Oliveira-Sousa]{ Alexandre N. Oliveira-Sousa}
\thanks{Supported by FAPESP grants 2017/21729-0 and 2018/10633-4 and by CAPES grant PROEX-9430931/D}

\address[A. N. Carvalho, P. Lappicy, E. M. Moreira, A. N. Oliveira-Sousa]{Instituto de Ci\^{e}ncias Ma\-te\-m\'{a}\-ti\-cas e de Computa\c{c}\~{a}o\\ 
Universidade de S\~{a}o Paulo-Campus de S\~{a}o Carlos \\	Caixa Postal 668, S\~{a}o Carlos SP, Brazil}
\email[A. Carvalho]{andcarva@icmc.usp.br}
\email[P. Lappicy]{lappicy@usp.br}
\email[E. Moreira]{estefani@usp.br}
\email[A. Oliveira-Sousa]{alexandrenosousa@gmail.com}

\subjclass[2020]{35B42, 37L45, 37D10, 37L25}

\keywords{Non-autonomous equations; Inertial manifolds; Hyperbolicity}

\begin{document}

%%%%%%%%%%%%%%%%%%%%%%%%%%%%%%%%%%%%%%%%%%%%%%%%%%%%%%%%%%%

\begin{abstract}
Inertial manifold theory, saddle point property and exponential dichotomy have been treated as different topics in the literature with different proofs. 
As a common feature, they all have the purpose of `splitting' the space to understand the dynamics.
We present a unified proof for the inertial manifold theorem, which as a local consequence yields the saddle-point property with a fine structure of invariant manifolds and the roughness of exponential dichotomy. In particular, we use these tools in order to establish the hyperbolicity of certain global solutions for non-autonomous parabolic partial differential equations.

\end{abstract}

\maketitle

\contentsline {section}{\tocsection {}{1}{Introduction}}{1}{section.1}%
\contentsline {section}{\tocsection {}{2}{Invariant Manifolds and their Stable Manifolds}}{5}{section.2}%
\contentsline {subsection}{\tocsubsection {}{2.1}{Invariant Manifolds}}{5}{subsection.2.1}%
\contentsline {subsection}{\tocsubsection {}{2.2}{Stable Manifold of an Invariant Manifold}}{12}{subsection.2.2}%
\contentsline {subsection}{\tocsubsection {}{2.3}{Local Inertial Manifolds}}{16}{subsection.2.3}%
\contentsline {section}{\tocsection {}{3}{Applications}}{18}{section.3}%
\contentsline {subsection}{\tocsubsection {}{3.1}{The Saddle Point Property}}{18}{subsection.3.1}%
\contentsline {subsection}{\tocsubsection {}{3.2}{Fine Description Within Invariant Manifolds }}{19}{subsection.3.2}%
\contentsline {subsection}{\tocsubsection {}{3.3}{Roughness of Exponential Dichotomy}}{21}{subsection.3.3}%
\contentsline {subsection}{\tocsubsection {}{3.4}{Hyperbolic Solutions for a Nonautonomous PDE}}{24}{subsection.3.4}%
\contentsline {section}{\tocsection {}{}{References}}{26}{section*.4}%

\section{Introduction}

%\textcolor{red}{PL: I suggest two things. First, we should give numbers to ALL equations. Second, instead of using equation+split, we should use subequations+align}

Inertial manifold theory, saddle point property and exponential dichotomy have been, so far, treated as different topics in the literature with distinct proofs. However, as a common feature, they all have the purpose of `splitting' the space in order to understand the dynamics.
The main goal of this paper is to provide a unified treatment to all these phenomena, including a broad range of applications in the autonomous and non-autonomous framework.

The inertial manifold theory has the purpose of reducing the relevant dynamics to an invertible dynamical system (often finite dimensional) and space splitting is implicitly established by the exponential attraction of the reduced dynamics. The proof presented here is inspired in the work of Henry \cite{HE} which in turn draws its inspiration from the book of Hale \cite{Hale-ODE}. The idea goes back to the work of Lyapunov \cite{Lyapunov} and Pliss \cite{Pliss} for the reduction principle when the splitting occurred at zero. The terminology of inertial manifold was introduced in the monumental work of Foias, Sell and Temam in \cite{FST88}, and there has been a tremendous effort in understanding these objects ever since. See \cite{Zelik} and references therein for a recent account. In particular, see \cite{KS} for inertial manifolds in the non-autonomous context. A central idea in the present paper is the introduction of the idea of an stable manifold of an inertial manifold, characterised by a graph, in contrast with the known results that prove only that the inertial manifold is exponentially attracting. 

%On one hand, the inertial manifold theory has a global nature, and its purpose is to reduce the dimension of the relevant dynamics to `finite' dimensions (or, more generally, to invertible dynamical systems), and space splitting is implicitly established by the exponential attraction of the reduced dynamics.
%On one hand, the inertial manifold theory has a global nature, and its purpose is to reduce the relevant dynamics to an invertible dynamical systems (often in finite dimensions), and space splitting is implicitly established by the exponential attraction of the reduced dynamics. The concept of inertial manifold was introduced in the monumental work of Foias, Sell and Temam in \cite{FST88}, and there has been a tremendous effort in understanding these objects ever since. See \cite{Zelik} and references therein for a recent account. In particular, see \cite{KS} for inertial manifolds in the non-autonomous context. A central idea in the present paper is the introduction of the idea of an stable manifold of an inertial manifold, characterized by a graph, in contrast with the known results that prove only that the inertial manifold is exponentially attracting. 

The saddle-point property has a local nature, and it arises as a nonlinear version of its linear exponential dichotomy counterpart, which splits the space into invariant linear subspaces that contain either exponential expansion or attraction.
Its study goes back to Perron, Massera and Sch\"{a}ffer \cite{perron,MS1,MS2,MS3} and many developments have been achieved ever since. See \cite{DK,Coppel,Coppel2,Scha,S-Y,YL,sacker-sell1,sacker-sell2,sacker-sell3,BV-Book,HE}.
In this direction, our results include: simple proofs of the saddle-point property and its robustness under perturbations; a fine description within the stable and unstable manifolds that allows us identify the possible and preferred directions of approximation (forwards and backwards in time) through an analysis of the growth/decay rates of solutions (new in the non-autonomous context); and an application to obtain that certain global solutions of non-autonomous parabolic partial differential equations are hyperbolic.
%
%Moreover, it is tied up with the concept of hyperbolicity, which is fundamental for the theory of dynamical systems. 

In order to properly introduce our results, some terminology is needed.
Consider the following semilinear differential initial value problem in a Banach space $(X,\|\cdot \|)$,
\begin{equation}\label{NLeq:glin}
\begin{split}
&\dot u = A(t)u + f(t,u),\ t>\tau,\\
&u(\tau)=u_0\in X,
\end{split}
\end{equation} 
where the map $f:\R\times X\to X$ is continuous, $f(t,0)=0$, for all $t\in \R$  
and uniformly Lipschitz in the second variable with Lipschitz constant $\ell>0$, i.e., $\|f(t,u)-f(t,\tilde{u})\|\leqslant \ell\|u-\tilde{u}\|$ for any $(t,u),(t,\tilde{u})\in \R\times X$. Assume that the family of linear operators $\{A(t):t\in \mathbb{R}\}$ (not necessarily bounded) defines a linear evolution process $\{L(t,\tau): t\geqslant \tau\}\subset\mathcal{L}(X)$, i.e., for each $(\tau,u_0)\in \mathbb{R}\times X$, the `solution' of the following linear problem,
\begin{equation}
\begin{split}
& \dot{u}=A(t)u, \ \ t\geqslant \tau,\\
& u(\tau)=u_0\in X,
\end{split}
\end{equation}
is given by $u(t,\tau,u_0)=L(t,\tau)u_0$, for $t\geqslant \tau$,  $L(t,t)=Id_X, L(t,s)L(s,\tau)=L(t,\tau)$, $t\geqslant s \geqslant \tau$ and $[\tau,\infty)\ni t\mapsto L(t,\tau)u_0\in X$ is continuous, for all $(\tau,u_0)\in \R\times X$. 

With this, solutions of \eqref{NLeq:glin} define a nonlinear evolution process $\{T(t,\tau): t\geqslant \tau\} \subset \mathcal{C}(X)$ given by the variation of constants formula, that is,
\begin{equation}\label{nonlinEP}
T(t,\tau)u=L(t,\tau)u+\int_\tau^t L(t,s)f(s,T(s,\tau)u)\,ds, \quad t\geqslant \tau, u\in X.
\end{equation}
%We assume that \eqref{NLeq:glin} defines a nonlinear evolution process $\{T(t,s):t\geqslant s\}$. 

Note that the hypothesis that the nonlinearity $f(t,\cdot)$ is globally Lipschitz and satisfy $f(t,0)=0$ are not restrictive for the intended analysis and are chosen in order to simplify the calculations, as we discuss in the upcoming remarks.
\begin{remark}\label{Rmk:shift}
The hypothesis that $f(t,0)=0$ can always be achieved, as one can translate any global solution to zero as follows.
In fact, assuming that $f(t,\cdot):X\to X$ is Frech\'{e}t differentiable and that if $u_*:\R\to X$ is a global solution of \eqref{NLeq:glin} and $\R\ni t\mapsto B(t):=f_u(t,u_*(t))\in \mathcal{L}(X)$ is strongly continuous, for a solution $u$ of \eqref{NLeq:glin}, we we note that $v(t):=u(t)-u_*(t)$ satisfies
%\begin{equation}\label{NLeq:glin0}
%\begin{split}
%&\dot v = A(t)v + f(t,u_*(t)+v)-f(t,u_*(t)),\ t>s\\
%&v(s)=v_0\in X
%\end{split}
%\end{equation}
%If $g$ is continuous and continuously differentiable with respect to the second variable, we can write
\begin{equation}\label{Eq_sl0}
\begin{split}
&\dot v = A(t)v +B(t)v+ g(t,v),\ t>\tau,\\
&v(\tau)=v_0\in X,
\end{split}
\end{equation}
where  $g(t,v):=f(t,u_*(t)+v)-f(t,u_*(t))-f_v(t,u_*(t))v$ and $v_0:=u_0-u_*(\tau)$. We observe that $g(t,0)=0\in X$ and $g_v(t,0)=0\in \mathcal{L}(X)$. In particular, when we wish to consider the local behaviour in a small tubular neighborhood of $u_*(\cdot)$ in $X$, it will be often possible to say that $g$ is Lipschitz with small constant $\ell>0$ in a small neighborhood of $0\in X$. In this case, it is possible to extend $g$ outside such small neighborhood of $0\in X$ such that it is globally Lipschitz with the same Lipschitz constant $\ell>0$, see for instance \cite[page 631]{Carvalho-Langa-07}.
\end{remark}
\begin{remark}
The hypothesis that $f(t,\cdot)$ is globally uniformly Lipschitz in the second variable is actually achieved after we cut-off a nonlinearity outside a region where it is Lipschitz, in order to describe the  behaviour of specific global solutions in that region.
Indeed, after the procedure that shifts global solutions to zero in Remark \ref{Rmk:shift}, one can often obtain a cut-off nonlinearity with sufficiently small $\ell>0$ on some neighborhood of zero, since $g(t,0)=0$ and $g_u(t,0)=0$. See \cite[Theorem 6.1]{Carvalho-Langa-07} for more details.
This cut-off operation will also be applied to bounded neighborhoods  of invariant sets, and in particular, bounded neighboorhoods of the global attractor.
\end{remark}
Let us define an inertial manifold for a nonlinear evolution process.
%
%\begin{definition}
%    Let $\{T(t,\tau): t\geqslant \tau\}$ be an evolution process and      $\{F(t): t\in \mathbb{R}\}$ be a family of nonempty subsets of $X$. We say that $\{F(t): t\in \mathbb{R}\}$ is
%    \begin{enumerate}
%        \item \emph{Forward attracting}, if $dist_H (T(t,\tau) B, F(t))\to 0$, as $t\to +\infty$, for any bounded set $B\subseteq X$. % and $s\in \mathbb{R}$,
        %\begin{equation*}
        %    dist_H (T(t,\tau) B, F(t))\stackrel{t\to +\infty}{\longrightarrow} 0.
        %\end{equation*}
%        \item \emph{Pullback attracting}, if $dist_H (T(t,\tau) B, F(t))\to 0$, as $\tau\to -\infty$, for any bounded set $B\subseteq X$. % and $t\in \mathbb{R}$,
        %\begin{equation*}
        %   dist_H (T(t,\tau) B, F(t))\stackrel{\tau\to -\infty}{\longrightarrow} 0.
        %\end{equation*}
%    \end{enumerate}
%Moreover, the set $\{F(t): t\in \mathbb{R}\}$ is called \textbf{exponentially attracting} if there are constants $\upsilon>0$ and $K:=K(B)>0$, for any bounded set $B\subseteq X$, such that $dist_H (T(t,\tau) B, F(t))\leqslant K(B)e^{-\upsilon(t-\tau)}$. \textcolor{red}{for which $t$? all $t$? or $t\geqslant \tau$? then why do we define forward and pullback attracting anyways?}
%\end{definition}
%
%
\begin{definition}\label{def-inertial-manifold}
A family $\{\mathcal{M}(t): t\in \mathbb{R}\}\subseteq X$ is called an \textbf{invariant manifold} for the evolution process $\{T(t,\tau): t\geqslant \tau\}$, if 
	\begin{enumerate}
		\item $\mathcal{M}(t)$ is a \emph{Lipschitz manifold}, for each $t\in \mathbb{R}$.
		\item $\{\mathcal{M}(t): t\in \mathbb{R}\}$ is \emph{invariant}, i.e., $T(t,\tau)\mathcal{M}(\tau)=\mathcal{M}(t)$, for $t\geqslant \tau$.
		\item $\{\mathcal{M}(t): t\in \mathbb{R}\}$ is \emph{exponentially dominated}, i.e.,
		\begin{enumerate}
            \item \emph{Forward exponentially dominated}: for any $\tau\in\mathbb{R}$ and bounded set $U\subseteq X$, there are $t_*{\color{black}\geqslant \tau},\delta_1\in \R,K>0$ such that $dist_H (T(t,\tau) U, \mathcal{M}(t))\leqslant Ke^{-\delta_1(t-\tau)} $, for  $t\geqslant t_*$.
            \item \emph{Pullback exponentially dominated}: for any $t\in\mathbb{R}$ and bounded set $U\subseteq X$, there are $\tau_*{\color{black}\leqslant t},\delta_2\in \R,K>0$ such that $dist_H (T(t,\tau) U, \mathcal{M}(t))\leqslant Ke^{-\delta_2(t-\tau)} $, for $\tau\leqslant \tau_*$.
            \end{enumerate}
		%\item There exists a $\delta>0$, such that for each $x\in X$, there is a 
		%$K=K(x)>0$ which
		%\begin{equation*}
		%dist_H(S(t,s)x,\mathcal{M}(t)) \leqslant K e^{-\delta(t-s)}, \  %t\geqslant s.
		%\end{equation*}
	\end{enumerate}
If $\delta_1>0$ and $\delta_2>0$, $\{\mathcal{M}(t): t\in \mathbb{R}\}\subseteq X$ is \emph{exponentially attracting} (pullback and forward) and it is called an \textbf{inertial manifold}.
\end{definition}
We prove in Theorem \ref{NLTum} the existence of an invariant manifold and its exponential attraction for the evolution process given by \eqref{nonlinEP}. We also  provide the characterization of the invariant manifold as a graph over a linear manifold. Moreover, in Theorem \ref{NLTsm}, we prove that the existence of a stable manifold of an invariant manifold in a general abstract setting, which again is given by a graph. 
Note that our hypothesis rely on the magnitude of the ratio between the exponential gap of the linear evolution process and the Lipschitz constant of the non-linearity, in contrast with the abstract conditions in \cite{KS}.
The main tool of our proof is an idea extracted from the existence of an invariant manifold contained in \cite[Chapter 6]{HE}. 
% which leads to the construction of the stable manifold of an exponentially attracting inertial manifold.
%Note that usually, the inertial manifold $\{\mathcal{M}(t): t\in \mathbb{R}\}$ and its stable manifold will be graphs over fixed (time independent) linear subspaces, which 
This leads to a reduction of \eqref{NLeq:glin} to an invertible flow, often on a finite dimensional subspace.
%\textcolor{red}{In particular, whenever the evolution process $\{T(t,\tau): t\geqslant \tau\}$ has a backwards bounded (i.e., the set $\cup_{t\leqslant 0}\mathcal{A}(t)$ is bounded) pullback attractor $\{\mathcal{A}(t): t\in \mathbb{R}\}$ and  an inertial manifold $\{\mathcal{M}(t): t\in \mathbb{R}\}$ then $\mathcal{A}(t) \subset \mathcal{M}(t)$ for all $t\in \R$.} 
As a consequence, we also obtain that the invariant manifold possess the asymptotic phase property in Corollary \ref{COR_ASYMPH}.

In particular, we  use a cut-off procedure in a neighborhood of the uniform attractor associated with \eqref{NLeq:glin}, which  yield in Theorem \ref{th-existence-local-inertial-manifolds-skew-products} a localized version of the inertial manifold. In order to achieve such `localization', we  use the skew-product semiflow associated with \eqref{NLeq:glin}.
%
%\par \textcolor{lightgreen}{We assume in Theorem \ref{NLTum} that the nonlinearity $f$ is globally Lipschitz. However, in the applications we usually have that $f$ is Lipschitz on bounded subsets of $X$. Thereforne, in Section \ref{sec-local-inertial-manifolds} we explain how to obtain an inertial manifold with this weaker assumption on $f$. The idea is to use a cut-off procedure in a neighborhood of an attracting set, which  be a uniform attractor associated with \eqref{NLeq:glin}. For this, we assume that \eqref{NLeq:glin} is associated with an skew product semiflow, in a partiuclar case where $A(t)\equiv A$. %We prove that the associated uniform attractor is stable, see Lemma \ref{lemma-uniform-attractor-stable}, which allow us to do a cut-off procedure in its neighborhood. 
%Then, assuming that the gap "$\gamma-\rho$" is suitable large, we prove the existence of a ``local inertial manifold'', see Theorem \ref{th-existence-local-inertial-manifolds-skew-products}. We apply this abstract result in Subsection \ref{sec:parabolicPDEs}}.

As particular cases of our study of the invariant manifold theory, we obtain the following:
\begin{itemize}

\item[(i)] The saddle point property (i.e. the characterization of the local stable and unstable manifolds as graphs of Lipschitz maps) in Corollary \ref{NLTexp_SPP}. Furthermore, from Theorem \ref{NLTum}, we obtain that the unstable manifold is exponentially attracting and we obtain explicit exponential growth/decay rates within the unstable/stable manifolds.
   
%\item[(ii)] \textcolor{lightgreen}{The roughness of exponential dichotomy for evolution process (possibly non-invertible), see Corollary \ref{NLTexp_dicho_roughness}. We prove a characterization of the linear stable and unstable manifolds of the perturbed problem as graphs of linear maps, leading to the robustness of hyperbolicity under perturbations. Our new proof is simpler than the monumental work of Henry in \cite{HE} because does not require a prior proof of admissibility in the discrete case. See also \cite{BV2} for a proof of robustness of nonuniform exponential dichotomy without using discrete case, which includes the case of (uniform) exponential dichotomies. For an introduction of nonuniform exponential dichotomies see \cite{BV-Book}.}

\item[(ii)] The roughness of exponential dichotomy for (possibly non-invertible) evolution processes in Corollary \ref{NLTexp_dicho_roughness}. Our new proof is simpler than the monumental work of Henry in \cite{HE} because it does not require a prior proof of admissibility in the discrete case.
See also \cite{BV1,BV2} for a different treatment of the robustness without passing through discrete case. However, our proof has the advantage of characterizing the linear stable and unstable manifolds of the perturbed problem as graphs of linear maps. %, leading to the robustness of hyperbolicity under perturbations. 
%For an introduction of nonuniform exponential dichotomies see \cite{BV-Book}.
%\textcolor{red}{Previously, we also cited \cite{BV-Book}... Do we need to cite it?}

%\textcolor{red}{The roughness of exponential dichotomy for the linearization around hyperbolic global solutions of autonomous and non-autonomous nonlinear differential equations (i.e., %the characterization of the linear stable and unstable manifolds of the perturbed problem as graphs of linear maps, leading to 
%    the robustness of hyperbolicity under perturbations) in Corollary \ref{NLTexp_dicho_roughness}. Our new proof on the roughness of exponential dichotomy in Corollary \ref{NLTexp_dicho_roughness} is simpler than the monumental work of Henry in \cite{HE} and does not require a prior proof of admissibility in the discrete case. We emphasize that it holds for the non-invertible case. See also \cite{BV1,BV2,BV-Book}.} \textcolor{red}{talk about the work of barreira}

\item[(iii)] A fine structure and linear approximation within invariant manifolds in Corollary \ref{CORFinestructure}, which extends the results \cite[Lemma 2.2]{BrunovskyFiedler86} and \cite[Lemma 6]{Angenent86} to certain non-autonomous equations. In particular, we are able to compare different asymptotic growth/decay rates within the unstable/stable manifolds of equation \eqref{NLeq:glin}. 
This enables us to identify the possible and the preferable directions of approximation (forwards and backwards) of the global hyperbolic solutions or hyperbolic equilibria.

\item[(iv)] {\color{black}If the problem under consideration has an invariant manifold, we may first project it into a finite dimensional subspace, to reduce the study of hyperbolicity or normal-hyperbolicity of a global solution to a finite dimensional problem, and then we use our specific knowledge of the finite dimensional problem to achieve it.} The hyperbolicity of certain global solutions, which are not obtained by means of perturbations  hyperbolic equilibria, for a non-autonomous one-dimensional scalar parabolic problem with localized large diffusion diffusion and reaction term given by $f(t,u)=u-\beta(t)u^3$ is considered in Corollary \ref{HypNonautEq}.  
\end{itemize}

\begin{remark}
As hyperbolicity is characterized by exponential dichotomy in the non-autonomous context and also for some autonomous infinite dimensional problems for which spectral separation does not necessarily induce a dichotomic behavior, see \cite[page 150]{Vesentini-02}.
\end{remark}

We believe that the results in this paper can be extended to the quasilinear or fully nonlinear case (given that some results on the saddle point property are already available in that case (see \cite[Section 9.1.2]{lunardi},\cite{yagi,Mielke91}). We also believe that the available results on the saddle point property for quasilinear and fully nonlinear problems can be improved requiring less regularity of the nonlinearity $f$ through the ideas introduced in this paper. This will be the subject of a future work.

A crucial notion used throughout the paper is the one of exponential splitting, that allows one to split a linear evolution process in two linear invariant subspaces related to growth/decay.
\begin{definition}\label{NLDexp-splitting}
A linear evolution process $\{L(t,\tau): t\geqslant \tau\}\subset {\mathcal L}(X)$ has {\bfseries exponential splitting}, with constant $M\geqslant 1$, exponents $\gamma,\rho\in \mathbb{R}$, with  $\gamma>\rho$,  and a family of projections $\{Q(t): \ t\in \R\}\subset {\mathcal L}(X)$, if 
\begin{enumerate}
\item[$i)$] $Q(t)L(t,\tau)=L(t,\tau)Q(\tau)$, for all $t\geqslant\tau$, 
\item[$ii)$]$L(t,\tau): Im(Q(\tau)) \to Im(Q(t))$ is an isomorphism, with inverse denoted by $L(\tau,t)$,
\item[$iii)$] the following estimates hold
\begin{equation}\label{NLEed-est}
\begin{split}
 \|L(t,\tau)Q(\tau)\|_{{\mathcal L}(X)}&\leqslant M e^{-\rho(t-\tau)},
\quad t\leqslant \tau,
\\
 \|L(t,\tau)(I-Q(\tau))\|_{{\mathcal L}(X)}&\leqslant M
e^{-\gamma(t-\tau)},\quad t\geqslant \tau.
\end{split}
\end{equation}
\end{enumerate}
In particular, if %$\{L(t,\tau): t\geqslant \tau\}$ has an exponential splitting with
$\gamma=-\rho$, then we say that $\{L(t,\tau): t\geqslant \tau\}$ has {\bfseries exponential dichotomy} with constant $M\geqslant 1$, exponent $\gamma>0$ and a family of  projections $\{Q(t): \ t\in \R\}$.
%with constant $M\geqslant 1$ and exponents $-\gamma=-\omega<-\rho=\omega=-\gamma<\beta$, $\gamma=\omega>0$, and projections $\{Q(t): \ t\in \R\}\subset {\mathcal L}(X)$ $($see \cite[Theorem I.15]{HES}$\,)$.

For $f(t,\cdot)\in C^1(X)$, a global solution $u_*:\mathbb{R}\to X$ of \eqref{NLeq:glin} is called {\bfseries hyperbolic} if the linear evolution process %$\{L_f(t,\tau): t\geqslant \tau \}\subset {\mathcal L}(X)$ 
given by
%\begin{equation}\label{EL-f}
$L_*(t,\tau):=L(t,\tau)+\int_\tau^t L(t,s) D_u f(s,u_*(s))
L_*(s,\tau)\,d s$
%\end{equation}
has exponential dichotomy.
\end{definition}
\begin{remark}
Note that the linear process $\{L(t,\tau): t\geqslant \tau\}\subset {\mathcal L}(X)$ has exponential splitting, with constant $M\geqslant 1$, exponents $\gamma,\rho\in \mathbb{R}$, with  $\gamma>\rho$, $\gamma>0$, and family of projections $\{Q(t): \ t\in \R\}\subset {\mathcal L}(X)$ if, and only if,
$\{e^{(\gamma -\frac{\gamma-\rho}{2})(t-\tau)}L(t,\tau): t\geqslant \tau\}\subset {\mathcal L}(X)$ has exponential dichotomy, with constant $M\geqslant 1$, exponent $\frac{\gamma-\rho}{2}>0$ and family of projections $\{Q(t): \ t\in \R\}\subset {\mathcal L}(X)$.
\end{remark}

Our aim with this paper is to prove the robustness of exponential splitting under linear and nonlinear perturbations and to show that this encompasses several topics that so far have been considered as different subjects in the literature, such as roughness of exponential dichotomy, the saddle point property, existence of invariant manifolds and the fine structure of solutions within the stable/unstable manifolds. All of this is achieved through the introduction of the stable manifold of an exponentially attracting invariant manifold.

The remaining of the paper is organized as follows. 
In Section \ref{Sin_man}, we present the main result that proves the existence of an invariant manifold in Theorem \ref{NLTum} and its stable manifold in Theorem \ref{NLTsm}, in case that $\gamma,\rho$ are real numbers in Definition \ref{NLDexp-splitting}. % with sufficiently large exponential gap $\gamma-\rho$. 
In Section \ref{Applications} we prove four consequences of the invariant manifold theorem: (i) the saddle point property in Corollary \ref{NLTexp_SPP}, (ii) a fine structure within stable and unstable manifolds in Corollary \ref{CORFinestructure}, (iii) the roughness of exponential dichotomies in Corollary \ref{NLTexp_dicho_roughness}, and (iv) the hyperbolicity of certain global solutions of a non-autonomous parabolic differential equation in Corollary \ref{HypNonautEq}. 
%In all consequences, the exponential gap $\gamma-\rho$ occurs at the imaginary axis such that $-\rho=\gamma>0$. 
The first three consequences occur when the Lipschitz constant $\ell>0$ is sufficiently small, and the last occurs when $\gamma-\rho$ is sufficiently large.

\section[Invariant Manifolds and Their Stable Manifolds]{Invariant Manifolds and their Stable Manifolds}\label{Sin_man}

\subsection{Invariant Manifolds}
The main result of this section gives sufficient conditions to ensure the existence of an invariant manifold. This condition basically states that the exponential gap $\gamma-\rho$ must be large when compared with the lipschitz constant $\ell$ of $f$. If $\gamma>0$ the invariant manifold will be an exponentially attracting inertial manifold.
\begin{theorem}\label{NLTum} 
Suppose that the linear evolution process $\{L(t,\tau): t\geqslant \tau\}$ has exponential splitting, with constant $M\geqslant 1$, exponents $\gamma>\rho$ and a family of projections $\{Q(t): \ t\in \R\}$. 
%\textcolor{red}{PL: my opinion is that we should call these projections $P(t)$, and the one that was $P(t)$ before, we should call $P_{\Sigma^*}(t)$.}
%
%If there is a constant $\kappa >0$ \textcolor{red}{PL: can we change $\kappa$ for something else, as it might be the laplacian} such that
%\begin{equation}\label{NLeq:inequalities}
%\begin{split}
%&\gamma-\rho-2LM(1+\kappa)>0\\
%&I_\Sigma:=\left[\frac{M^2L (1+M) }{\gamma-\rho-2LM(1+\kappa)} \right]= \nu<1, \quad  I_\eta:=\frac{M^2L (1+\kappa )}{\gamma-\rho-2LM(1+\kappa)} = \kappa\\
%&\delta:=\gamma-ML\kappa-\frac{M^2L^2(1+\kappa)(1+M)}{\gamma-\rho-ML(1+\kappa)}>0 
%\end{split}
%\end{equation}
If $f:\R\times X\to X$ is continuous, $f(t,0)=0$, $f(t,\cdot):X\to X$ is Lipschitz continuous with Lipschitz constant $\ell>0$, for all $t\in \R$, and
\begin{equation}\label{NLeq:inequalities}
    \frac{\gamma-\rho}{\ell} > \max \{ %{\color{black}2M(1+\kappa)^*} , 
    M^2+2M+\sqrt{8M^3}, 3M^2+2M\}, 
\end{equation}
then there is a continuous function 
\begin{equation}
\begin{aligned}\label{NLeq:prop_sigma}
    \Sigma^*:\mathbb{R}\times X &\to X\\
    (t,u)&\mapsto \Sigma^*(t,u)%=\Sigma^*(t,Q(t)u)=(I-Q(t))\Sigma^*(t,u),
\end{aligned}
\end{equation}
such that $\Sigma^*(t,u)=\Sigma^*(t,Q(t)u)=(I-Q(t))\Sigma^*(t,u)$ and
$\Sigma^*(t, 0)=0$, for all $t\in\mathbb{R}$. In addition $\Sigma^*(t,\cdot):X\to X$ Lipschitz continuous with Lipschitz constant $\kappa=\kappa(\gamma,\rho,\ell,M)>0$, for all $t\in \R$, that is, {\color{black}$\|\Sigma^*(t, u)-\Sigma^*(t, \tilde{u})\|\leqslant \kappa \|u-\tilde{u}\|$, for all $(t,u),(t,\tilde{u})\in \R\times X$}. 

Moreover, the graph of $\Sigma^*(t,.)$, for each $t\in \mathbb{R}$, given by
\begin{equation}
    \mathcal{M}(t):=\lbrace u\in X: u=q+\Sigma^*(t,q), q\in Im(Q(t)) \rbrace ,
\end{equation}
%\begin{equation}
%    \mathcal{M}(t):=\lbrace q+\Sigma(t,q)\in X: q\in Im(Q(t)) \rbrace,
%\end{equation}
yields an invariant manifold $\{\mathcal{M}(t): t\in \mathbb{R}\}$ for the evolution process $\{T(t,\tau): t\geqslant \tau\}$ given by \eqref{nonlinEP}. In other words, it is invariant and satisfies the following properties, according to the nonlinear projection $P_{\Sigma^*}(t)u:=Q(t)u+\Sigma^*(t,Q(t)u)$ onto $\mathcal{M}(t)$, for $(t,u)\in\mathbb{R}\times X$, 
\begin{itemize}
    \item[(i)] \emph{$\{\mathcal{M}(t): t\in \mathbb{R}\}$ has controlled growth}: for any $(\tau,u)\in \R\times X$, 
    \begin{equation}\label{invSigmaNL}
    \|T(t,\tau)P_{\Sigma^*}(\tau)u\| \leqslant M(1+\kappa)e^{-(\rho+\ell M (1+\kappa))(t-\tau)}\|P_{\Sigma^*}(\tau)u\|,\quad t\leqslant \tau.
    \end{equation}

    \item[(ii)]  \emph{$\{\mathcal{M}(t): t\in \mathbb{R}\}$ satisfies}: for any $(\tau,u)\in \R\times X$,
    \begin{equation}\label{expattSigmaNL}
    \| T(t,\tau)u-P_{\Sigma^*}(t)T(t,\tau)u\| \leqslant M\| (I-P_{\Sigma^*}(\tau))u\| e^{-\delta (t-\tau)}, \quad t\geqslant \tau,
    \end{equation}
    where  $ \delta:={\color{black}\gamma-M\ell -\frac{ M^2\ell^2(1+\kappa)(1+M)}{\gamma-\rho-\ell M (1+\kappa)}}$. If $ \delta>0$, $\{\mathcal{M}(t): t\in \mathbb{R}\}$ is an inertial manifold.
    
\end{itemize}
\end{theorem}
%
%\textcolor{red}{PL: it seems to me that $\delta>0$ in \eqref{deltaNL}  is not a hypothesis at all. since we already cut the spectrum: inside the inertial manifold we have either growth or decay, but outside the inertial manifold we MUST have decay by spectrum splitting, which means that we must prove that $\delta$ is positive, but we do not have to assume anything.}
%
\begin{remark}
When $\gamma$ is a sufficiently large positive number or the Lipschitz constant $\ell$ is sufficiently small, the condition that $\delta>0$ is satisfied. Note that the existence of the invariant manifold is independent of the sign of $\gamma$ and it only depends on the ratio  $\frac{\gamma-\rho}{\ell}$ between the exponential separation $\gamma-\rho$  and the Lipschitz constant $\ell$. In particular, $\kappa\to 0$ if $\frac{\gamma-\rho}{\ell} \to +\infty$.
\end{remark}
\proof
The proof is divided into two parts. First, we show that there is a function $\Sigma^*$ yielding the graph of the invariant manifold, as desired. Second, we show that this graph is exponentially dominated.

For the first part, given $\kappa>0$, consider the following complete metric space,
\begin{equation}\label{FPSigma}
\mathcal{LB}_\Sigma(\kappa) := \left\{ \Sigma\in C(\mathbb{R}\!\times\! X\!,\! X) :
\begin{array}{c}
 \ \sup_{t\in \R} \| \Sigma(t,u)\!-\!\Sigma(t,\tilde{u})\| \!\leqslant\!\kappa\|u\!-\!\tilde{u}\|,  \\
\Sigma(t,0)\!=\!0, \ \Sigma(t,u)\!=\!\Sigma(t,Q(t)u)= (I-Q(t))\Sigma(t,u),\ \forall t\in \R
\end{array}
\right\}.
\end{equation}
%\begin{equation}\label{FPSigma}
%\begin{split} 
%\mathcal{LB}_\Sigma(\kappa)
%\!=\!\Big\{\Sigma\!:\!\mathbb{R}\!\times\! X\!\to\! X&: \Sigma\hbox{ is continuous, } \Sigma(t,0)\!=\!0,\ 
%\| \Sigma(t,x)\!-\!\Sigma(t,x')\| \!\leqslant\!\kappa\|x\!-\!x'\|,  \\
%&\hbox{and } \Sigma(t,x)\!=\!(I-Q(t))\Sigma(t,Q(t)x)\!=\!\Sigma(t,Q(t)x),\ \forall t\in \R\Big\}.
%\end{split}
%\end{equation} 
with the metric $|\!|\!|\Sigma-\tilde{\Sigma}|\!|\!|:={\displaystyle\sup_{t\in \R}\sup_{u\neq 0} \frac{\|\Sigma(t,u)- \tilde{\Sigma}(t,u)\|}{\|u\|}}$.

We are looking for $\Sigma \in \mathcal{LB}_\Sigma(\kappa)$ with the property that, if $(\tau,\eta)\in \mathbb{R}\times X$, then a solution $u$ of 
\eqref{NLeq:glin} with initial data $u(\tau)=Q(\tau)\eta+\Sigma(\tau,Q(\tau)\eta)\in X$ can be decomposed as $u(t)=q(t)+p(t)$, where $p(t)=\Sigma(t,q(t))$ for all $t\in\mathbb{R}$. Thus, $q$ and $p$ must satisfy
\begin{subequations}\label{def:pq}
\begin{align}
q(t)&=L(t,\tau)Q(\tau)\eta +\int_\tau^t L(t,s)Q(s)f(s,q(s)+\Sigma (s,q(s)))ds, \ t\leqslant \tau,\label{solq(t)}\\
p(\tau)&=L(\tau,t)(I-Q(t))p(t) +\int_{t}^{\tau} L(\tau,s)(I-Q(s))f(s,q(s)+\Sigma(s,q(s))ds, \ t\leqslant \tau. \label{solp(t)}
\end{align}
\end{subequations}
%\textcolor{red}{I suggest we change the role of $\tau$ and $s$ throughout this proof. Moreover, beforehand we had $t\geqslant s$, and now we have $t\leqslant \tau$. It seems to me that you are solving $y(\tau)$ forwards, so, we should keep $t\geqslant s$, and write $y(t)$, whereas we should write $q(s)$, no? } \textcolor{lightgreen}{AS: I think in this situation the notation should not be confusing at all. $\tau$ is fixed real number, not initial or final initial time. We want to  get rid of $t$ in the equation for y to obtain the operator G, so it does not matter the previews notation ( t for final time and $s$ fo initial time).The important thing is that $\tau $ is fixed, and we need to solve the equation on $x$ for $t\leqslant \tau$ to define the operator $G$.} \textcolor{red}{PL: I understand this, but if you write $t\geqslant s$, and you fix $t$, so what we are consistent with the previous notation, you can write $q(s)$ and $y(t)$ above, and we have the \emph{same} meaning as the above equations, except with a different notation which is consistent with the previous one, that $t\geqslant s$; and in the equation for $y(t)$, we get rid of $s$. Hence, we fix $t$, we solve the equation $q(s)$ for $s\leqslant t$, and then we define the operator $G$ by getting rid of the 'parameter' $s$ in the equation for $y(t)$.}
%
\par First, let us control the growth of $q(t)$. Since 
$\{L(t,s): t\geqslant s\}$
has exponential splitting, $f(t,0)=0$ and $f$ and $\Sigma$ are Lipschitz with respective constants $\ell$ and $\kappa$, we obtain
\begin{equation}\label{eq-estamative-for-q}
\|q(t)\|
\leqslant  Me^{-\rho(t-\tau)}\|\eta\|  +\int_t^\tau \ell Me^{-\rho(t-s)}(1+\kappa) \|q(s)\|ds, \quad t\leqslant \tau. 
\end{equation}
\begin{comment}
For $\psi(t)=e^{\rho(t-\tau)}\|q(t)\| $ we have 
\begin{equation*} 
\psi(t) \leqslant M\|\eta\|  +\int_t^\tau LM(1+\kappa)\psi(s)ds.
\end{equation*}
\end{comment}
Then, by Gr\"onwall's Lemma,
\begin{equation}\label{eq-boundedness_for_Q(t)z(t)}
\|q(t)\| \leqslant Me^{(\rho+ M\ell(1+\kappa))(\tau-t)}\|\eta\|, \quad t\leqslant
\tau.
\end{equation}
Heuristically, since we wish that $p(t)=\Sigma(t,q(t))$ for 
$\Sigma\in \mathcal{LB}_\Sigma(\kappa)$,
the growth in equation \eqref{eq-boundedness_for_Q(t)z(t)} implies that the limit $e^{-\gamma (\tau-t)}\|p(t)\|\to 0$, as $t\to -\infty$. 
Thus, due to the exponential splitting of
$\{L(t,\tau):t\geqslant \tau\}$, the first term in \eqref{solp(t)} goes to zero as $t\to -\infty$, yielding
%\textcolor{red}{PL: add some extra heuristic arguments that allows us to formally pass this limit.}
%\textcolor{red}{PL: I'm confused $p(\tau)$ corresponds to the 'stable' part. if we let $t\to -\infty$, this should not be bounded, no? I think this sentence is still a bit weird...}
%\textcolor{lightgreen}{I forget to ask Alexandre C, we can ask on saturday, if I cannot remeber why. But in the case of Exponential Dichotomy I am able to explain, in fact the solutions p(t) goes to zero.}
\begin{equation}
    p(\tau)=\int_{-\infty}^{\tau} L(\tau,s)(I-Q(s))f(s,q(s)+\Sigma(s,q(s))ds.
\end{equation}
Hence, to rigorously obtain $\Sigma\in \mathcal{LB}_\Sigma(\kappa)$ that satisfies $p(\tau)=\Sigma(\tau,Q(\tau)\eta)$, it is equivalent to find a fixed point of the following map,
\begin{equation}\label{NLinvman7}
G(\Sigma)(\tau,\eta ) := \int_{-\infty}^\tau L(t,s)(I-Q(s))f(s,q(s)+\Sigma(s,q(s))) ds.
\end{equation}
Next, we show that $G:\mathcal{LB}_\Sigma(\kappa)\to \mathcal{LB}_\Sigma(\kappa)$ is a well defined contraction in the complete metric space $\mathcal{LB}_\Sigma(\kappa)$. 
%
%The main idea of the proof below constitutes in showing that $G$ is a well-defined map and more importantly, it is a contraction, which we embark now.
%

Let $\eta,\tilde{\eta}\in X$,
$\Sigma,\tilde{\Sigma}\in \mathcal{LB}_\Sigma(\kappa)$ with corresponding solutions $q(t),\tilde{q}(t)$ of \eqref{solq(t)}.
%For $\Sigma '\in \mathcal{LB}_\Sigma(\kappa)$ and $\tilde{\eta} \in X$, we consider $x'$ satisfying
%\begin{equation*}
%\begin{split}
%&q(t)=U(t,\tau)Q(\tau)\eta +\int_\tau^t U(t,s)Q(s)f(s,q(s)+\Sigma (s,q(s)))ds\ \hbox{ and }\\
%\tilde{q}(t) = U(t,\tau)Q(\tau)\tilde{\eta} +\int_\tau^t U(t,s)Q(s)f(s,\tilde{q}(s)+\tilde{\Sigma} (s,\tilde{q}(s)))ds, \qquad t\leqslant \tau.
%\end{split}
%\end{equation*}
%
%Then,
%\begin{equation*}
%q(t)\!-\!\tilde{q}(t)\!=\!U(t,\tau)Q(\tau)(\eta\!-\!\tilde{\eta})\!+\!\int_\tau^t\!\! U(t,s)Q(s)[f(s,q(s)\!+\!\Sigma (s,q(s)))\!-\!f(s,\tilde{q}(s)\!+\!\tilde{\Sigma}(s,\tilde{q}(s)))]ds.
%\end{equation*}
%And
%
Thus, for $t\leqslant \tau$, 
\begin{equation*}
\begin{split}
\|q(t)-\tilde{q}(t) \| \leqslant & M e^{{\rho}(\tau-t)}\|\eta\!-\!\tilde{\eta}\|\\
&+M\int_t^\tau e^{{\rho}(s-t)}\|f(s,q(s)+\Sigma (s,q(s)))-f(s,\tilde{q}(s)
+\tilde{\Sigma}(s,\tilde{q}(s)))\|ds \\
\leqslant & Me^{{\rho}(\tau-t)} \|\eta\!-\!\tilde{\eta}\| \\
&+ \ell M\int_t^\tau e^{-{\rho}(t-s)}
\left(\|q(s)-\tilde{q}(s) \|+\|\Sigma(s,q(s) )-\tilde{\Sigma}(s,\tilde{q}(s))\| \right) ds \\
\leqslant & M e^{{\rho}(\tau-t)}\|\eta\!-\!\tilde{\eta}\|\\
&+\ell M\!\!\int_t^\tau\!\! e^{{\rho}(s-t)}\!\left(\|\Sigma(s,q(s))\!-\!\tilde{\Sigma}(s,q(s) )\|
\!+\!(1\!+\!\kappa)\|q(s)-\tilde{q}(s)\| \right) ds \\
\leqslant & M e^{{\rho}(\tau-t)}\|\eta\!-\!\tilde{\eta}\| \\
&+ \ell M\int_t^\tau  e^{{\rho}(s-t)}\left((1+\kappa )\|q(s)-\tilde{q}(s)\|
+|\!|\!|\Sigma -\tilde{\Sigma}|\!|\!|\|q(s)\|\right) ds.
\end{split}
\end{equation*}
Then, due to \eqref{eq-boundedness_for_Q(t)z(t)},
\begin{equation}\label{eq-boundedness-xx'}
\begin{split}
\|q(t)-\tilde{q}(t) \| 
\leqslant &  M e^{{\rho}(\tau-t)}\|\eta\!-\!\tilde{\eta}\|\!+\! \ell M(1\!+\!\kappa)\!\!\int_t^\tau \!\! 
e^{{\rho}(s-t)}\|q(s)\!-\!\tilde{q}(s)\|ds\\
&+\ell M^2\|\eta\|{ |\!|\!|\Sigma -\tilde{\Sigma}|\!|\!|}\!\!
\int_t^\tau\!\! e^{(\rho+ M\ell(1+\kappa))(\tau-s)} e^{{\rho}(s-t)}ds\\
%\leqslant &  M e^{{\rho}(\tau-t)}\|\eta\!-\!\tilde{\eta}\|\!+\! \ell M(1\!+\!\kappa)\!\!\int_t^\tau\!\! 
%e^{{\rho}(s-t)}\|q(s)\!-\!\tilde{q}(s)\|ds\\
%&+\ell M^2\|\eta\|{ |\!|\!|\Sigma %-\tilde{\Sigma}|\!|\!|}
%e^{(\rho+ M\ell(1+\kappa))(\tau-t)}\int_t^\tau %e^{\ell M(1+\kappa)(t-s)}ds,\\
\leqslant &  M e^{{\rho}(\tau-t)}\|\eta\!-\!\tilde{\eta}\|\!+\! \ell M(1\!+\!\kappa)\!\!\int_t^\tau\!\! 
e^{{\rho}(s-t)}\|q(s)\!-\!\tilde{q}(s)\|ds\\
&+\frac{M\|\eta\|}{(1+\kappa)}{ |\!|\!|\Sigma -\tilde{\Sigma}|\!|\!|}e^{(\rho+ M\ell(1+\kappa))(\tau-t)},\\
\end{split}
\end{equation}
and, by Gr\"onwall's Lemma
\begin{equation}\label{eq-boundedness-x-x'}
    \|q(t)-\tilde{q}(t) \| \leqslant M\left[\|\eta -\tilde{\eta}\| +\frac{\|\eta\|}{1+\kappa} { |\!|\!|\Sigma -\tilde{\Sigma}|\!|\!|}\right]
e^{(\rho+2 M\ell(1+\kappa))(\tau-t)}, \quad t\leqslant \tau.
\end{equation}
%
%Making $\phi(t)=e^{(\rho+LM(1+\kappa))(t-\tau)}\|q(t)-\tilde{q}(t)\|$,
%\begin{equation*}
%\begin{split}
%\phi(t)\leqslant M\left[\|\eta -\tilde{\eta}\| +{\color{lightgreen}\frac{\|\eta\|}{1+\kappa}} { |\!|\!|\Sigma -\tilde{\Sigma}|\!|\!|}\right]+M\, L\,(1+\kappa)\int_t^\tau \phi(s) ds .
%\end{split}
%\end{equation*}
%
Finally, we now discuss bounds of the function $G$. Indeed, equations \eqref{eq-boundedness_for_Q(t)z(t)} and \eqref{eq-boundedness-x-x'} imply
\begin{equation*}
\begin{split}
\|G(\Sigma )(\tau,\eta ) - &G(\tilde{\Sigma} )(\tau, \tilde{\eta})\| \\
& \leqslant M\int_{-\infty}^\tau  e^{-\gamma(\tau - s)
}\|f(s,q(s)+\Sigma (s,q(s)))-f(s,\tilde{q}(s)
+\tilde{\Sigma}(s,\tilde{q}(s)))\|_{X }ds \\
%\leqslant & M\int_{-\infty}^\tau  e^{-\gamma (\tau - s)} \left( L\|\Sigma(s,q(s) )-\tilde{\Sigma}(s,\tilde{q}(s))\| + L \|q(s) -\tilde{q}(s)\| \right)ds \\
& \leqslant  \ell M \int_{-\infty}^\tau  e^{-\gamma (\tau - s)
} \left((1+\kappa )\|q(s) -\tilde{q}(s)\| + { |\!|\!|\Sigma -\tilde{\Sigma}|\!|\!| \|q(s)\|}\right)  ds \\
& \leqslant  \ell M^2 (1+\kappa ) \left[\|\eta -\tilde{\eta}\| +\frac{\|\eta\|}{1+\kappa} { |\!|\!|\Sigma -\tilde{\Sigma}|\!|\!|}\right] \int_{-\infty}^\tau  e^{-(\gamma-\rho-2 M\ell(1+\kappa)) (\tau - s)} ds
 \\
& \quad + \ell M^2\|\eta\| { |\!|\!|\Sigma -\tilde{\Sigma}|\!|\!|} \int_{-\infty}^\tau  e^{-(\gamma-\rho- M\ell(1+\kappa)) (\tau - s)}
    ds \, 
\end{split}
\end{equation*}
Due to \eqref{NLeq:inequalities} and upcoming choice of $\kappa$, we obtain that $\gamma-\rho-2\ell M(1+\kappa)\geqslant0$ and the above integrals are convergent. Thus,
\begin{equation*}
\begin{split}
\|G(\Sigma )(\tau,\eta ) &- G(\tilde{\Sigma} )(\tau, \tilde{\eta})\|\\ 
\leqslant & \frac{\ell M^2 (1+\kappa )}{\gamma-\rho-2\ell M(1+\kappa)}
\left[\|\eta -\tilde{\eta}\| + \frac{\|\eta\|}{1+\kappa} { |\!|\!|\Sigma -\tilde{\Sigma}|\!|\!|}\right] \\
&+ \frac{\ell M^2\|\eta\|}{\gamma-\rho-\ell M(1+\kappa)}  { |\!|\!|\Sigma -\tilde{\Sigma}|\!|\!|}   \\
%%
%\leqslant & \frac{M^2L (1+D )}{\gamma-\rho-2LM(1+D)}\|\eta -\tilde{\eta}\|
%\\
%&+\left[\frac{M^2L }{\gamma-\rho-2LM(1+D)}+ \frac{M^2L}{\gamma-\rho-LM(1+D)} \right]  { |\!|\!|\Sigma -\tilde{\Sigma}|\!|\!|\|\eta\|}  \\
%%
\leqslant & \frac{\ell M^2 (1+\kappa )}{\gamma-\rho-2\ell M(1+\kappa)}\|\eta -\tilde{\eta}\|
+\frac{2\ell M^2 }{\gamma-\rho-2\ell M(1+\kappa)}  { |\!|\!|\Sigma -\tilde{\Sigma}|\!|\!|\|\eta\|},
\end{split}
\end{equation*}
where the denominators are positive, due to \eqref{NLeq:inequalities}.
%%
%\begin{equation*}
%\begin{split}
%\leqslant & \frac{\ell M^2 (1+D )}{\gamma-\rho-2\ell M(1+D)}
%\left[\|\eta -\tilde{\eta}\| + \frac{\|\eta\|}{1+D} { |\!|\!|\Sigma -\tilde{\Sigma}|\!|\!|}\right] \\
%&+ \frac{\ell M^2\|\eta\|}{\gamma-\rho-\ell M(1+D)}  { |\!|\!|\Sigma -\tilde{\Sigma}|\!|\!|}   \\
%%
%\leqslant & \frac{M^2L (1+D )}{\gamma-\rho-2LM(1+D)}\|\eta -\tilde{\eta}\|
%\\
%&+\left[\frac{M^2L }{\gamma-\rho-2LM(1+D)}+ \frac{M^2L}{\gamma-\rho-LM(1+D)} \right]  { |\!|\!|\Sigma -\tilde{\Sigma}|\!|\!|\|\eta\|}  \\
%%
%\leqslant & \frac{\ell M^2 (1+D )}{\gamma-\rho-2\ell M(1+D)}\|\eta -\tilde{\eta}\|
%+\left[\frac{2\ell M^2 }{\gamma-\rho-2\ell M(1+D)} \right]  { |\!|\!|\Sigma -\tilde{\Sigma}|\!|\!|\|\eta\|} 
%\end{split}
%\end{equation*}
Consequently,
\begin{equation}\label{NLeq:propG}
\|G(\Sigma )(\tau, \eta ) - G(\tilde{\Sigma} )(\tau,\tilde{\eta})\|\leqslant \kappa \|\eta -
\tilde{\eta}\| + \nu|\!|\!|\Sigma -\tilde{\Sigma}|\!|\!|\|\eta\|,
\end{equation} 
in case that
\begin{subequations}\label{nuNL}
\begin{align}
    \frac{\ell M^2 (1+\kappa )}{\gamma-\rho-2\ell M(1+\kappa)} &\leqslant \kappa,\label{LipsG}\\
    \frac{2\ell M^2  }{\gamma-\rho-2\ell M(1+\kappa)} &<1.\label{contrG}
\end{align}
\end{subequations}
Equation \eqref{LipsG} can be rewritten as $2M\kappa^2+(M^2 + M - (\gamma-\rho)/\ell) \kappa+M^2\leqslant 0$. This can be seen as a quadratic polynomial (in $\kappa$), which admits two real roots due to \eqref{NLeq:inequalities}, given by\footnote{Note that $\kappa_\pm$ given by \eqref{Delta^*} satisfy $\kappa_-\to 0$ and $\kappa_+\to +\infty$, if either $(\gamma-\rho) \to +\infty$ or $\ell \to 0$.}
\begin{equation}\label{Delta^*}
    \kappa_\pm:= \frac{\frac{\gamma-\rho}{\ell}-M^2-2M\pm \sqrt{(\frac{\gamma-\rho}{\ell}-M^2-2M)^2-8M^3}}{4M}.
\end{equation}
Moreover, the condition $(\gamma-\rho)/\ell > M^2+2M+\sqrt{8M^3}$ in \eqref{NLeq:inequalities} implies that $(\gamma-\rho)/\ell > M^2+2M$ and thus $\kappa_+>\kappa_->0$. Thus, \eqref{LipsG} is satisfied for any $\kappa \in [\kappa_-,\kappa_+]$.
Equation \eqref{contrG} holds true for $\kappa_-$, due to $(\gamma-\rho)/\ell > 3M^2+2M$ in \eqref{NLeq:inequalities}. Moreover, we can isolate $\kappa$ in \eqref{contrG}, and thus this inequality is satisfied for any $\kappa < \kappa _*:=(\gamma-\rho)/(2M\ell) -M-1$.
Therefore, both conditions \eqref{nuNL} are satisfied for any $\kappa \in [\kappa_-,\min\{\kappa_+,\kappa_*\})$.
Note that $\kappa_-<\kappa_*$, due to the hypothesis $(\gamma-\rho)/\ell> 3M^2+2M$ in \eqref{NLeq:inequalities} and hence the interval $[\kappa_-,\min\{\kappa_+,\kappa_*\})$ is not empty.

Consequently, inequality \eqref{NLeq:propG} with $\Sigma=\tilde{\Sigma}$ implies that
the image of the map $G$ lies in $\mathcal{LB}_\Sigma(D)$, and thereby it is well-defined.
Similarly, inequality \eqref{NLeq:propG} with $\eta=\tilde{\eta}$ yields that 
$G$ is a contraction. 
Therefore, the map $G$ has a unique fixed point, $G(\Sigma^*)=\Sigma^*$. This establishes the existence of the invariant manifold and its invariance. 
Furthermore, $\Sigma^*$ being Lipschitz with constant $\kappa>0$ and $\Sigma(t,0)=0$, together with \eqref{eq-estamative-for-q}, implies the growth estimate \eqref{invSigmaNL} within the invariant manifold. This completes the first part of the proof.

%Proof of \eqref{invSigmaNL} is well defined.
%consider $q_0\in Im(Q(\tau))$ and the solution 
%$q(t)$ of 
%\begin{equation*}
%	q(t)=L(t,\tau)Q(\tau)q_0+\int_\tau^t L(t,r)Q(r)f(r,q(r)+\Sigma^u(r,q(r)) ) \, dr, \ \ t\geqslant \tau.
%\end{equation*}
%This defines a curve $q(t)+\Sigma(t,q^*(t)$, $t\in \mathbb{R}$. By the same reasoning of \eqref{NLinvman7} we have that
%\begin{equation*}
%	\Sigma^u(t,q(t)) := \int_{-\infty}^t L(t,s)(I-Q(s))f(s,q(s)+\Sigma^u(s,q(s))) ds, \ \ t\in \mathbb{R}.
%\end{equation*}
%Thus $q(t)+\Sigma(t,q^*(t)$, $t\in \mathbb{R}$ defines a global solution of $\{T(t,\tau): t\geqslant \tau \}$ passing through 
%$q_0+\Sigma^u(\tau,q_0)$ at time $\tau$, with 
%$\Sigma^u(t,q(t))\to 0$ as $t\to 0$. 
%From \eqref{eq-boundedness_for_Q(t)z(t)}
%\begin{equation}
%\|q(t)\|\leqslant Me^{(\gamma- M\ell(1+\kappa))(t-\tau)}\|q_0\|, \ \  t\leqslant \tau.
%\end{equation}
%Hence, $q(t)\to 0$ as $t\to -\infty$

We now embark in the second part of the proof. For $(t,u)\in \mathbb{R}\times X$, define the nonlinear projection
$P_{\Sigma^*}(t)u:=Q(t)u+\Sigma^*(t,u)$. We show that
$\mathcal{M}(t)=\{Im(P_{\Sigma^*}(t)): t\in \mathbb{R}\}$ is has the property that any solution satisfies \eqref{expattSigmaNL} (exponential attraction if $\delta>0$), %, i.e., there are positive constants $\delta,K>0$ such that
%\begin{equation}
%\|T(t,\tau)\eta - P(t)T(t,\tau)\eta \| \leqslant K e^{-\delta t}\|(I-P(\tau))\eta\|,
%\end{equation}
%for any $\eta\in X$ and $t\geqslant \tau$.
%\begin{equation*}
%q(t)=U(t,\tau)Q(\tau)\eta +\int_{\tau}^t U(t,s)Q(s)f(s,q(s)+p(s))ds
%\end{equation*}
and thus we wish to bound the variable $\xi  (t) := T(t,\tau)u - P_{\Sigma^*}(t)T(t,\tau)u$ for any $\eta\in X$ and $t\geqslant \tau$. 
Note that $\xi(t)=p(t)-\Sigma^*(t,q(t))$ due to the definitions in \eqref{def:pq}.

%Recall that $q(t):=Q(t)T(t,\tau)\eta$ and $p(t):=(I-Q(t))T(t,\tau)\eta$ in equations \eqref{solq(t)},\eqref{solp(t)}, and thus  $\xi(t)=p(t)-\Sigma^*(t,q(t))$.
Define $q^*(s,t)$, for $s\leqslant t$, as
\begin{equation}\label{defq*}
    q^*(s,t) :=L(s,t)q(t)+ \int_{t}^s L(s,r)Q(r)f(r,q^*(r,t)+\Sigma^* (r,q^*(r,t)))dr.
\end{equation}
Since $f,\Sigma^*$ are Lipschitz with respective constants $\ell,\kappa>0$, we obtain
\begin{equation}
\begin{split}
\|q^*(s,t)- q(s)\|%&= \Big\| \int_{t}^s L(s,r)Q(r)f(r,q^*(r,t)+\Sigma^* (r,q^*(r,t)) -f(r,q(r)+p(r))dr\Big\|, \\
&\leqslant M \int_s^t e^{\rho(r-s)}\|f(r,q^*(r,t)+\Sigma^* (r,q^*(r,t)))-f(r,q(r)+p(r))\| dr \\
&\leqslant  M\ell \int_s^t e^{\rho(r-s)} ( \|\Sigma^* (r,q^* (r ,t)) -
p(r)\| +
\| q^* (r ,t) - q(r)\|)dr, \\
&\leqslant M\ell \int_s^t  e^{\rho(r-s)}( \| \Sigma^* (r,q(r)) -
p(r)\| +
(1+\kappa )\| q^* (r ,t) - q(r)\|)dr, \ s\leqslant t. 
%&\leqslant ML\int_s^t  e^{\rho(r-s)} ( \|\xi (r )\|+ (1+\kappa ) \| q^* (r ,t) -q(r)\| )dr %.
\end{split}
\end{equation}

%Note that  $\xi(r)=p(r)-\Sigma^*(r,q(r))$, 
Hence, by Gronwall's Lemma and definition of $\xi$,
\begin{equation}\label{eq-3-14}
\|q^* (s,t)-q(s) \|\leqslant M\ell \int_s^t e^{(\rho+M\ell (1+\kappa)) (r-s)}\| \xi  (r)\| dr, \ \ s\leqslant t.
\end{equation}

Also, for $s\leqslant \tau \leqslant t$, we obtain %$\|q^* (s,t)-q^* (s,\tau)\|$. 
\begin{equation*}
\begin{split}
\|&q^* (s,t) -q^* (s,\tau)\|
\leqslant \|L(s,\tau)Q(\tau)[q^* (\tau, t) - q(\tau)]\| \\
&  \quad +\| {\ds\int_{\tau}^s} L(s,r)Q(r)[f(r,q^*(r,t)+\Sigma^* (r,q^*(r,t)))-f(r,q^*(r,\tau)+\Sigma^* (r,q^*(r,\tau))]\|dr\\
&  \leqslant   M^2\ell \, e^{\rho(\tau-s)} \int_{\tau}^t e^{(\rho+M\ell (1+\kappa)) (r-{\tau})}\| \xi  (r)\| dr + M\ell (1+\kappa ) {\ds\int_s^{\tau}}  e^{\rho(r-s)}
\|q^* (r,t) - q^* (r,\tau)\| dr,
\end{split}
\end{equation*}
and by Gr\"onwall's Lemma
\begin{equation}\label{eq-3-15}
\|q^* (s,t) -q^* (s,\tau)\|\leqslant M^2\ell \,  \int_{\tau}^t e^{(\rho+M\ell (1+\kappa))(r-s)}\| \xi  (r)\| dr.
\end{equation}
Now, we use these inequalities to estimate $\|\xi(t)\|$. Note that
\begin{equation*}
\begin{aligned}
\xi (t) &-L(t,\tau)(I-Q(\tau))\xi (\tau) = p (t) -L(t,\tau)p(\tau)\!-\! \Sigma^*(t,q(t)) \!+\!L(t,\tau)\Sigma^* (\tau,q(\tau))\\
&=\!\!{\ds\int_{\tau }^t} L(t,s)(I\!-\!Q(s))f(s,q(s)+p(s))ds\!-\!\! \int_{-\infty}^t\!\!\!\!
L(t,s)(I\!-\!Q(s))f(s,\Sigma^* (s, q^*(s,t))\!+\!q^* (s,t))ds \\
& \ \ \ +\int_{-\infty}^{\tau} L(t,s)(I-Q(s))f(s,q^* (s,\tau)+\Sigma^*(s,q^*(s,\tau)))ds \\
&={\ds\int_{\tau}^t} L(t,s)(I-Q(s))[f(s,q(s)+p(s))-f(s,q^* (s,t))+\Sigma^*(s,q^*(s,t)))]ds\\
& \ \ \ -{\ds \int_{-\infty}^{\tau}} L(t,s)(I-Q(s))[f(s, q^*(s,t)+\Sigma^* (s, q^*(s,t)))-f(s,q^* (s,\tau)+\Sigma^* (s, q^*(s,\tau))) ]ds.
\end{aligned}
\end{equation*}
    Thus, using \eqref{eq-3-14} and \eqref{eq-3-15},
we obtain
\begin{equation*}
\begin{split}
 \|\xi (t) &-L(t,\tau)(I-Q(\tau))\xi (\tau)\|\\
\leqslant & M\ell {\ds\int_{\tau}^t}  e^{-\gamma(t-s)} \left( \|p(s) - \Sigma^* (s,q^* (s,t))\| + \|q(s) - q^* (s,t)\| \right)ds \\
& + M\ell (1+\kappa )\int_{-\infty}^te^{-\gamma(t-s)}  \|q^* (s,\tau) - q^* (s,t)\| ds \\
\leqslant & M\ell {\ds\int_{\tau}^t} e^{-\gamma(t-s)}\|\xi  (s)\|ds +M^2\ell^2(1+\kappa){\ds\int_{\tau}^t}  
  e^{-\gamma (t-r)}\| \xi  (r)\|\int_\tau^re^{-(\gamma-\rho-M\ell (1+\kappa)) (r-s)} ds  dr \\
&+M^3\ell^2 (1+\kappa ){\ds  
 \int_\tau^te^{-\gamma(t-r)}e^{-(\gamma-\rho-M\ell (1+\kappa)) (r-\tau)}  \| \xi  (r)\| \int_{-\infty}^{\tau}} e^{-\gamma(\tau-s)}e^{(\rho+M\ell (1+\kappa)) (\tau-s)}ds dr\\
 \leqslant & M\ell {\ds\int_{\tau}^t} e^{-\gamma(t-s)}\|\xi  (s)\|ds +\frac{M^2\ell^2(1+\kappa)}{\gamma-\rho-M\ell (1+\kappa)} 
  {\ds\int_{\tau}^t}e^{-\gamma (t-r)}\| \xi  (r)\|  dr \\
&+\frac{M^3\ell^2 (1+\kappa )}{\gamma-\rho-M\ell (1+\kappa)}
{\ds\int_\tau^t } \| \xi  (r)\|e^{-\gamma(t-r)} dr
\end{split}
\end{equation*}
%
%Hence
%\begin{equation*}
%\begin{split}
%\|\xi (t) &-L(t,\tau)(I-Q(\tau))\xi (\tau)\|\leqslant  M\ell  {\ds\int_{\tau}^t} e^{-\gamma(t-s)}\|\xi  (s)\| \, ds  \\
%%
%%&
%%+ M^2\ell^2(1+\kappa){\ds\int_{\tau}^t} e^{-\gamma (t-r)}  \int_s^t e^{-(\gamma -\rho-M\ell (1+\kappa)) (r-s)}\| \xi  (r)\| dr ds\\
%%%
%& + M^2\ell^2(1+\kappa){\ds\int_{\tau}^t} 
%\| \xi  (r)\|e^{-\gamma(t-r)}\int_{\tau}^r e^{-(\gamma-\rho-M\ell (1+\kappa)) (r-s)}   \, ds \ dr\\
%& + M^3\ell^2 (1+\kappa )e^{-(\gamma-\rho-M\ell(1+\kappa))(t-\tau)}
%\int_{\tau}^t e^{\rho (r-{t})}\| \xi  (r)\|
%{\ds \int_{-\infty}^{\tau}}  e^{-(\gamma-\rho-M\ell (1+\kappa))(\tau-s)}
%  \, ds \, dr\\
%\leqslant & M\ell {\ds\int_{\tau}^t} e^{-\gamma(t-r)}\|\xi  (r)\|dr \\
%& + \frac{M^2\ell^2(1+\kappa)}{\gamma-\rho-M\ell(1+\kappa)}{\ds\int_{\tau}^t} 
%e^{-(\gamma-ML(1+\kappa)) (t-r)} \| \xi  (r)\| \, dr\\
%& + \frac{M^3\ell^2 (1+\kappa )}{\gamma-\rho-M\ell(1+\kappa)}e^{-(\gamma-\rho-M \ell(1+\kappa))(t-\tau)}
% \int_{\tau}^t e^{\rho (r-{t})}\| \xi  (r)\| dr,
%\end{split}
%\end{equation*}
and we have that
\begin{equation*}
%& \leqslant M\ell {\ds\int_{\tau}^t} e^{-\gamma(t-r)}\|\xi  (r)\|dr +\frac{M^2\ell^2(1+\kappa)(1+M)}{\gamma-\rho-M\ell(1+\kappa)}{\ds\int_{\tau}^t} 
 %e^{-(\gamma-M\ell(1+\kappa)) (t-r)} \|\xi(r)\|dr\\
 \|\xi (t) -L(t,\tau)(I-Q(\tau))\xi (\tau)\| \leqslant   \left[M\ell+\frac{M^2\ell^2(1+\kappa)(1+M)}{\gamma-\rho-M\ell(1+\kappa)}\right] {\ds\int_{\tau}^t} 
 e^{-\gamma (t-r)} \|\xi(r)\|dr.\\
\end{equation*}
Thus, 
\begin{equation*}
\|\xi (t)\| \leqslant M e^{-\gamma(t-\tau)} \|\xi  (\tau)\|
+  \left[M\ell +\frac{M^2\ell^2(1+\kappa)(1+M)}{\gamma-\rho-M\ell (1+\kappa)}\right]{\ds\int_{\tau}^t} 
 e^{-\gamma (t-r)} \|\xi(r)\|dr
\end{equation*}
and
\begin{equation*}
e^{\gamma t}\|\xi (t)\| \leqslant  M e^{\gamma\tau} \|\xi  (\tau)\|+ \left[M\ell +\frac{M^2\ell^2(1+\kappa)(1+M)}{\gamma-\rho-M\ell (1+\kappa)}\right] {\ds\int_{\tau}^t} 
 e^{\gamma r} \|\xi(r)\|dr.
\end{equation*}
By Gr\"onwall's Lemma, we obtain the bound in equation \eqref{expattSigmaNL}.
%$$ \|\xi (t)\| \leqslant M e^{-\delta (t-\tau)} \|\xi  (\tau)\| $$
%where
%$$\delta:=\gamma-M\ell \kappa-\frac{M^2\ell^2(1+\kappa)(1+M)}{\gamma-\rho-M\ell (1+\kappa)}>0.$$
%The proof is now complete. 
\cqd

\begin{remark}
Continuity of the invariant (inertial if $\delta>0$) manifold with respect to perturbations on $f$, and thus continuity with respect to parameters, follows as usual; see \cite{Carvalho-Langa-07} for a detailed proof in the case that $\rho<0$.
In contrast with the usual regularity issues for stable/unstable manifolds, when the $C^k$ regularity of the vector field implies the $C^k$ regularity of these invariant manifolds, the regularity of invariant (inertial if $\delta>0$) manifolds is more delicate, see \cite{CLS92,KZ}. However, such regularity can be achieved by stronger exponential gap conditions.
%Also, we expect that the inertial manifold is of same regularity of the vector field $f$, i.e., $f\in C^k$ implies $\mathcal{M}(t)\in C^k$ by means of an appropriate cut-off. See  \textcolor{red}{[add reference.. ]}. \textcolor{red}{ this is not true... the paper of zelik says is $f\in C^2$, the IM is only holder}
Moreover, regularity for the invariant (inertial if $\delta>0$) manifold in the non-autonomous case can be achieved in a similar manner as in \cite[Chapter 4]{BCL}.
These considerations are also valid for the upcoming results.
\end{remark}

Note that, if $\delta>0$, the exponentially attracting property \eqref{expattSigmaNL} of the inertial manifold in Theorem \ref{NLTum} shows that any solution in phase-space converges to some set within the inertial manifold. %, which may not necessarily consist of solutions.
If, in addition, we assume that the projections associated with the dichotomy have finite rank, then any solution converges to a solution within the inertial manifold, and thereby the inertial manifold also possess the asymptotic phase property described below in a more general context in which $\delta$ does not need to be positive.

\begin{corollary}\label{COR_ASYMPH}
Consider the hypothesis of Theorem \ref{NLTum} and that $\mathrm{dim}\, Im (Q(\tau))<\infty$. Then, for any $u_0\in X$, there exists a global solution $\bar{q}(t)+\Sigma^*(t,\bar{q}(t))\in \mathcal{M}(t)$ of \eqref{NLeq:glin} and a constant $c>0$ such that
%$$
%\|T(t,t_0)u_0-\bar{q}(t)-\Sigma^*(t,\bar{q}(t))\|_X\leqslant c e^{-\gamma (t-t_0)} \|(I-Q(t_0))u_0-\Sigma^*(t_0,\bar{q}(t_0)))\|_X.
%$$
\begin{equation}
\begin{split}
\|T(t,\tau)u_0-\bar{q}(t)-\Sigma^*(t,\bar{q}(t))\|\leqslant c e^{-\delta  (t-\tau)} \|(I-Q(\tau))u_0-\Sigma^*(\tau,\bar{q}(\tau)))\|.
\end{split}
\end{equation}
\end{corollary}

\proof 
Since $f,\Sigma^*$ are Lipschitz with respective constants $\ell,\kappa>0$, we obtain 
\begin{equation*}
\begin{split}
\|q^*(s,t)- q(s)\|
%&= \Big\| \int_{t}^s L(s,r)Q(r)f(r,q^*(r,t)+\Sigma^* (r,q^*(r,t)) -f(r,q(r)+p(r))dr\Big\|, \\
&\leqslant M \int_s^t e^{\rho(r-s)}\|f(r,q^*(r,t)+\Sigma^* (r,q^*(r,t)))-f(r,q(r)+p(r))\| dr \\
%&\leqslant  M\ell \int_s^t e^{\rho(r-s)} ( \|\Sigma^* (r,q^* (r ,t)) -
%p(r)\| +
%\| q^* (r ,t) - q(r)\|)dr, \\
&\leqslant M\ell \int_s^t  e^{\rho(r-s)}(\| \xi  (r)\|+
(1+\kappa )\| q^* (r ,t) - q(r)\|)dr,\\
%&\leqslant M\ell \int_s^t e^{\rho (r-s)}\| \xi  (r)\| dr{\color{black} e^{M\ell (1+\kappa) (t-s)}},
%&\leqslant ML\int_s^t  e^{\rho(r-s)} ( \|\xi (r )\|+ (1+\kappa ) \| q^* (r ,t) -q(r)\| )dr %.
\end{split}
\end{equation*}
hence
\begin{equation}
\begin{split}
\|q^*(s,t)- q(s)\|%&= \Big\| \int_{t}^s L(s,r)Q(r)f(r,q^*(r,t)+\Sigma^* (r,q^*(r,t)) -f(r,q(r)+p(r))dr\Big\|, \\
%&\leqslant M \int_s^t e^{\rho(r-s)}\|f(r,q^*(r,t)+\Sigma^* (r,q^*(r,t)))-f(r,q(r)+p(r))\| dr \\
%&\leqslant  M\ell \int_s^t e^{\rho(r-s)} ( \|\Sigma^* (r,q^* (r ,t)) -
%p(r)\| +
%\| q^* (r ,t) - q(r)\|)dr, \\
%&\leqslant M\ell \int_s^t  e^{\rho(r-s)}( \| \Sigma^* (r,q(r)) -
%p(r)\| +
%(1+\kappa )\| q^* (r ,t) - q(r)\|)dr,\\
&\leqslant M\ell \int_s^t e^{(\rho+M\ell (1+\kappa) ) (r-s)}\| \xi  (r)\| dr ,
%&\leqslant ML\int_s^t  e^{\rho(r-s)} ( \|\xi (r )\|+ (1+\kappa ) \| q^* (r ,t) -q(r)\| )dr %.
\end{split}
\end{equation}
and that
\begin{equation}
\|\xi(r)\| \leqslant M\| \xi(s)\| e^{-\left(\gamma-M\ell -\frac{ M^2\ell^2(1+\kappa)(1+M)}{\gamma-\rho-\ell M (1+\kappa)}\right) (r-s)}, \quad r\geqslant s.
\end{equation}
and
\begin{equation}\label{bar{q}est-AC}
\begin{split}
\|q^*(s,t)- q(s)\|
&\leqslant M^2\ell \int_s^t  e^{-\left(\gamma-\rho-M\ell (2+\kappa)-\frac{ M^2\ell^2(1+\kappa)(1+M)}{\gamma-\rho-\ell M (1+\kappa)}\right) (r-s)} dr \| \xi(s)\|,
\end{split}
\end{equation}

%%%%%%%%%%%%%%%%%%%%%%%%%%%%%%%
%

Since $\{q^*(\tau,t):t\geqslant \tau\}$ is bounded in $Im(Q(\tau))$, there is a sequence $t_n\to \infty$ and $\bar{q}(\tau) \in  Im(Q(\tau))$ such that $q^*(\tau,t_n)\stackrel{n\to \infty}{\longrightarrow} \bar{q}(\tau)$. Then, for every $s\geqslant \tau$, we have that $q^*(s,t_n)\stackrel{n\to \infty}{\longrightarrow}\bar{q}(s) $ and thereby $\bar{q}$ defines a solution in the invariant manifold, i.e.,
\begin{equation}
\begin{split}
&\bar{q}(s) :=L(s,\tau)\bar{q}(\tau)+ \int_{\tau}^s L(s,r)Q(r)f(r,\bar{q}(r)+\Sigma^* (r,\bar{q}(r)))dr, \ s\leqslant \tau\\
&
\Sigma^*(s,\bar{q}(s)) := \int_{-\infty}^s L(\tau,\theta)(I-Q(\theta))f(\theta,\bar{q}(\theta)+\Sigma(\theta,\bar{q}(\theta))) d\theta,\ s\leqslant \tau.
\end{split}
\end{equation}

From \eqref{bar{q}est-AC}, we deduce that
\begin{equation}
\|\bar{q}(s)- q(s)\|\leqslant  \frac{M^2\ell \| \xi(\tau)\|}{\gamma-\rho-M\ell (2+\kappa)-\frac{ M^2\ell^2(1+\kappa)(1+M)}{\gamma-\rho-\ell M (1+\kappa)}}   e^{-\left(\gamma-M\ell -\frac{ M^2\ell^2(1+\kappa)(1+M)}{\gamma-\rho-\ell M (1+\kappa)}\right) (s-\tau)}.
\end{equation}
The remaining part of the proof follows directly from Theorem \ref{NLTum}.
%\textcolor{red}{what remains to be proved?}
\cqd

\subsection[Stable Manifold of an Invariant Manifold]{Stable Manifold of an Invariant Manifold}\label{NLSum-stable}

\begin{theorem}\label{NLTsm}
Suppose that the linear evolution process $\{L(t,\tau): t\geqslant \tau\}$ has exponential splitting, with constant $M\geqslant 1$, exponents $\gamma>\rho$ and a family of projections $\{Q(t): \ t\in \R\}$. 

If $(\gamma-\rho)/\ell$ satisfies \eqref{NLeq:inequalities}, then there is a continuous function 
\begin{equation}
\begin{aligned}\label{NLeq:prop_theta}
    \Theta^*:\mathbb{R}\times X &\to X\\
    (t,u)&\mapsto \Theta^*(t,u),
\end{aligned}
\end{equation}
such that $\Theta^*(t,u)=\Theta^*(t,(I-Q(t))u)=Q(t)\Theta^*(t,u)$, 
and $\Theta^*(t, 0)=0$ for all $t\in\mathbb{R}$, which is uniformly Lipschitz with constant $\kappa=\kappa(\gamma,\rho,\ell,M)>0$ , i.e., $\|\Theta^*(t, u)-\Theta^*(t, \tilde{u})\|\leqslant \kappa \|u-\tilde{u}\|$ for all $(t,u),(t,\tilde{u})\in \R\times X$. 

Moreover,  if  $P_{\Theta^*}(t)u:=\Theta^*(t,(I-Q(t))u)+(I-Q(t))u$, for all $(t,u)\in \R\times X$, the family given by
\begin{equation}
\left\{Im(P_{\Theta^*}(t)): t\in\R\right\}:=\{\left\{ P_\Theta^*(t,u): u\in X\} : t\in\R \right\},
\end{equation}
is positively invariant such that
\begin{equation}\label{DecayStM}
    \|T(t,\tau)P_{\Theta^*}(\tau)u\| \leqslant M(1+\kappa)e^{-(\gamma- M\ell(1+\kappa))(t-\tau)}\|P_{\Theta^*}(\tau)u\|, \quad t\geqslant \tau, \ u\in X,
\end{equation}
 and{\color{black}
\begin{equation}\label{attraction_stable}
\|u- P_{\Theta^*}(\tau)u \| \leqslant M e^{\hat{\delta}(t- \tau)}\|(I-P_{\Theta^*}(\tau))T(t,\tau)u\|, \quad t\geqslant \tau, u\in X,
\end{equation}
where $\hat{\delta}=\rho + M\ell +\frac{M^2\ell^2(1+\kappa)(1+M)}{\gamma-\rho-M\ell (1+\kappa)}$.

Furthermore, if $\gamma- M\ell(1+\kappa)>0$, $\left\{Im(P_{\Theta^*}(t)): t\in\R\right\}$ is the stable manifold of the inertial manifold  $\left\{Im(P_{\Sigma^*}(t)): t\in\R\right\}$.
}
\end{theorem}

\proof
Given $\kappa>0$ consider the complete metric space
\begin{equation}\label{FPTheta}
\begin{split} 
\mathcal{LB}_\Theta(\kappa)=\Big\{ &\Theta\in C(\mathbb{R}\times X, X) :
  \| \Theta(t,u)-\Theta(t,\tilde{u})\| \leqslant\kappa\|u-\tilde{u}\|,\ \Theta(t,0)=0,  \\
&\Theta(t,u)=\Theta(t,(I-Q(t))u)\in Im(Q(t)),\ \forall (t,u),(t,\tilde{u})\in \R\times X\Big\}.
\end{split}
\end{equation} 
with the metric $|\!|\!|\Theta-\tilde{\Theta}|\!|\!|={\displaystyle\sup_{t\in \R}\sup_{u\neq 0} \frac{\|\Theta(t,u)- \tilde{\Theta}(t,u)\|}{\|u\|}}$.

Next we outline the heuristic procedure that will establish the way of proving that the invariant manifold is given as a graph of a map in $\mathcal{LB}_\Theta(\kappa)$.
We are looking for $\Theta \in \mathcal{LB}_\Theta(\kappa)$ with the property that, if $(\tau,\eta)\in \mathbb{R}\times X$, then a solution $u$ of 
\eqref{NLeq:glin}, with initial data $u(\tau)=\Theta(\tau,(I-Q(\tau))\eta)+(I-Q(\tau))\eta\in X$, can be decomposed as $u(t)=q(t)+p(t)$, where $q(t)=\Theta(t,p(t))$ for all $t\geqslant \tau$. Thus, $q$ and $p$ must satisfy, for $ t\geqslant \tau$,
\begin{subequations}\label{def:pq2}
\begin{align}
q(t)&=L(t,\tau)Q(\tau)\eta +\int_\tau^t L(t,s)Q(s)f(s,p(s)+\Theta (s,p(s)))ds, \\
p(t)&=L(t,\tau)(I-Q(\tau))\eta +\int_{\tau}^{t} L(t,s)(I-Q(s))f(s,p(s)+\Theta (s,p(s)))ds.
\end{align}
\end{subequations}
It follows that
\begin{equation*}
\begin{split}
\|p(t)\|&\leqslant \|L(t,\tau)(I-Q(\tau))\eta\| +\int_\tau^t \|L(t,s)(I-Q(s))f(s,p(s)+\Theta(s,p(s)))\|ds\\
&\leqslant Me^{-\gamma(t-\tau)}\|\eta\| +\int_\tau^t M \ell e^{-\gamma(t-s)}(1+\kappa)\|p(s)\|ds.
\end{split}
\end{equation*}

Using Grownwall's inequality,
\begin{equation*}
\begin{split}
\|p(t)\| \leqslant Me^{-(\gamma-M \ell (1+\kappa))(t-\tau)}\|\eta\| .
\end{split}
\end{equation*}

From this and from the fact that $q(t)=\Theta(p(t))$, we conclude that 
\begin{equation*}
\begin{split}
\|L(\tau,t)Q(t)q(t)\|= \| L(\tau,t)\Theta(p(t))\|
\leqslant 
%Me^{\rho(t-\tau)}\kappa M e^{-(\gamma-M \ell (1+\kappa))(t-\tau)}\|\eta\|= 
\kappa M^2e^{-(\gamma-\rho-M \ell (1+\kappa))(t-\tau)}\|\eta\|.
\end{split}
\end{equation*}

Applying $L(\tau,t)Q(t)$ to \eqref{def:pq2}, using that  $\Theta(\tau,(I-Q(\tau))\eta)=Q(\tau)\eta$ and making $t\to \infty$ we have 
$$
0=\Theta(\tau,(I-Q(\tau))\eta) +\int_\tau^\infty L(\tau,s)Q(s)f(s,p(s)+\Theta (s,p(s)))ds. 
$$

Inspired by this we define the operator $\tilde{G}: \mathcal{LB}_\Theta(\kappa)\to \mathcal{LB}_\Theta(\kappa)$ by
\begin{equation}\label{Gtilde}
\tilde{G}(\Theta)(\tau,\eta ) =-\int_\tau^\infty L(\tau,s)Q(s)f(s,p(s)+\Theta(s,y (s))) ds,\  (\tau,\eta)\in \R\times X.
\end{equation}

%We will prove that $\tilde{G}(\cdot)$ is a well defined map from $ \mathcal{LB}_\Theta(\kappa)$  into itself and it is a contraction.
%
%
%
%Let $\Theta \in \mathcal{LB}_\Theta(\kappa)$ and for $t\geqslant\tau$ and $\eta\in X$, consider
%\begin{equation}
%p(t):=L(t,\tau)(I-Q(\tau))\eta +\int_\tau^t L(t,s)(I-Q(s))f(s,p(s)+\Theta (s,p(s)))ds,
%\end{equation}
%and define the following operator
%\begin{equation}\label{Gtilde}
%\tilde{G}(\Theta)(\tau,\eta ) :=-\int_\tau^\infty L(\tau,s)Q(s)f(s,p(s)+\Theta(s,p (s))) ds.
%\end{equation}

The fact that $\tilde{G}$ is a well-defined contraction is similar to Theorem \ref{NLTum}, and we refrain from giving a proof. 
Hence $\tilde{G}$ admits a unique fixed point $\Theta^*\in \mathcal{LB}_\Theta(\kappa)$ satisfying the desired properties.

{\color{black}

We now embark in the proof of \eqref{attraction_stable}. For any $(\tau,\eta)\in X$ and $t\geqslant \tau$.
\begin{equation*}
p(t)=L(t,\tau)(I-Q(\tau))\eta +\int_{\tau}^t L(t,s)(I-Q(s))f(s,q(s)+p(s))ds
\end{equation*}
and thus we wish to bound the variable $\eta  (t) := T(t,\tau)u - P_{\Theta^*}(t)T(t,\tau)u$ for any $u\in X$ and $t\geqslant \tau$. 
Note that $\eta(t)=q(t)-\Theta^*(t,p(t))$ due to the definitions in \eqref{def:pq2}.

Define $p^*(s,t)$, for $s\geqslant t$, as
\begin{equation}
    p^*(s,t) :=L(s,t)p(t)+ \int_{t}^s L(s,r)(I-Q(r))f(r,\Theta^* (r,p^*(r,t))+p^*(r,t))dr.
\end{equation}
Since $f,\Theta^*$ are Lipschitz with respective constants $\ell,\kappa>0$, we obtain
\begin{equation}
\begin{split}
\|p^*(s,t)- p(s)\|&\leqslant M \int^s_t \!e^{-\gamma(s-r)}\|f(r,\Theta^* (r,p^*(r,t))+p^*(r,t))\!-\!f(r,q(r)+p(r))\| dr \\
%&\leqslant  M\ell \int^s_t e^{-\gamma(s-r)} ( \|q(r)-\Theta^* (r,p^* (r ,t)) \| +
%\| p^* (r ,t) - p(r)\|)dr, \\
&\leqslant M\ell \!\int^s_t\!  e^{-\gamma(s-r)}( \| q(r)-\Theta^* (r,p(r)) \| \!+\!
(1\!+\!\kappa )\| p^* (r ,t) - p(r)\|)dr,\\
\end{split}
\end{equation}
and, by Gr\"onwall's Lemma,
\begin{equation}
\begin{split}
\|p^*(s,t)- p(s)\|&\leqslant M\ell \int^s_t e^{-(\gamma-M\ell (1+\kappa))(s-r)}\| \eta  (r)\| dr.\label{eq-3-14-1}
\end{split}
\end{equation}
}

{\color{black}

Also, for $s\geqslant t \geqslant \tau$, we obtain 
\begin{equation*}
\begin{split}
&\|p^* (s,\tau) -p^* (s,t)\|
\leqslant \|L(s,t)(I-Q(t))[p^* (t,\tau) - p(t)]\| \\
& +\| {\ds\int_{t}^s} L(s,r)(I-Q(r))[f(r,\Theta^* (r,p^*(r,\tau))+p^*(r,\tau))-f(r,\Theta^* (r,p^*(r,t))+p^*(r,t))]\|dr\\
&  \leqslant   M^2\ell e^{-\gamma(s-t)} \int^t_\tau e^{-(\gamma-M\ell (1+\kappa))(t-r)}\| \eta  (r)\| dr  + M\ell (1+\kappa ) {\ds\int^s_{t}}  e^{-\gamma(s-r)}
\|p^* (r,\tau) - p^* (r,t)\| dr,
\end{split}
\end{equation*}
and again by Gr\"onwall's Lemma,
\begin{equation}\label{eq-3-15-1}
\|p^* (s,\tau) -p^* (s,t)\|\leqslant M^2\ell \,  \int_\tau^t e^{-(\gamma -M\ell (1+\kappa) )(s-r)}\| \eta  (r)\| dr.
\end{equation}
Now, we use these inequalities to estimate $\|\eta(\tau)\|$. Note that
\begin{equation*}
\begin{aligned}
&\eta (\tau) -L(\tau,t)Q(t)\eta (t) = q (\tau) -L(\tau,t)q(t)\!-\! \Theta^*(\tau,p(\tau)) \!+\!L(\tau,t)\Theta^* (t,p(t))\\
%&=\!\!{\ds\int_{t }^\tau} L(\tau,s)Q(s)f(s,q(s)+p(s))ds+\int_\tau^\infty L(\tau,s)Q(s)f(s,\Theta^* (s, p^*(s,\tau))+p^* (s,\tau))ds \\
%&-\int_t^\infty L(\tau,s)Q(s)f(s,\Theta^*(s,p^*(s,t))+p^* (s,t))ds \\
&=\int_t^\tau L(\tau,s)Q(s)[f(s,q(s)+p(s))-f(s,\Theta^*(s,p^*(s,\tau))+p^* (s,\tau))]ds\\
&+\int_t^\infty L(\tau,s)Q(s)[f(s, \Theta^* (s, p^*\!(s,\tau))+p^*\!(s,\tau))\!-\!f(s,\Theta^* (s, p^*\!(s,t))+p^* \!(s,t)) ]ds.
\end{aligned}
\end{equation*}

    Thus, using \eqref{eq-3-14-1} and \eqref{eq-3-15-1},
we obtain
\begin{equation*}
\begin{split}
&\|\eta (\tau) -L(\tau,t)Q(t)\eta (t)\|\leqslant M\ell {\ds\int_\tau^{t}}  e^{-\rho(\tau-s)} \left( \|q(s) - \Theta^* (s,p^* (s,\tau))\| + \|p(s) - p^* (s,\tau)\| \right)ds \\
& + M\ell (1+\kappa )\int_t^{\infty}e^{-\rho(\tau-s)}  \|p^* (s,t) - p^* (s,\tau)\| ds \\
&\leqslant  M\ell \!{\ds\int_\tau^t}\! e^{-\rho(\tau-s)}\|\eta  (s)\|ds +M^2\ell^2(1\!+\!\kappa)\!{\ds\int_\tau^t} \! e^{-(\gamma-\rho-M\ell (1+\kappa) )(s-\tau)} \int^s_\tau e^{-(\gamma-M\ell (1+\kappa))(\tau-r)}\| \eta  (r)\| drds \\
&+M^3\ell^2 \,  (1+\kappa ){\ds \int^{\infty}_{t}}  
e^{-(\gamma-\rho-M\ell (1+\kappa) )(s-\tau)} \int_\tau^{t}e^{-(\gamma -M\ell (1+\kappa) )(\tau-r)}\| \eta  (r)\| dr ds.
\end{split}
\end{equation*}
Hence
\begin{equation*}
\begin{split}
\|\eta (\tau) -L(\tau,t)Q(t)\eta (t)\|&\leqslant  M\ell \!{\ds\int_\tau^t}\! e^{-\rho(\tau-s)}\|\eta  (s)\|ds   \\
&+ \frac{M^2\ell^2(1+\kappa)(1+M)}{\gamma-\rho-M\ell(1+\kappa)}{\ds\int_\tau^t} 
e^{-(\gamma-\rho-M\ell(1+\kappa)) (\tau-r)}e^{-\rho (\tau-r)} \| \eta  (r)\| \, dr\\
\end{split}
\end{equation*}
and we have that
\begin{equation*}
 \|\eta (\tau) -L(\tau,t)Q(t)\eta (t)\| \leqslant   \left[M\ell+\frac{M^2\ell^2(1+\kappa)(1+M)}{\gamma-\rho-M\ell(1+\kappa)}\right] {\ds\int_\tau^t} 
 e^{-\rho(\tau-r)} \|\eta(r)\|dr.
\end{equation*}
Thus, 
\begin{equation*}
\|\eta (\tau)\| \leqslant M e^{-\rho(\tau-t)} \|\eta  (t)\|+  \left[M\ell +\frac{M^2\ell^2(1+\kappa)(1+M)}{\gamma-\rho-M\ell (1+\kappa)}\right]{\ds\int_\tau^t} 
e^{-\rho(\tau-r)} \|\eta(r)\|dr.
\end{equation*}
%and
%\begin{equation*}
%e^{\rho t}\|\eta (\tau)\| \leqslant  M e^{\rho t} \|\eta  (t)\|+ \left[M\ell +\frac{M^2\ell^2(1+\kappa)(1+M)}{\gamma-\rho-M\ell (1+\kappa)}\right] {\ds\int_\tau^t} 
% e^{\rho r} \|\eta(r)\|dr.
%\end{equation*}
By Gr\"onwall's Lemma, we obtain the bound in \eqref{attraction_stable}.

}

\cqd

%{\color{blue}
%\begin{remark}
%Note that, for $\gamma>0$, Theorem \ref{NLTsm}  implies that for each global solution in the inertial manifold we have a complementary direction to the inertial manifold which is contained in its stable manifold.
%\end{remark}
%}

Note that Theorem \ref{NLTsm} does not follow from the exponentially dominated property \eqref{expattSigmaNL} in Theorem \ref{NLTum}. Indeed, the proof follows closely the first part of the proof of the existence of the invariant manifold, and estimates \eqref{DecayStM} and \eqref{expattSigmaNL} are closely related to the estimates \eqref{invSigmaNL} and \eqref{attraction_stable}. 
\begin{remark} 
Theorem \ref{NLTsm} (when $\delta>0$) provides a complementary direction of the inertial manifold that is contained in the stable manifold of $u\equiv 0$, which consisting of special solutions that decay to the trivial equilibrium $u\equiv 0$ with exponential decay rate $-(\gamma- M\ell(1+\kappa))<0$ according to \eqref{DecayStM}.
In fact, for any global bounded solution within the inertial manifold (e.g. hyperbolic equilibrium or a normally hyperbolic periodic orbit), one can shift such solution to zero according to Remark \ref{Rmk:shift}, and hence obtain a complementary direction to the inertial manifold that its contained in the stable manifold of the aforementioned object. Moreover, in the case $\hat{\delta}>0$, solutions outside this manifold grow apart  exponentially from it, according to \eqref{attraction_stable}.
%Hence, any global bounded solutions within the inertial manifold have a respective stable manifold with solutions that converge to rate $\gamma-\ell M(1+\kappa)>0$, 
%This is in contrast with the decay rate $\delta>0$ of any point in phase-space (which may not be in the stable manifold of any object within the inertial manifold) towards the inertial manifold in \eqref{deltaNL}. 
\end{remark}

\begin{remark} \label{rmk:translateto0}
The above results in Theorems \ref{NLTum} and \ref{NLTsm} are proved for general operators $A(t)$ that generate a linear evolution process and for nonlinearities $f(t,.):X\to X$, and consequently the are well suited to be applied for hyperbolic PDEs. However, these results can also be replicated for parabolic PDEs, in the case that $-A$ is a sectorial operator, and more importantly, if $f:\mathbb{R}\times X^\alpha\to X$, where $X^\alpha$ denotes the fractional power space with $\alpha \in (0,1)$. 
Since these results do not bring new insights in the proof, we refrain from giving a detailed proof. However, we mention some differences.
As a preparation, note that the analytic semigroup $e^{At}$ satisfies $\|e^{At}Qu\|_{X^\alpha}\leqslant M e^{-\rho t}\|u\|_{X^\alpha}$ and $\|e^{At}Qu\|_{X^\alpha}\leqslant N e^{-\rho t}\|u\|_X$ respectively for $t\leqslant 0$ and $t<0$.
Therefore, following an analogous proof of Theorem \ref{NLTum}, the equations \eqref{nuNL} actually become
\begin{subequations}\label{nuNLpar}
\begin{align}
\frac{\ell M^2 (1+\kappa )\Gamma(1-\alpha)}{[\gamma-\rho-2\ell N (1+\kappa)]^{1-\alpha}}&\leqslant \kappa, \label{nuNLpar1}\\ 
\frac{2\ell M^2 \Gamma(1-\alpha)}{[\gamma-\rho-2\ell N (1+\kappa)]^{1-\alpha}} & < 1,\label{nuNLpar2}
\end{align}
\end{subequations}
where $\Gamma(\cdot)$ is the gamma function.
There are different ways to achieve \eqref{nuNLpar} by appropriately choosing $\kappa>0$. We present a similar choice as in Theorem \ref{NLTum} which fits our purposes in the upcoming sections.
%First, note that in order for the above fractions to be well-defined, the denominator needs to be positive, which yields an upper bound on $\kappa>0$,
%\begin{equation}
%    \kappa<\kappa_+:=\frac{\gamma-\rho}{2M\ell}-1.    
%\end{equation}
Equation \eqref{nuNLpar2} yields an upper bound on the Lipschitz constant according to 
\begin{equation}
    \kappa<\kappa_*:=\frac{\gamma-\rho-[2M^2\ell \Gamma(1-\alpha)]^{\frac{1}{1-\alpha}}}{2N\ell}-1.
\end{equation}
Note that $\kappa_*<(\gamma-\rho)/(2N\ell)-1$ and thus the denominators in \eqref{nuNLpar} are positive.
Moreover, we rewrite \eqref{nuNLpar1} as $\ell M^2 (1+\kappa )\Gamma(1-\alpha) \leqslant \kappa [\gamma-\rho-2\ell N (1+\kappa)]^{1-\alpha}$. Note that an upper bound on $\kappa\leqslant \kappa_+$, for some $\kappa_+\in \mathbb{R}_+$, is able to give a lower bound on the right-hand side of this inequality, $ \kappa [\gamma-\rho-2\ell N (1+\kappa_+)]^{1-\alpha} \leqslant \kappa [\gamma-\rho-2\ell N (1+\kappa)]^{1-\alpha}$, and if such lower bound is also an upper bound for the left-hand side of the desired inequality, $ \ell M^2 (1+\kappa )\Gamma(1-\alpha)\leqslant \kappa [\gamma-\rho-2\ell N (1+\kappa_+)]^{1-\alpha} $, then we obtain the desired result. %Note this procedure gets rid of the fractional power of $\kappa>0$.
However, since the upper bound $\kappa<\kappa_*$ does not yield a meaningful result, we consider the following upper bound,
\begin{equation}
    \kappa\leqslant\kappa_{+}:=\frac{1}{2}\frac{\gamma-\rho-[2M^2\ell \Gamma(1-\alpha)]^{\frac{1}{1-\alpha}}}{2N\ell}-1,
\end{equation}
and thereby the desired inequality yields a lower bound on $\kappa>0$,
\begin{equation}
    \kappa\geqslant\kappa_-:=\frac{2^{1-\alpha}M^2\ell\Gamma(1-\alpha)}{\left[\gamma-\rho+(2M^2\ell \Gamma(1-\alpha))^{\frac{1}{1-\alpha}}\right]^{1-\alpha}}.
\end{equation}
Thus for any $\kappa \in [\kappa_-,\kappa_+)$, we obtain a well-defined contraction and prove Theorem \ref{NLTum} in the parabolic case.
%Note that $\kappa_{-}\to 0$ and $\kappa_{+}\to \infty$ if either $(\gamma-\rho)\to\infty$ or $\ell\to 0$, and thereby the interval $ [\kappa_-,\kappa_+)$ is not empty for sufficiently large $(\gamma-\rho)\to\infty$ or sufficiently small $\ell\to 0$.
%
%Thus one can then modify the proof, in case that the hypothesis \eqref{NLeq:inequalities} is replaced by the following hypothesis
%\begin{equation}
%\frac{\gamma-\rho}{\ell [\gamma-\rho]^{1-\alpha}}>c(M),
%\end{equation}
%for some constant $c(M)\mathbb{R}_+$.
{\color{black}Similarly, modifying the proof accordingly, the estimate given in \eqref{expattSigmaNL} will be satisfied with exponent}
\begin{equation}
     \delta_{par}:=\gamma-\rho-2N\ell (1+\kappa)- \left[ 2 \Gamma(1-\alpha) \ell M \left(1+\frac{\ell (1+2M)N (1+\kappa)}{\gamma-\rho- 2\ell N (1+\kappa)}\right) \right]^\frac{1}{1-\alpha}.
\end{equation}
{\color{black} In applications the interplay between the size of the gap $\gamma-\rho$ and the constant $N$ will play a major role.} See \cite{Mik} for sharp conditions that guarantee the existence of inertial manifolds. 
\end{remark}

\subsection{Local Inertial Manifolds}\label{sec-local-inertial-manifolds}

\par The aim of this subsection is to obtain a local version of the inertial manifold in Theorem \ref{NLTum} in the particular case that the diffusion $A(t)\equiv A$ and the non-linearity $f(t,\cdot)$ is Lipschitz in a neighborhood of an attracting set, which guarantees the existence of a local inertial manifold only in this neighborhood. For this purpose, we will need a stability property of the attracting set, which is natural when dealing with autonomous dynamical systems. Thus we will use the notion of skew product semiflow associated with our problem.

\par First, we recall some basic notions of skew product semiflows. Let $\mathcal{C}(\mathbb{R}\times X,X)$ be the space of all continuous functions 
$g:\mathbb{R}\times X\to X$. Let $\Sigma$ be a subset of  $\mathcal{C}(\mathbb{R}\times X,X)$ which is translation invariant, that is, for any  $g\in \Sigma$ we have that $\theta_tg\in \Sigma$ for all $t\in \R$, where $\theta_t (s,x)=g(t+s,x)$, for all $(s,x)\in \R\times X$.
Assume that  $\Sigma$ is endowed with a metric which makes it a compact metric space and that $f\in \Sigma$.

%We define the {\bfseries  continuous time-shift operator} $\theta_t$ of $g\in \mathcal{C}(\mathbb{R}\times X,X)$ as $\theta_t g(\cdot,\cdot):=g(t+\cdot,\cdot)$,
%%\begin{equation*}
%%\theta_t g(\cdot,\cdot):=g(t+\cdot,\cdot),  
%%\end{equation*}
%for each $t\in \mathbb{R}$.
%
%Consider the vector field of equation \eqref{NLeq:glin} given by $h(t,u):=Au+f(t,u)$, for $(t,u)\in \mathbb{R}\times X$. 
Define the set of all $t$-translations of $f$ by $\mathcal{B}_f:=\{\theta_tf: t\in \mathbb{R}\}$, and denote by $\mathcal{B}$ the closure of $\mathcal{B}_f$ with the metric of $\Sigma$, 
%\begin{equation*}
%\mathcal{B} = \hbox{ the closure of } \mathcal{B}_h \hbox{ in } \mathcal{C}(\mathbb{R}\times X,X),
%\end{equation*}
which is known as the {\bfseries Hull of $f$}. 
The hull $\mathcal{B}$ is compact.
%see %\cite[Theorem 14]{Sell67} {\color{black} (neste caso $X=\R^n$ \'{e} localmente compacto)}  and  
%\cite[Chapter V]{Chepyzhov-Vishik-02}. %{\color{black} este caso \'{e} mais abstrato e n\~{a}o especifica a topologia)}.

\begin{comment}

{\color{blue} POSSIBLE SOLUTION!!!}

{\color{black} 
1) Vamos ter que considerar uma que o espa\c{c}o de s\'{\i}mbolos seja um espa\c{c}o m\'{e}trico compacto. 

2) Imagino que isto seja solucionado da seguinte maneira. 

2.1) Assuma que existem espacos $Y,Z$ tais que $Z\subset\subset X\subset\subset Y$ que $f: \R\times Y\to Z$ seja uniformemente cont\'{\i}nua em limitados de $\R\times Y$ e que leve limitados em limitados. 

2.2) Considere a topologia da converg\^{e}ncia uniforme nas partes compactas de $C(\R\times Y,X)$ como espa\c{c}o de s\'{\i}mbolos e depois tome o Hull com respeito a esta topologia.

2.3) Em seguida considere o subsepa\c{c}o $X$ de $Y$ que pode ser escrito como uni\~{a}o cont\'{a}vel de subconjuntos compactos de $Y$. 

2.4) As restri\c{c}\~oes dos elementos de $C(\R\times Y,X)$ a $X$, com esta topologia nos dar\'{a} um Hull que ser\'{a} compacto como conseq\"{u}\^{e}ncia do Teorema de Ascoli-Arzel\'{a}.}

\end{comment}

%{\color{blue} ALTERNATIVE SOLUTION!!!}

%{\color{black}Seguir o contexto mais abstrato que est\'{a} descrito em \cite[Chapter V]{Chepyzhov-Vishik-02}.}

%$f\in \mathcal{C}(\mathbb{R}\times X,X)$ is bounded and uniformly continuous in $t$, then $\mathcal{B}$ is a compact metric space, with the induced topology of $\mathcal{C}(\mathbb{R}\times X,X)$

We assume that for each $b\in \mathcal{B}$ the problem
\begin{equation}
\dot{u}=Au+b(t,u), \ \ t\geqslant 0,
\end{equation} 
has a unique solution given by $\varphi(t,b)u_0\in X$ for each $t\geqslant 0$, with initial condition $u(0)=u_0\in X$ and 
$\varphi:\mathbb{R}\times \mathcal{B}\times X\to X$ defines a {\bfseries co-cycle}, i.e., 
$\varphi(0,b)=I_X$ and $\varphi(t+\tau,b)=\varphi(t,\theta_\tau b)\varphi(\tau ,b)$, for $t,\tau \geqslant 0$. 
In particular, since $\{T(t,\tau): t\geqslant \tau \}$ is the evolution process associated to \eqref{NLeq:glin}, we have that $T(t,\tau)=\varphi(t-\tau,\theta_\tau f)$ for any $t\geqslant \tau$.
%\begin{equation*}
%    T(t,\tau)=\varphi(t-\tau,\theta_\tau h), \ \ t\geqslant \tau.
%\end{equation*}

\par The continuous time-shift operator $\theta_t$ and the co-cycle $\varphi$ define the \textbf{skew-product semigroup} $\{\Pi(t): t\geqslant 0\}$ on the phase-space $\mathbb{X}:=X\times \mathcal{B}$, given by
\begin{equation}
\Pi(t)(u,b):=(\varphi(t,b)u,\theta_t b), 
\end{equation}
where $(u,b)\in X\times \mathcal{B}$, $t\geqslant 0$ and the space $\mathbb{X}$ endowed with metric $d_\mathbb{X}(\cdot,\cdot)=\max\{\|\cdot\|_X, d_\mathcal{B}(\cdot)\}$.

A compact set $\mathcal{A}\subseteq X$ is called the \textbf{uniform attractor} for the co-cycle $\varphi$, if for any bounded set $U\subseteq X$, we have that $\lim_{t\to +\infty} \sup_{b\in \mathcal{B}} dist_H(\varphi(t,b)U, \mathcal{A})=0$. Similarly, a compact set $\mathbb{A}\subseteq \mathbb{X}$ is called the \textbf{global attractor} for the skew-product semi-flow $\{\Pi(t): t\geqslant 0\}$, if it is invariant and attracts all bounded sets of $\mathbb{X}$.
The global attractor for the skew-product semi-flow is related to the uniform attractor of the co-cycle, when both exist, according to $\mathcal{A}:=P_X\mathbb{A}$, where $P_X:\mathbb{X}\to X$ is the projection onto $X$. 
	
%For any subset $B\subseteq X$ and $\epsilon>0$ we define $O_\epsilon(B):=\{x\in X: d(x,B)<\epsilon\}$ as the \textbf{$\epsilon$-neighborhood of $B$}, where $d(x,B):=\inf_{b\in \mathcal{B}} \|x-y\|_X$ for $y\in B$.

%\begin{definition}
%	Let $(\varphi,\theta)$ be the co-cycle that defines the skew product semi-flow $\{\Pi(t): t\geqslant 0\}$ and $\mathcal{A}$ be a compact subset of $X$. We say that	$\mathcal{A}$ is a \textbf{uniform attractor} for $\varphi$, if  for any bounded set $B\subseteq X$ we have that
%	\begin{equation*}
%	\lim_{t\to +\infty} \sup_{b\in \mathcal{B}} dist_H(\varphi(t,b)B, \mathcal{A})=0.
%	\end{equation*}
%\end{definition}

%\begin{remark}
%	For any skew-product semi-flow $\{\Pi(t): t\geqslant 0\}$ with a global attractor 	$\mathbb{A}$, the set $\mathcal{A}:=P_X\mathbb{A}$ is the uniform attractor for $\varphi$, where $P_X:\mathbb{X}\to X$ is the projection onto $X$. 
%\end{remark}

\par We now prove that the uniform attractor satisfies certain stability property. 

\begin{lemma}\label{lemma-uniform-attractor-stable}

{\color{black} 
Assume that $ \mathcal{B}$ is as above}. Let $\{\Pi(t): t\geqslant 0\}$ be the associated skew-product semi-flow in $\mathbb{X}$ and assume that it has a global attractor $\mathbb{A}$. Then for a given $\epsilon>0$, there exists 
$\epsilon_*\in (0,\epsilon)$ such that
\begin{equation}\label{eq-uniform-attractor-stable}
\bigcup_{t\geqslant 0} \bigcup_{ b\in \mathcal{B} }
\varphi(t,b) O_{\epsilon_*}(\mathcal{A}) \subset O_\epsilon(\mathcal{A}),
\end{equation}
where $O_\epsilon(\mathcal{A})\subseteq X$ is an $\epsilon$-neighborhood of the uniform attractor $\mathcal{A}$.
%where $\mathcal{I}$ is the global attractor of $\{\theta_t:t\in \mathbb{R}\}$.
\end{lemma}
\proof
    Since $\mathbb{A}$ is the global attractor for $\{\Pi(t): t\geqslant 0\}$, then  given $\epsilon>0$, there exists 
 	$\epsilon_*\in (0,\epsilon)$ such that
 	\begin{equation}\label{eq-stability-global-attractor}
 	\bigcup_{t\geqslant 0} \Pi(t) O_{\epsilon_*}(\mathbb{A}) \subset O_\epsilon(\mathbb{A}).
 	\end{equation}
    where $O_\epsilon(\mathbb{A})\subseteq \mathbb{X}$ is an $\epsilon$-neighborhood of the global attractor $\mathbb{A}$.
	For the proof of this fact, see for instance \cite{Hale-88}.
 	Next, towards a contradiction, suppose that \eqref{eq-uniform-attractor-stable} does not hold for these $\epsilon,\epsilon_*$. Hence there exists $t_0\geqslant 0$, $b\in \mathcal{B}$, and 
	$u\in O_{\epsilon_*}(\mathcal{A})$ such that
	$\varphi(t,b)u\notin O_\epsilon(\mathcal{A})$. 
	Since $\mathbb{A}=\mathcal{A}\times \mathcal{B}$, we have 
	$(u,b)\in O_{\epsilon_*}(\mathbb{A})$, which from \eqref{eq-stability-global-attractor} we have that
	$\Pi(t)(u,b)\in O_\epsilon(\mathbb{A})$, for all $t\geqslant 0$. In particular, 	$d(\varphi(t,b)u,\mathcal{A})<\epsilon$, for all $t\geqslant 0$, which contradicts the assumption that 
	$\varphi(t,p)u\notin O_\epsilon(\mathcal{A})$, and the proof is complete.
\cqd

\medskip
\begin{comment}
that is, there exists 
		$\Sigma:\mathbb{R}\times O_\epsilon(\mathcal{A})\to X$ which $u\mapsto\Sigma(t,u)$ is uniformly Lipschitz uniformly in $t$ such that 
\end{comment}

We are ready to define and prove the existence of local inertial manifolds for equation \eqref{NLeq:glin}. Note that instead of requiring that $f$ is globally Lipschitz on $X$, as in Theorem \ref{NLTum}, we now suppose that it is Lipschitz on bounded sets.
\begin{definition}\label{loc:IM}
	Let $\epsilon>\epsilon_*>0$ be such that 
	$T(t,\tau) O_{\epsilon_*}(\mathcal{A}) \subset O_\epsilon(\mathcal{A})$, for every $t\geqslant \tau$. 
	A family of subsets $\{\mathcal{M}_{\text{loc}}(t): t\in \mathbb{R}\}\subseteq X$ is said to be a \textbf{local inertial manifold} if
	\begin{enumerate}
		\item $\mathcal{M}_{\text{loc}}(t)$ is a Lipschitz manifold for each $t\in \mathbb{R}$.
		\item $T(t,\tau)[ O_{\epsilon_*}(\mathcal{A})\cap\mathcal{M}_{\text{loc}}(\tau)]\subset O_{\epsilon}(\mathcal{A})\cap \mathcal{M}_{\text{loc}}(t)$, $t\geqslant \tau$.
		\item $\{\mathcal{M}_{\text{loc}}(t): t\in \mathbb{R}\}$ is exponentially attracting in the sense of Definition \ref{def-inertial-manifold}, part 3.
	\end{enumerate}
\end{definition}

\begin{theorem}\label{th-existence-local-inertial-manifolds-skew-products}
Consider the equation \eqref{NLeq:glin} in which $f:\mathbb{R}\times X\to X$ is Lipschitz in bounded sets uniformly in $t\in\mathbb{R}$, i.e., $\|f(t,u)-f(t,\tilde{u})\|_{X}\leqslant \ell \|u-\tilde{u}\|_X$, for any $(t,u),(t,\tilde{u})$ within bounded subsets of $X$. 
Suppose that the linear evolution process $\{L(t,\tau): t\geqslant \tau\}$ has exponential splitting, with constant $M\geqslant 1$, exponents $\gamma>\rho$ with $\gamma>0$, and a family of projections $\{Q(t): \ t\in \R\}$.
Assume also that the associated skew-product semigroup $\{\Pi(t): t\geqslant 0\}$ has a global attractor $\mathbb{A}$ and uniform attractor $\mathcal{A}$.
    % paragrafo aqui
    
If $\gamma$ and $\gamma-\rho$ are sufficiently large, % for some Lipschitz constant $\ell>0$ in the bounded neighborhood $ O_{\epsilon}(\mathcal{A})$, % of the global attractor $\mathcal{A}$, 
then there exists a function $\Sigma_{\text{loc}}:\mathbb{R}\times O_{\epsilon}(\mathcal{A})\to X$ such that  $\Sigma_{\text{loc}}(t,u)\in Im(I-Q(t))$, which is uniformly Lipschitz with constant $\kappa>0$, %=\kappa(\gamma,\rho,\ell,M)>0$ , 
i.e., $\|\Sigma_{\text{loc}}(t, u)-\Sigma_{\text{loc}}(t, \tilde{u})\|_{X}\leqslant \kappa \|u-\tilde{u}\|_{X}$ for all $(t,u),(t,\tilde{u})\in  O_{\epsilon}(\mathcal{A})$. 

Furthermore, the local graph of $\Sigma_{\text{loc}}(t,.)$, for each $t\in \mathbb{R}$, given by 
%\begin{equation*}
%\mathcal{M}(t)=\{Q(t)u+\Sigma(t,Q(t)u): u\in O_{\epsilon_0}(\mathcal{A}) \}, \ \ t\in \mathbb{R}.
%\end{equation*}
\begin{equation}\label{loc:graph}
    \mathcal{M}_{\text{loc}}(t):=\lbrace u\in X: u=q+\Sigma_{\text{loc}}(t,q), q\in Q(t)O_{\epsilon}(\mathcal{A}) \rbrace ,
\end{equation}
%\begin{equation*}
%\mathcal{M}(t)=\{Q(t)u+\Sigma(t,Q(t)u): u\in O_{\epsilon_0}(\mathcal{A}) \}, \ \ %t\in \mathbb{R}.
%\end{equation*}
yields a local inertial manifold $\{\mathcal{M}_{\text{loc}}(t): t\in \mathbb{R}\}$ for 
$\{T(t,\tau): t\geqslant \tau\}$.
\begin{comment}
\textcolor{red}{PL: local?} inertial manifold for $\{T(t,\tau): t\geqslant \tau\}$.  \textcolor{red}{PL: what is $T(t,s)$ here? this is not true. If $f$ is discontinuous there is NO evolution process at all.} \textcolor{lightgreen}{This is not true, Caratheodory Differential equations for instance, the function $f$ is not continuous. $f$ being monotone and bounded is enough to prove global well-posedness. Some papers of Longo , Obaya Novo, deal with these pathological situations. However I agreed in making the statatement more clear. Put just Lipschitiz in bounded sets and in a Remark Explain that the only needed assumption is Lipschitz in a neighborhood of the uniform attractor.}

	$\widehat{\mathcal{M}}$ pullback attracts exponentially every bounded 
	subset of $X$. More precisely, for each bounded set $B$ and $t\in \mathbb{R}$, there exists 
	$s_0=s_0(B,t)\leqslant t$ such that
	$T(t,s)B\subset O_{\epsilon_0}(\mathcal{A})$, $s\leqslant s_0$ and 
	$$
	\| T(t,s)b - (Q(t)T(t,s)b+\Sigma (t,T(t,s)b))\|_{X} \leqslant C(B)e^{-\delta (t-s)}, \quad s\leqslant s_0, \  b\in B.
	$$
	In particular,  Note that, $C(B)$ does not depend of the initial time $s$.
	\item $\widehat{\mathcal{M}}$ forward attracts exponentially every bounded subset of $X$.
	\item The family 
	$\{\mathcal{M}(t)\cap O_{\epsilon_0'}(\mathcal{A}): t\in \mathbb{R}\}$ satisfies
	\begin{equation*}
	T(t,\tau)[\mathcal{M}(\tau)\cap O_{\epsilon_0'}(\mathcal{A})]
	\subset \mathcal{M}(t).
	\end{equation*}
\end{comment}
\end{theorem}

\proof
Denote by $\ell>0$ the Lipschitz constant of $f$ in the bounded neighborhood $ O_{\epsilon}(\mathcal{A})$. % of the global attractor $\mathcal{A}$, 
Let $\tilde{f}:\mathbb{R}\times X\to X$ be a Lipschitz extension of $f$ out
of $O_{\epsilon}(\mathcal{A})$, which coincides with $f$ on $O_{\epsilon}(\mathcal{A})$ and has global Lipschitz constant $\ell>0$. 
Conditions of Theorem \ref{NLTum} are satisfied for $\tilde{f}$ and $\gamma-\rho$ sufficiently large, and thus there exists an inertial manifold $\widetilde{\Sigma}\in \mathcal{LB}_\Sigma(\kappa)$ for the equation \eqref{NLeq:glin}, with $f$ replaced by $\tilde{f}$.
%Since $h(t,\cdot)=A(t)+f(t,\cdot)$, we have $T(t,\tau)=\varphi(t-\tau,\theta_\tau h)$,  $t\geqslant \tau$,\textcolor{red}{I would delete this sentence before this point... there is no argument here, is there?} 

We then consider the restriction of the global inertial manifold $\widetilde{\Sigma}$ to the neighborhood $O_{\epsilon}(\mathcal{A})$, i.e., define, $\Sigma_{\text{loc}}:=\widetilde{\Sigma}|_{\mathbb{R}\times O_{\epsilon}(\mathcal{A})}$ and we obtain the sets defined in \eqref{loc:graph}. 
In order to prove $\mathcal{M}_{loc}(
\cdot)$ yields a local inertial manifold, note that it is a Lipschitz graph, since $\widetilde{\Sigma}$ is Lipschitz. Moreover, due to Lemma \ref{lemma-uniform-attractor-stable}, there is a 
$\epsilon_*\in (0,\epsilon)$ such that $T(t,\tau) O_{\epsilon_*}(\mathcal{A}) \subset O_\epsilon(\mathcal{A})$ for $t\geqslant \tau$.
Together with the fact that $\widetilde{\Sigma}$ yields a global invariant set, this implies the second condition in definition \ref{loc:IM}.

Lastly, we prove that the graph of the local inertial manifold in \eqref{loc:graph} is exponentially attracting. 
%\begin{equation*}
%\mathcal{M}(t):=\{Q(t)u+\Sigma(t,Q(t)u): u\in O_{\epsilon_0}(\mathcal{A}) \}.
%\end{equation*}
Given $u\in O_{\epsilon_*}(\mathcal{A})$, we have
$T(t,\tau)u\in O_{\epsilon}(\mathcal{A})$ for $t\geqslant \tau$.
Since $\Sigma_{\text{loc}}$ coincides with $\widetilde{\Sigma}$ on the neighborhood $O_{\epsilon}(\mathcal{A})$, which is exponentially attracting according to \eqref{expattSigmaNL}, we obtain,
\begin{equation}\label{eq-general-local-inertial-manif-skewproduct}
\begin{split}
    \| T(t,\tau)u - P_{\Sigma_{\text{loc}}}(t) T(t,\tau)u \|_{X} 
    &\leqslant M\| (I-P_{\Sigma_{\text{loc}}}(\tau))u\|_{{X}} e^{-\delta (t-\tau)}\\
    &\leqslant M(1+M+\kappa)e^{-\delta (t-\tau)}\|u\|_X, \qquad \forall\, t\geqslant \tau,
%\| T(t,\tau)u_0 - (Q(t)T(t,\tau)u_0+\Sigma (t,T(t,\tau)u_0))\|_{X} \leqslant M\| (I-Q (\tau))u_0 - \Sigma (\tau,u_0)\|_{{X}} e^{-\delta (t-\tau)}.
\end{split}
\end{equation}
for any $u\in O_{\epsilon_*}(\mathcal{A})$, where $P_{\Sigma_{\text{loc}}}(\tau)u:= Q(\tau)u+\Sigma_{\text{loc}} (\tau,u)$.
Note that $\delta>0$ in \eqref{expattSigmaNL} for sufficiently large $\gamma$.

We now prove that the family $\{\mathcal{M}_{\text{loc}}(t): t\in \mathbb{R}\}$ exponentially pullback attracts any bounded subset of $O_{\epsilon_*}(\mathcal{A})$.
Since $\mathcal{A}$ is bounded and 
$\mathcal{A}=\cup_{t\in \mathbb{R}} \mathcal{A}(t)$, we have that $\mathcal{A}(t)\subset \mathcal{M}_{\text{loc}}(t)$ for each $t\in \mathbb{R}$, and consequently  $\{\mathcal{M}_{\text{loc}}(t): t\in \mathbb{R}\}$ pullback attracts %(not necessarily in exponential fashion) 
bounded sets of $X$. %, since $\mathcal{A}(t)$ are the pullback attractors.
Let us show that the rate of convergence is indeed exponential.

Since $\mathcal{A}$ is an uniform attractor, for each bounded subset $U\subseteq X$, there exists $t_*=t_*(U, \epsilon_*)\geqslant 0$ such that $T(t_*+\tau,\tau)U\subset O_{\epsilon_*}(\mathcal{A})$, for every $\tau\in \mathbb{R}$. 
Let $t\in \mathbb{R}$ be fixed and choose any $\tau\leqslant\tau_*:=t-t_*$, then from Lemma \ref{lemma-uniform-attractor-stable}, we have that $T(t,\tau)U\subset O_{\epsilon}(\mathcal{A})$, for every $\tau\leqslant \tau_*$.
 Thus, for any $u_0\in U$, the evolution $T(t_*+\tau,\tau)u_0\in O_{\epsilon_*}(\mathcal{A})$, for every $\tau\leqslant \tau_*$. Thus applying the inequality \eqref{eq-general-local-inertial-manif-skewproduct} at the point $T(t_*+\tau,\tau)u_0$, we obtain
\begin{equation}
\begin{split}
    \| T(t,\tau)u_0 - P_{\Sigma_{\text{loc}}}(t)T(t,\tau)u_0\|_{X} &\leqslant  M(1+M+\kappa) \|T(t_0+\tau,\tau)u_0\|_X  e^{\delta t_0} e^{-\delta (t-\tau)}\\
    &\leqslant C e^{-\delta (t-\tau)}, \qquad \forall\, \tau\leqslant \tau_*,
\end{split}
\end{equation}
%for all $\tau\leqslant \tau_*$ 
for any $u_0\in U$, where $C>0$ is some constant independent of $\tau$.
%\textcolor{red}{i dont undestand this ending.. how do you get above? how to go to eq below? here the RHS still depends on $\tau$ (ie, $T(t_0+\tau,\tau)u$, so, how to pass on the limit?)} \textcolor{lightgreen}{$T(t_0+\tau,\tau)u$ is in a bounded set for every $\tau \leqslant \tau_0$ .}
%
Thus,
\begin{equation}
dist_H(T(t,\tau)U,\mathcal{M}_{\text{loc}}(t))\leqslant C e^{-\delta(t-\tau)}, \qquad \forall \,\tau\leqslant \tau_*,
\end{equation}
which proves that $\{\mathcal{M}_{\text{loc}}(t): t\in \mathbb{R}\}$ is pullback exponentially attracting.
Analogous arguments yield forward exponentially attraction, and consequently exponential attraction. 
%Therefore $\{\mathcal{M}_{\text{loc}}(t): t\in \mathbb{R}\}$ is a local inertial manifold for $\{T(t,\tau): t\geqslant \tau\}$ and the proof is complete.
\cqd
%
\begin{comment}
{\color{lightgreen}
\begin{remark}
To obtain a positively invariant local inertial manifold, one can proceed as follows:

Let $\widetilde{\mathcal{M}}_{\epsilon}= \Pi_X (\widetilde{\mathcal{M}} )\cap O_{\epsilon}$, where
\begin{equation*}
    \widetilde{\mathcal{M}}=\{(\tau,Q(\tau)+\widetilde{\Sigma}(\tau,Q(\tau)u))\in \mathbb{R}\times X: u\in X \}.
\end{equation*}
Then the positively invariant local inertial manifold is 
\begin{equation*}
    \mathcal{M}=\{(\tau,u)\in \mathbb{R}\times X: (t,T(t,\tau)u)\in \widetilde{\mathcal{M}}_{\epsilon}, \ \ \forall t\geqslant \tau \},
\end{equation*}
\end{remark}
}
\textcolor{red}{PL: I do not understand the above remark...}
\end{comment}
%
\begin{comment}
\textcolor{red}{i would delete the below rmk... it follows from thm 2.1, no?}
%
\begin{remark}
    Note that, from the proof of Theorem \ref{th-existence-local-inertial-manifolds-skew-products}, the pullback attractor 
    $\{\mathcal{A}(t):t\in \mathbb{R}\}$ for $\{T(t,\tau): t\geqslant \tau\}$ 
    satisfies
	$\mathcal{A}(t)\subset \mathcal{M}(t)$ for every $t\in \mathbb{R}$.
%	and $\mathcal{A}\subset\bigcup_{t\in \mathbb{R}}\mathcal{M}(t)$.
\end{remark}
\end{comment}

\section{Applications}\label{Applications}

\subsection{The Saddle Point Property}\label{SPP}

We now obtain the saddle point property as an immediate consequence of Theorems \ref{NLTum} and \ref{NLTsm}. 

We respectively define the {\bfseries unstable} and {\bfseries stable sets} of a hyperbolic global solution $u_*$ of \eqref{NLeq:glin} as
\begin{subequations}
\begin{align}
W^u(u_*)&:= \left\{ (\tau,u_0)\in \mathbb{R}\times X :
\begin{array}{c}
 \text{ there is a solution } u:(-\infty,\tau]\to X \text{ such that }  \\
u(\tau)=u_0 \text{ and } \lim_{t\to -\infty} \|u(t)-u_*(t)\|_X=0
\end{array}
\right\}, \label{defWu} \\
W^s(u_*)&:= \left\{ (\tau,u_0)\in \mathbb{R}\times X :
\begin{array}{c}
 \text{ there is a solution } u:[\tau,\infty)\to X \text{ such that }  \\
u(\tau)=u_0 \text{ and } \lim_{t\to +\infty} \|u(t)-u_*(t)\|_X=0
\end{array}
\right\}. 
\end{align}
\end{subequations}
%
%\begin{definition}
%	The unstable manifold of a hyperbolic solution $\xi^*$ is the set
%	\begin{eqnarray*}
%		W^u(\xi^*)=\bigg\{(\tau,z)\in \mathbb{R}\times X: \hbox{ there is a backward solution } \zeta \\
%		\hbox{ satisfying } \zeta(\tau)=z \hbox{ and such that }
%		\lim_{t\to -\infty} \|\zeta(t)-\xi^*(t)\|_X=0
%		\bigg\}.
%	\end{eqnarray*}
%\end{definition}
%
%\begin{definition}
%	The stable manifold of a hyperbolic solution $\xi^*$ is the set
%\begin{eqnarray*}
%		W^s(\xi^*)=\bigg\{(\tau,z)\in \mathbb{R}\times X: \hbox{ there is a forward solution } \zeta \\
%		\hbox{ satisfying } \zeta(\tau)=z \hbox{ and such that }
%		\lim_{t\to +\infty} \|\zeta(t)-\xi^*(t)\|_X=0
%		\bigg\}.
%	\end{eqnarray*}
%\end{definition}
%

The existence, continuity and local characterization as graphs of the unstable and stable manifolds for non-autonomous hyperbolic asymptotic profiles was proved in \cite{Carvalho-Langa-07}. We obtain their results as a corollary of Theorems \ref{NLTum} and \ref{NLTsm}.
%In fact, we obtain that the unstable and stable manifolds of $\xi^*$ are given the graph of the maps 
%\begin{equation*}
%\begin{aligned}
%\mathbb{R}\times X\ni (\tau,u) \mapsto \Sigma^u(t,Q(t)u)\in Im(I-Q(t)), \\
%\mathbb{R}\times X\ni (\tau,u) \mapsto \Theta^s(t,u-Q(t)u) \in Im(Q(t)).
%\end{aligned}
%\end{equation*}
We only construct the unstable and stable manifolds for the trivial equilibria $u_*=0$, as any other global hyperbolic solution can be shifted to zero. See Remark \ref{Rmk:shift}.

%%%%%%%%%%%%%%%%%%%%%%%%%%%%%%%%%%
\begin{comment}
i.e., the points of the unstable manifold will be those of the form
\begin{equation*}
    \{(t,Q(t)z+\Sigma^u(t,Q(t)u)): (t,u): \mathbb{R}\times X, \ \ u \hbox{ small,}
\end{equation*}
and the points of the stable manifold will be those of the form
\begin{equation*}
    \{(t,\Theta^s(t,(I-Q(t))u))+(I-Q(t))u): (t,u): \mathbb{R}\times X, \ \ u \hbox{ small,}.
\end{equation*}
\end{comment}
%%%%%%%%%%%%%%%%%%%%%%%%

\begin{corollary}\label{NLTexp_SPP} 
Suppose that the linear evolution process $\{L(t,\tau): t\geqslant \tau\}$ has exponential dichotomy, with constant $M\geqslant 1$, exponent $\gamma>0$ and %$\rho=-\gamma$, and 
a family of projections $\{Q(t): \ t\in \R\}$.% and $f:\mathbb{R}\times X\to X$ with $Lip(f(t,\cdot)\leqslant \ell$. 

Suppose that $\ell>0$ is sufficiently small,
then there are continuous functions $\Sigma^u\in \mathcal{LB}_\Sigma(\kappa)$ and $\Theta^s\in \mathcal{LB}_\Theta(\kappa)$ such that the unstable and stable manifolds of $u_*=0$ are given by
\begin{subequations}
\begin{align}
W^u(0)&=\{(\tau,u)\in \R\times X: u=Q(\tau)
u+\Sigma^{u}(\tau, Q(\tau)u)\}, \label{eq-unstable-graph}  \\
W^s(0)&=\{(\tau,u)\in \R\times X: u=\Theta^s(\tau,(I-Q(\tau))u)+(I-Q(\tau))u \}.
\end{align}
\end{subequations}
Moreover, solutions within the unstable (resp. stable) manifold exponentially decay to zero backwards (resp. forwards) in time, according to \eqref{invSigmaNL} and \eqref{DecayStM}.
\end{corollary}

\proof 
For $\ell>0$ sufficiently small, the condition \eqref{NLeq:inequalities} is satisfied and $\delta>0$, and thus we obtain the graph of $\Sigma^*$ from Theorem \ref{NLTum}. We now prove that the unstable set $W^u(0)$ defined in \eqref{defWu} coincides with the graph of $\Sigma^u:=\Sigma^*$. On one hand, the graph of $\Sigma^u$ is contained in the unstable set by \eqref{invSigmaNL}. On the other hand, any solution $z:(-\infty,t]\to X$ which backwards converges to zero satisfies, from \eqref{expattSigmaNL},
\begin{equation*}
 \| z(t)-P_{\Sigma^*}(t)z(t)\|=\|(I-Q(t))z(t)-\Sigma^u(t,Q(t)z(t)) \| \leqslant M\| (I-P_{\Sigma^*}(\tau))z(\tau)\| e^{-\delta (t-\tau)}, \quad t\geqslant \tau.
\end{equation*}
Since $\delta>0$,
we obtain that $(I-Q(t))z(t)=\Sigma^u(t,Q(t)z(t))$ for all $t\in \mathbb{R}$ as $\tau \to -\infty$, and thus any element in the unstable set lies in the graph of $\Sigma^u$. 
%Note that $\delta>0$ in \eqref{deltaNL} for sufficiently small $\ell>0$.
The case of stable manifold is analogous applying Theorem \ref{NLTsm}.  
%
%From Theorem \ref{NLTum} applied for $-\rho=\gamma$ and $\ell>0$, there exists $\Sigma^*:=\Sigma^u\in \mathcal{LB}(\kappa)$ such that the family  $\{Q(\tau)u + \Sigma^{u}(\tau, Q(\tau)u): \ \tau\in \R\}$  is invariant for $\{T(t,\tau): t\geqslant \tau\}$ and satisfies \eqref{invSigmaNL} and \eqref{expattSigmaNL}.
%
%Let $z: \mathbb{R}\to X$ be a global solution of $\{T(t,\tau): t\geqslant \tau\}$ bounded as $t\to -\infty$, then from 
%\eqref{expattSigmaNL} we see that
%\begin{equation*}
%    \|(I-Q(t))z(t)-\Sigma^u(t,Q(t)z(t))\|_X\leqslant e^{-\delta(t-\tau)}      \|(I-P(\tau))z(\tau)\|_X, \ \ t\geqslant \tau.
%\end{equation*}
%Letting $\tau \to -\infty$ we obtain that $(I-Q(t))z(t)=\Sigma^u(t,Q(t)z(t))$, for all $t\in \mathbb{R}$.
%In particular, 
%\begin{equation}
%    W^u(0)\subset \{(\tau,u)\in \R\times X: u=Q(\tau) u+\Sigma^{u}(\tau, Q(\tau)u), \ \tau\in \R\}.
%\end{equation}
%
%Reciprocally, if $\eta\in Im(P_{\Sigma^u}(\tau))$, there exists $\zeta(t)=T(t,\tau)\eta$, for $t\in \mathbb{R}$ is well defined and from  \eqref{invSigmaNL}
%\begin{equation*}
%    \|\zeta(t)\|_{X} \leqslant M(1+\kappa)e^{(\gamma-\ell M (1+\kappa))(t-\tau)}\|u\|_{X},\quad t\leqslant \tau.
%    \end{equation*}
% Thus, $\zeta(t)\to 0$ as $t\to -\infty$ and the proof is complete.
\cqd

%Now, we consider the case $\xi^*\neq 0$, by the change of variables $v=u-\xi^*$, we translate all the dynamics around  $\xi^*$ to the zero solution, so that we consider 
%\begin{equation*}
%\dot v = A(t)v +B(t)v+ g(t,v),\ t>\tau,
%\end{equation*}
%where $B(t):=f_u(t,\xi^*(t))$ and $g(t,v):=f(t,\xi^*(t)+v)-f(t,\xi^*(t))-f_v(t,\xi^*(t))v$.  We observe that $g(t,0)=0\in X$ and $g_v(t,0)=0\in \mathcal{L}(X)$. In particular, in a small neighborhood $V$ of $0\in X$, $g$ is Lipschitz with constant with small constant $\ell>0$.  In this case it is possible to extend $g$ outside $V$ in such way that this extension is globally Lipschitz with the same Lipschitz constant $\ell>0$, see for instance \cite[page 631]{Carvalho-Langa-07}. At this scenario, we obtain the existence of \textit{local unstable} and  {local stable} manifolds. \textcolor{red}{this is already in Remark 3.1}

\begin{remark}
Corollary \ref{NLTexp_SPP} assumes that the nonlinearity $f(t,\cdot)$ is globally uniformly (in $t$) Lipschitz, and thereby yields global unstable and stable manifolds. However, in general, we have that the nonlinearity is only locally Lipschitz on a neighborhood of some equilibrium (here assumed to be zero after suitable translation), $O_\epsilon(0)$, and thus by means of a cut-off outside this region, we obtain maps, $\Sigma^u:\mathbb{R}\times O_{\epsilon}(0)\to X$ and $\Theta^s:\mathbb{R}\times O_{\epsilon}(0)\to X$, whose graphs contain the local unstable and stable manifolds, $W^u_{\text{loc}}(0)$ and $W^s_{\text{loc}}(0)$, respectively. 
Indeed, due to the respective growth bounds \eqref{invSigmaNL} and \eqref{DecayStM}, there is a neighborhood $O_{\epsilon_*}(0)$ contained in the cut-off region $O_{\epsilon}(0)$ such that %backward or forward 
orbits in the smaller neighborhood $O_{\epsilon_*}(0)$ do not leave the bigger neighborhood $O_{\epsilon}(0)$ in the appropriate time direction.
See the proof of Theorem \ref{th-existence-local-inertial-manifolds-skew-products} for more details on such cut-off of the nonlinearity to obtain local graphs, bearing in mind that after the procedure of shifting a solution to zero and obtaining a nonlinearity $g(t,u)$ in Remark \ref{Rmk:shift}, one can construct a cut-off function with sufficiently small $\ell>0$ on some neighborhood $O_{\epsilon}(0)$, since $g(t,0)=0$ and $g_u(t,0)=0$. See \cite[Theorem 6.1]{Carvalho-Langa-07} for more details.
\end{remark}

\subsection{Fine Description Within Invariant Manifolds }\label{APPconv}

We describe a finer growth and decay structure within invariant manifolds obtained in Section \ref{SPP}, in case of an additional exponential splitting. This allows the comparison between two different growth (resp. decay) rates within the unstable (resp. stable) manifold, which dictates the directions and the preferred directions along which solutions may approach the equilibria backwards in time (resp. forwards in time). We then apply this result to the example of asymptotically autonomous PDEs, which generalizes the well-known case of autonomous equations, see \cite[Lemma 2.2]{BrunovskyFiedler86} and \cite[Lemma 6]{Angenent86}. %As a standing hypothesis in this Section, we consider only autonomous diffusion, $A(t)\equiv A$, in equation \eqref{NLeq:glin}.

Often, exponential splitting for a linear evolution process $\{L(t,\tau): t\geqslant \tau\}$  occur with several different exponents and families of projections. The splitting given by the different projections are well behaved (nested) as it is proved in the next lemma.
\begin{lemma}\label{lemma_projections} 
Suppose that the linear evolution process $\{L(t,\tau): t\geqslant \tau\}$ has exponential splitting, with constant $M\geqslant 1$, exponents $\gamma>\rho$, and family of projections $\{Q(t): \ t\in \R\}$.
Assume further that $\{L(t,\tau): t\geqslant \tau\}$ has another exponential splitting with constant $M_*\geqslant 1$, %of the negative spectrum of the linear operator $A$, i.e., the hypothesis of Theorem \ref{NLTum} are also valid for some 
with exponents $\gamma_* >\rho_*$ such that $\gamma > \gamma_*$ and family of projections $\{Q_*(t): \ t\in \R\}$.
Then $Im(Q_*(t))\subset Im(Q(t))$ and $\mathrm{Ker}(Q(t))\subset \mathrm{Ker}(Q_*(t))$.
\end{lemma}

\proof For $\mu=\min\{\gamma,\frac{\gamma_*-\rho_*}{2}\}$, we have that
\begin{subequations}
    \begin{align}
    Im(Q_*(\tau))&=\{ u\in X:\ e^{( \gamma_* -\mu)(t-\tau)}L(t,\tau)u\in X \hbox{ is defined for all } t\leqslant \tau \hbox{ and is bounded} \},\\
    Im(Q(\tau))&=\{ u\in X:\ e^{(\gamma -\mu)(t-\tau)} L(t,\tau)u\in X \hbox{ is defined for all } t\leqslant \tau \hbox{ and is bounded} \},\\
    \mathrm{Ker}(Q_*(t))&=\{ u\in X:\ e^{(\gamma_* -\mu)(t-\tau)}L(t,\tau)u\in X  \hbox{ is bounded for all } t\geqslant \tau\},\\
    \mathrm{Ker}(Q(t))&=\{ u\in X:e^{(\gamma -\mu)(t-\tau)}L(t,\tau)u\in X  \hbox{ is bounded for all } t\geqslant \tau\},
    \end{align}
\end{subequations}
see \cite[Theorem 7.12]{Carvalho-Langa-Robinson-13}.
Clearly $Im(Q(\tau))\supset Im(Q_*(\tau))$ and $\mathrm{Ker}(Q(t))\subset \mathrm{Ker}(Q_*(t))$.\cqd

%{\color{black} We do not need $\gamma>0$. We only need that solutions are globally defined in the intertial manifold to prove the result below.
%Correct the result for as it is it only talks about solutions in the manifolds and we actually want to discuss solutions that are not necessarily in the manifolds}.

%\eject

\begin{corollary}\label{CORFinestructure} 
Suppose that the linear evolution process $\{L(t,\tau): t\geqslant \tau\}$ has exponential splitting, with constant $M\geqslant 1$, exponents $\gamma>\rho$, and family of projections $\{Q(t): \ t\in \R\}$.
Assume further that $\{L(t,\tau): t\geqslant \tau\}$ has another exponential splitting with constant $M_*\geqslant 1$, with exponents $\gamma_*>\rho_*$ such that
$\rho\geqslant \gamma_*$
 and family of projections $\{Q_*(t): \ t\in \R\}$.
%Assume further that the positive spectrum $\sigma_+(A)$ admits a further splitting at $\gamma_*\in (\gamma,\beta)$ with projections $\{Q_*(t): \ t\in \R\}$.

If $\ell>0$ is sufficiently small, then there are graphs corresponding to the fast and slow submanifolds represented by $W^{u}_{\text{fast}}(0),W^{u}_{\text{slow}}(0)$ $(W^{u}_{\text{fast}}(0)(t)$ is tangent to $Im(Q_*(t))$ and $(W^{u}_{\text{slow}}(0)(t)$ is tangent to $Ker(Q_*(t)Q(t)))$ of the invariant manifold represented by $W^u(0)$ $(W^u(0)(t)$ is tangent to $Im(Q(t)))$ such that %which are respectively tangent to $Im(Q(t)),Im(Q_*(t))\subseteq Im(Q(t)),Im(Q^c_*(t)):=Im(Q(t))\backslash Im(Q_*(t))$ such that
\begin{equation}
\lim_{t\to -\infty}\frac{\|(I-P_{slow}(\tau))T(\tau,t) u\|}{\|(I-P_{fast}(\tau))T(\tau,t) u\|}=0,
\end{equation}
for any $(\tau,u)\in W^u(0)\backslash W^{u}_{\text{fast}}(0)$, where $P_{fast}(\cdot)$ and $P_{slow}(\cdot)$ denote the respective nonlinear projections onto $W^{u}_{\text{fast}}(0)$ and $W^{u}_{\text{slow}}(0)$. %the. Note that $T(t,\tau)Q_*(\tau)u\in W^{u}_{\text{fast}}(0)$ and $T(t,\tau)Q^c_*(\tau)u\in W^{u}_{\text{slow}}(0)$. 
\end{corollary}

\begin{figure}[h]
 \begin{tikzpicture}[scale=0.5]

 %projections
 \draw[lightgray] (4.7,0) -- (-5,0) node[left] {\tiny{$_{Im(Q_*(t))}\!\!$}};
\draw[lightgray] (0,4.7) -- (0,-5) node[below] {\tiny{$^{Ker(Q_*(t)Q(t))}$}};

\draw[thin,lightbrown2,domain=-3:3,samples=100] plot ({0.5*\x^2},\x);
\draw[thin,lightbrown2,domain=-3:3,samples=100] plot (-{0.5*\x^2},\x);
 \draw[thin,lightbrown2,domain=-2.2:2.2,samples=100] plot ({\x^2},\x); 
 \draw[thin,lightbrown2,domain=-2.2:2.2,samples=100,] plot (-{\x^2},\x); 

 \draw[thick,lightgreen,domain=-5:5,samples=100,] plot (-{0.007*\x^3},\x) node[right] {\footnotesize{$_{\!\!W^u_{\!slow}(0)}$}}; 
\draw[thick,lightgreen,domain=- 5:5,samples=100] plot (\x,{0.007*\x^3})node[below] {\footnotesize{$_{\ \ \, W^u_{\!fast}(0)}$}};
\filldraw (0,0) circle (2.5pt);% node[anchor=south west]{\footnotesize{$0$}};
  
  \draw [-latex,lightgreen](-0.05,1.9) -- (-0.06,2);
  \draw [-latex,lightgreen](0.05,-1.9) -- (0.06,-2);

\draw [-latex,lightgreen](1.9,0.047) -- (2,0.054);
\draw [-latex,lightgreen](2.4,0.095) -- (2.5,0.105);
\draw [-latex,lightgreen](-1.9,-0.05) -- (-2,-0.054);
\draw [-latex,lightgreen](-2.4,-0.095) -- (-2.5,-0.105);

\draw [-latex,lightbrown2](1.9,1.95) -- (2.0,2.0);
\draw [-latex,lightbrown2](-1.9,-1.95) -- (-2.0,-2.0);
\draw [-latex,lightbrown2](1.9,-1.95) -- (2.0,-2.0);
\draw [-latex,lightbrown2](-1.9,1.95) -- (-2.0,2.0);

\draw [-latex,lightbrown2](1.92,1.37) -- (1.95,1.38);
\draw [-latex,lightbrown2](-1.92,1.37) -- (-1.95,1.38);
\draw [-latex,lightbrown2](1.92,-1.37) -- (1.95,-1.38);
\draw [-latex,lightbrown2](-1.92,-1.37) -- (-1.95,-1.38);

\end{tikzpicture}
\caption{The local dynamics inside the unstable manifold of $u_*\equiv 0$. 
Note that solutions in $W^u(0)\backslash W^u_{fast}(0)$ are tangent space of $W^u_{slow}(0)$ as $t\to -\infty$. 
} \label{FIGunstable}
\end{figure}
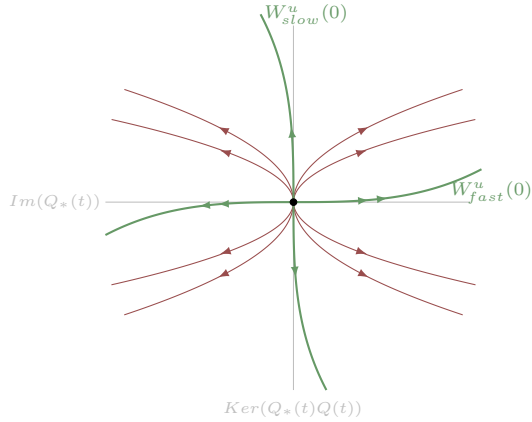

{\color{black}
\proof
For $(\tau, u_0)\in W^u(0)\setminus W^u_{fast}(0)$, we obtain from  \eqref{expattSigmaNL} and \eqref{attraction_stable} that
\begin{equation}\label{Ab_Eq1}
  \|(I-P_{fast}(t))T(t,\tau) u_0\| \leqslant Me^{-\delta_*(t-\tau)}\|(I-P_{fast}(\tau))u_0\|,\quad t\geqslant \tau.
\end{equation}
where $ \delta_*=\gamma_*-M\ell -\frac{ M^2\ell^2(1+\kappa)(1+M)}{\gamma_*-\rho_*-\ell M (1+\kappa)}$.
Similarly, from and \eqref{attraction_stable} we have that 

\begin{equation}\label{Ab_Eq2} 
    \|(I-P_{slow}(\tau))u_0\| \leqslant Me^{\hat{\delta}_*(t-\tau)}\|(I-P_{slow}(\tau))T(t,\tau)u_0\|, \quad t\geqslant \tau.
\end{equation}
where $\hat{\delta}_*=\rho_* + M\ell +\frac{M^2\ell^2(1+\kappa)(1+M)}{\gamma_*-\rho_*-M\ell (1+\kappa)}$ .

%{\color{red}
%
%\begin{equation}\label{Ab_Eq1}
%  \frac{1}{\|(I-P_{fast}(\tau))T(\tau,t)u_0\|}\leqslant \frac{Me^{-\delta_*(t-\tau)}}{\|(I-P_{fast}(t))u_0\|} ,\quad t\geqslant \tau.
%\end{equation}
%where $ \delta_*=\gamma_*-M\ell -\frac{ M^2\ell^2(1+\kappa)(1+M)}{\gamma_*-\rho_*-\ell M (1+\kappa)}$.
%Similarly, from and \eqref{attraction_stable} we have that 
%
%\begin{equation}\label{Ab_Eq2} 
%    \|(I-P_{slow}(\tau))T(\tau,t)u_0\| \leqslant Me^{\hat{\delta}_*(t-\tau)}\|(I-P_{slow}(\tau))u_0\|, \quad t\geqslant \tau.
%\end{equation}
%where $\hat{\delta}_*=\rho_* + M\ell +\frac{M^2\ell^2(1+\kappa)(1+M)}{\gamma_*-\rho_*-M\ell (1+\kappa)}$ .
%
%
%
%}

Therefore, the bounds \eqref{Ab_Eq1} and \eqref{Ab_Eq2}  applied to $u_0=T(\tau,t) u$ yield
\begin{equation}
\frac{\|(I-P_{slow}(\tau))T(\tau,t) u\|}{\|(I-P_{fast}(\tau))T(\tau,t) u\|}\leqslant M^2e^{-\bar{\delta}_*(t-\tau)} \frac{\|(I-P_{slow}(t))u_0\|}{\|(I-P_{fast}(t))u_0\|},\quad t\geqslant \tau.
\end{equation}
where $\bar{\delta}_*=\gamma_*-\rho_*-2M\ell -\frac{ 2M^2\ell^2(1+\kappa)(1+M)}{\gamma_*-\rho_*-\ell M (1+\kappa)}$.
The limit $\tau\!\to\!-\infty$ yields the desired claim, since $\bar{\delta}_*> 0$, for suitably small $\ell\!>\!0$.
\cqd

}

\begin{remark} 
Consider a small non-autonomous perturbation of a scalar semilinear parabolic equation in one spatial dimension,
\begin{equation}\label{na-pert}
\begin{split}
&u_t=u_{xx}+f_\nu(t,x,u),\ x\in (0,\pi),\\
&u(t,0)=u(t,\pi)=0,
\end{split}
\end{equation}
for $x\in (0,\pi)$ with Dirichlet boundary conditions, $\nu\in [0,1)$ and $f_\nu(t,x,0)=0$ for all $t\in\mathbb{R}$, $x\in (0,\pi)$ and $\nu\in [0,1)$. 
Assume that %$0\notin \sigma(\partial_{xx} + f'_{0}(x,0)\cdot I)$ and \textcolor{blue}{we dont need this... we already assume below that all equilibria are hyperbolic...},
\begin{equation}
\lim_{\nu\to 0} \sup_{t\in \R}\sup_{x\in [0,\pi]}\sup_{|u|\leqslant R}\{|f_\nu(t,x,u) -f_0(x,u)|+|\partial_uf_\nu(t,x,u) -\partial_uf_0(x,u)|\}=0  
\end{equation}
for all $R>0$.
Assume further that $f_{0}(x,u)$ satisfy dissipative growth conditions and thus the limiting autonomous equation with nonlinearity $f_0(x,u)$ possess a gradient global attractor. Suppose that all equilibria of the limiting autonomous equations are hyperbolic. If $\nu$ is suitably small, associated to each hyperbolic equilibria there is a hyperbolic global solution and the non-autonomous evolution process associated to \eqref{na-pert} is dinamically gradient, see \cite{BCCP}.

Due to Remark \ref{Rmk:shift}, without loss of generality, we analyse the local dynamics nearby $u_*\equiv 0$. Consider the linear operator $A:=\partial_{xx}+f '_0(x,0)I$, which is Sturm-Liouville. Hence its spectrum consists of simple eigenvalues $\{\lambda_k\}_{k\in\mathbb{N}}$ with corresponding eigenfunctions denoted by $\{\phi_k\}_{k\in\mathbb{N}}$. 
%
%In particular, since $0\notin\sigma(\partial_{xx}+f'_{0}(0)I)$, $u_*\equiv 0$ is a hyperbolic equilibrium for 
%\begin{equation}\label{na-ci}
%\begin{split}
%&u_t=u_{xx}+f_0(x,u),\ x\in (0,\pi),\\
%&u(0)=u(\pi)=0
%\end{split}
%\end{equation}
Let $n\in \N$ be the Morse index of the hyperbolic equilibrium $u_*\equiv 0$ and thus $\lambda_1,\ldots,\lambda_n>0$, whereas $\lambda_{n+1},\lambda_{n+2},\ldots<0$.

For $k\leqslant n$, consider the splitting $\sigma(A)\!:=\!\sigma_k^+ \cup \sigma_k^-$, where $\sigma_k^+\!=\!\{\lambda_1,...,\lambda_k\}$ and $\sigma_k^-\!=\!\{\lambda_{k+1},\lambda_{k+2},...\}$. Let $Q^0_k$ be the spectral projection associated to $\sigma_k^+$. For sufficiently small $\nu\in [0,1)$, there are projections $Q^\nu_k(t)$ (which are suitably close to the orthogonal projections $Q^0_k$) associated to the exponential splitting of the linear process described by
\begin{equation*}
\begin{split}
&u_t=u_{xx}+f'_0(x,0)u+[f'_\nu(t,x,0)-f'_0(x,0)]u, \ x\in (0,\pi)\\
&u(0)=u(\pi)=0.
\end{split}
\end{equation*}
After an appropriate change of coordinates, we can apply Corollary \ref{NLTexp_SPP} for any $k=1,...,n$, yielding the $k$-fast unstable manifold, $W^u_k(0)\subseteq W^u(0)$, which is tangent to $Im(Q_k^\nu(t))$.
Moreover, we can also apply Corollary \ref{CORFinestructure} to determine that solutions in $W^u(0)\backslash W^u_k(0)$ are tangent to the lower dimensional subspaces $Im(Q^\nu(t))\backslash Im(Q^\nu_k(t))$, where $Q^\nu(t)$ is the projection associated with the exponential dichotomy that splits the positive and negative eigenstructure. Due to the continuity of the projections we conclude that all solutions in the unstable manifold approach the equilibria in directions with are suitably close to the directions of the eigenfunctions of the corresponding linearization. 

As a consequence, solutions of equations with asymptotically autonomous nonlinearities approach the equilibria in the direction of an eigenfunction associated to the linearization of the limiting problem around the associated equilibria.

\end{remark}

%{\color{blue} I claim that this is a consequence of the inertial manifold theorem for, the eigenvalues are all simple and $f(0)=0$, $f'(0)=0$. This allow us to make the exponential splitting at any point between two eigenvalues, choose $\kappa>0$, then choose the neighborhood (and consequently $L$) so small that all conditions of the inertial manifold are satisfied. $\kappa$ arbitrary ensures that $\Sigma$ and $\Theta$ are zero and have zero derivatives at zero. One can proceed now by making the first splitting between $\lambda_1$ and $\lambda_2$, restricting the analysis to the graph of $\Theta$, one can make the splitting between $\lambda_2$ and $\lambda_3$ and so on.}

\subsection{Roughness of Exponential Dichotomy}\label{NLSrough}

In this section, we prove that the roughness of exponential dichotomy, i.e., that exponential dichotomies are preserved under small perturbations.

We first assume that the linear evolution process $\{ L(t,\tau): t\geqslant \tau\}$ associated to the problem
\begin{equation}\label{Leq:glin0}
\dot{u}=A(t)u, \ t\geqslant \tau, \  \  u(\tau)=u_0.
\end{equation} 
has exponential dichotomy with constant $M$ and exponent $\gamma>0$ and then consider the linear evolution process $\{ T(t,\tau): t\geqslant \tau\}$, associated to a perturbation of it, given by the following linear equation,
\begin{equation}\label{Leq:glin}
\dot{u}=A(t)u+B(t)u, \ t\geqslant \tau, \ \  u(\tau)=u_0.
\end{equation} 
where the map $t\mapsto B(t)\in \mathcal{L}(X)$ is strongly continuous for $t\in\mathbb{R}$ and ${\sup_{t\in \R}}\|B(t)\|_{\mathcal{L}(X)}\leqslant \ell $, for some suitably small $\ell>0$.
Recall that, as in \eqref{nonlinEP}, the evolution process $\{T(t,\tau): t\geqslant \tau\} \subset \mathcal{L}(X)$ associated to \eqref{Leq:glin} is given by
\begin{equation}\label{eq-linear-evo-proc-T-definition}
T(t,\tau)=L(t,\tau)+\int_\tau^t L(t,s)B(s)T(s,\tau)\,ds, \quad t\geqslant \tau.
\end{equation}
We wish to prove that \eqref{Leq:glin} has exponential dichotomy for suitably small $\ell$.

This result can be obtained by firstly applying Theorems \ref{NLTum} and \ref{NLTsm} in a linear setting, which are suitable in order to establish the existence of the linear invariant manifold and its stable manifold (see Corollary \ref{LCum}) and then apply it to \eqref{Leq:glin} with $\gamma>0$ and $\rho=-\gamma$. 
%\textcolor{red}{PL: add something about $A(t)$.}
%
%We can thereby obtain the stable manifold of a linear inertial manifold, which has been proven to exist in Corollary \ref{LCum}, if we consider solutions of the linear problem \eqref{Leq:glin}.

%\textcolor{red}{PL: to be honest, I'm not even sure we should write these corollaries below at all. they are just the theorems in the linear setting... in the proof of the roughness, we can simply mention that we apply the main theorems for the linear case, period.}

\begin{corollary}\label{LCum} 
Suppose that $\{L(t,\tau): t\geqslant \tau\}$ has exponential splitting with constant $M$, exponents $\gamma>\rho$ and family of projections $\{Q(t): \ t\in \R\}$ and that \eqref{NLeq:inequalities} is satisfied. The following holds:
\begin{itemize}
\item  There are maps $\Sigma^*,\Theta^*: \mathbb{R}\times X\to X$, 
$\Sigma^*(t,\cdot),\Theta^*(t,\cdot)\in \mathcal{L}(X)$ and $\|\Sigma^*(t, u)\|\leqslant \kappa \|u\|_{X}$, $\|\Theta^*(t, u)\|_{X} \leqslant \kappa \|u\|_{X}$ for all $(t,u)\in \R\times X$ and for some $\kappa=\kappa_\ell>0$;
\item The graph $\mathcal{G}(\Sigma^*)$ of  $\Sigma^*$ is an invariant family and  \eqref{expattSigmaNL} holds, the graph $\mathcal{G}(\Theta^*)$ of $\Theta^*$ is a positively invariant family;
\item The evolution process $\{T(t,\tau): t\geqslant \tau\}$ given by \eqref{eq-linear-evo-proc-T-definition} satisfies
\begin{subequations}\label{LCum-Inequalities}
\begin{align}
    \|T(t,\tau)P_{\Sigma^*} (\tau)\|_{\mathcal{L}(X)} &\leqslant M(1+\kappa) e^{-(\rho+ M\ell(1+\kappa))(t-\tau)},\quad t\leqslant \tau,\\
    \|T(t,\tau)P_{\Theta^*} (\tau)\|_{\mathcal{L}(X)} &\leqslant M(1+\kappa) e^{-(\gamma- M\ell(1+\kappa))(t-\tau)},\quad t\geqslant \tau ,
\end{align}
\end{subequations}
where $P_{\Sigma^*}(t)u\!:=\!Q(t)u\!+\!\Sigma^*(t,Q(t)u)$ and $P_{\Theta^*}(t)u\!:=\!\Theta^*(t,(I\!-\!Q(t))u)+(I\!-\!Q(t))u$, $t\in \mathbb{R}$.
\end{itemize}
\end{corollary}
\proof The proof is a direct consequence of Theorems \ref{NLTum} and \ref{NLTsm} in the case that $f(t,\cdot)$ is linear and uniformly (with respect to $t$) bounded. 
Note that, the linearity of $\Sigma(t,\cdot)$ follows since $f(t,\cdot)$ is linear, and thereby $G(\Sigma)$ given by equation \eqref{NLinvman7} is also linear. Consequently, the fixed point, $G(\Sigma^*)(t,u)=\Sigma^*(t,u)$, is linear. Similarly, $\tilde{G}$ in equation \eqref{Gtilde} is also linear and so is $\Theta^*$.
\cqd
%since the fixed point argument in Theorem \ref{NLTum} can be replicated in the following space consists of linear functions,
%\begin{equation*}
%\begin{split} 
%\mathcal{LB}_\Sigma(\kappa) 
%=\Big\{&\Sigma:\R\to \mathcal{L}(X): \  \Sigma\hbox{ is strongly continuous,}  \\
%&\Sigma(t,x)=(I-Q(t))\Sigma(t,Q(t)x)=\Sigma(t,Q(t)x) \hbox{ and } \| \Sigma(t,\cdot)\|_{\mathcal{L}(X)} \leqslant \kappa,\ \forall t\in \R\Big\}.
%\end{split}
%\end{equation*} 
%with metric $|\!|\!|\Sigma-\tilde{\Sigma}|\!|\!|={\displaystyle\sup_{t\in \R} \|\Sigma(t,\cdot)- \tilde{\Sigma}(t,\cdot)\|_{\mathcal{L}(X)}}$. This is a subsets of the previously chosen set in \eqref{FPSigma}. Moreover, if $f$ and $\Sigma(t,\cdot)$ are linear, so is $G(\Sigma)(t,\cdot )$.\cqd

%\proof The proof is a direct consequence of Theorem \ref{NLTsm} observing that the linearity of $\Theta(t,\cdot)$ can be archived from the uniqueness of the fixed points and simply choosing 
%\begin{equation*}
%\begin{split} 
%\mathcal{LB}_\Theta=\Big\{&\Theta:\R\to \mathcal{L}(X): \  \Theta\hbox{ is strongly continuous,}\\
%&\Theta(t,y)\!=\!Q(t)\Theta(t,(I-Q(t))y)\!=\!\Theta(t,(I-Q(t))y) \hbox{ and } \|\Theta(t,\cdot)\|_{\mathcal{L}(X)} \leqslant \kappa,\ \forall t\!\in\! \R\Big\}.
%\end{split}
%\end{equation*} 
%with the metric $|\!|\!|\Theta-\Theta'|\!|\!|={\displaystyle\sup_{t\in \R} \|\Theta(t,\cdot)- \Theta'(t,\cdot)\|_{\mathcal{L}(X)}}$ which are subsets of the previously chosen set and noting that if $f$ and $\Theta(t,\cdot)$ are linear, so is $G_\Theta(\Theta)(t,\cdot )$.\cqd

Next, we show the robustness of the exponential dichotomy. %, i.e., the persistence of the dichotomy of a linear evolution process $\{L(t,\tau): t\geqslant \tau\}$ generated by $A(t)$, in case of a bounded perturbation $B(t)$, which generates the linear evolution process denoted by $\{T(t,\tau): t\geqslant \tau\}$. 
%Moreover, we quantify how close the projections related to these dichotomies are, which depends on the size of the perturbation $B(t)$.

\begin{corollary}\label{NLTexp_dicho_roughness} 
Suppose that $\{L(t,\tau):t\geqslant \tau\}$ has exponential dichotomy with constant $M\geqslant 1$, exponent $\gamma>0$ and family of projections $\{Q(t):t\in \mathbb{R}\}$. 
Moreover, assume that $\sup_{t\in \mathbb{R}}\|B(t)\|_{\mathcal{L}(X)}\leqslant \ell$, where $\ell>0$ satisfies
\begin{equation}\label{sizeofperturb}
    \ell < \frac{2\gamma}{3M(M+1)}.
\end{equation}
%$\kappa=\kappa_\ell<1/2$  obtained in Corollary \ref{LCum}.
%\textcolor{red}{PL: there should be no $\kappa$ here... we should write the explicit condition on $\ell$ and how big this perturbation can be, i.e., $\ell < ...$}
%\textcolor{lightgreen}{AS: I believe the way that is written is ok, since $\kappa_l\to 0,$ as $\ell \to 0$, it is reasonable to make $\ell$ suitable small such that $\kappa_\ell <1/2$. Explicitly the condition on $\ell$ is not so nice, but I have the computations, if it is the case we can analyze and put in the theorem.}
%
Then $\{T(t,\tau):t\geqslant \tau\}$ has exponential dichotomy, i.e., 
there is a family of projections $\{Q_\ell (t):t\in \R\}$ such that $T(t,\tau):Im(Q_\ell (\tau))\to Im(Q_\ell (t))$ is an isomorphism for $t\geqslant \tau$, with inverse denoted by $T(\tau,t)$, and 
\begin{equation}\label{RoughDichoEq}
\begin{aligned}
    \|T(t,\tau)Q_\ell (\tau)\|_{\mathcal{L}(X)} &\leqslant M_\ell e^{\gamma_\ell (t-\tau)},\quad t\leqslant \tau\\
    \|T(t,\tau)(I-Q_\ell (\tau))\|_{\mathcal{L}(X)} &\leqslant M_\ell e^{-\gamma_\ell (t-\tau)},\quad t\geqslant \tau ,
\end{aligned}
\end{equation}
where $M_\ell:=M(1+\kappa_\ell)/(1-2\kappa_\ell)>1$ and $\gamma_\ell:=\gamma-\ell M(1+\kappa_\ell)>0$ for the Lipschitz constant $\kappa_\ell$ obtained in Corollary \ref{LCum}.
Moreover, 
\begin{equation}\label{NLcont_proj}
\sup_{t\in \mathbb{R}}\|Q(t)-Q_\ell (t)\|_{\mathcal{L}(X)}\leqslant \frac{2\kappa_\ell}{1-2\kappa_\ell}. %\frac{\ell M^2 (1+\kappa )}{\gamma-2M\ell (1+\kappa)}.
\end{equation} 
\end{corollary}

\begin{remark}
Note that Corollary \ref{NLTexp_dicho_roughness} explicitly provides the constant and exponent of the exponential dichotomy for $\{T(t,\tau): t\geqslant \tau\}$ depending of $\ell$. 
Moreover, the constant $\kappa_\ell$ in inequality \eqref{NLcont_proj} is given by smallest possible Lipschitz constant $\kappa_-$ explicitly obtained in \eqref{Delta^*}, which occurs in the equality of equation \eqref{nuNL}.
In particular, the hypothesis \eqref{sizeofperturb} holds in case that $\ell>0$ is suitably small. Furthermore, note that $\kappa_\ell \to 0, M_\ell\to M, \gamma_\ell \to \gamma, Q_l\to Q$ as $\ell \to 0$, uniformly in $t$. 
\end{remark}

\proof 
%The linearity of $\Sigma(t,\cdot)$ and of $\Theta(t,\cdot)$ can be achived from the uniqueness of the fixed points and simply choosing 
%\begin{equation*}
%\begin{split} 
%\mathcal{LB}_\Sigma(\kappa)
%=\Big\{\Sigma:\R\to \mathcal{L}(X): \  \Sigma&\hbox{ is strongly continuous,}\ \Sigma(t,x)=(I-Q(t))\Sigma(t,Q(t)x)=\Sigma(t,Q(t)x),  \\
%&\hbox{and } \| \Sigma(t,\cdot)\|_{\mathcal{L}(X)} \leqslant \kappa,\ \forall t\in \R\Big\}.
%\end{split}
%\end{equation*} 
%with the metric $|\!|\!|\Sigma-\tilde{\Sigma}|\!|\!|={\displaystyle\sup_{t\in \R} \|\Sigma(t,\cdot)- \tilde{\Sigma}(t,\cdot)\|_{\mathcal{L}(X)}}$
%and
%\begin{equation*}
%\begin{split} 
%\mathcal{LB}_\Theta=\Big\{\Theta:\R\to \mathcal{L}(X): \  &\Theta\hbox{ is strongly continuous,}\  \Theta(t,y)=Q(t)\Theta(t,(I-Q(t))y)=\Theta(t,(I-Q(t))y)  \\
%&\hbox{and } \| \Theta(t,\cdot)\|_{\mathcal{L}(X)} \leqslant \kappa,\ \forall t\in \R\Big\}.
%\end{split}
%\end{equation*} 
%with the metric $|\!|\!|\Theta-\Theta'|\!|\!|={\displaystyle\sup_{t\in \R} \|\Theta(t,\cdot)- \Theta'(t,\cdot)\|_{\mathcal{L}(X)}}$
%which are subsets of the previously chosen sets and noting that if $f$ and $\Theta(t,\cdot)$ and $\Sigma(t,\cdot)$ are linear, so are $G_\Theta(\Theta)(t,\cdot )$ and $G_\Sigma(\Sigma)(t,\cdot )$.
%
%\bigskip
%
%
%
%In order to prove the evolution process  $\{T(t,\tau): t\geqslant \tau\}$ has exponential dichotomy with projections $\{Q_\ell (t):t\in \mathbb{R}\}$, 
Consider the linear maps  $P_{\Sigma^*}(t)u:=Q(t)u+{\Sigma^*} (t,Q(t)u)$ and 
$P_{\Theta^*}(t)u:=(I-Q(t))u+ {\Theta^*} (t,(I-Q(t))u)$, for $(t,u)\in \mathbb{R}\times X$,
which were obtained in Corollary \ref{LCum}, where $\Sigma^*$ and $\Theta^*$ are bounded linear maps that do not exceed $\kappa_\ell>0$.

\par We will prove that 
$X=Im(P_{\Sigma^*}(t))\oplus Im(P_{\Theta^*}(t))$, for every $t\in \mathbb{R}$.
Equivalently, we show that, for each $(t,u)\in \mathbb{R}\times X$, the operator $\mathcal{I}_u(t)$, defined by 
\begin{equation}
\begin{aligned}\label{Irough}
    \mathcal{I}_u(t): X &\to X\\
    v &\mapsto \mathcal{I}_u(t) v:=u-{\Sigma^*}(t,v)-{\Theta^*}(t,v),
\end{aligned}
\end{equation}
 admits a unique fixed point on $X$. In fact, if that is the case, for each $(t,u)\in \R\times X$, there exists a unique $v_u\in X$ such that $\mathcal{I}_u(t)v_u=v_u$, %\textcolor{red}{PL: I would say the notation $v_u$ is not very good... maybe $v_*$?} 
that is,
\begin{equation}
    u-{\Sigma^*}(t,v_u)-{\Theta^*}(t,v_u)=v_u=Q(t)v_u+(I-Q(t))v_u,
\end{equation}
or
\begin{equation}
    u=Q(t)v_u+{\Sigma^*}(t,v_u)+(I-Q(t))v_u+{\Theta^*}(t,v_u)=P_{\Sigma^*}(t)v_u+P_{\Theta^*}(t)v_u,
\end{equation}
Which is the unique representation of $u$ as a sum of elements of
$Im(P_{\Sigma^*}(t))$ and $Im(P_{\Theta^*}(t))$ and proves the desired decomposition.

\begin{figure}[H]
\minipage{0.49\textwidth}\centering
    \begin{tikzpicture}[scale=1.25]
    %coordinates
    \draw[lightbrown2] (-1.5,0) -- (1.5,0) node[right] {\tiny{$Q(t)$}};
    \draw[lightbrown2] (0,-1.5) -- (0,1.5) node[above] {\tiny{$I\!-\!Q(t)$}};    

    %graphs    
    \draw[rotate=20,lightgray] (-1.5,0) -- (1.5,0) node[right] {\tiny{$\Sigma^*(t,\!\cdot)$}};
    \draw[rotate=20,lightgray] (0,-1.5) -- (0,1.5) node[left] {\tiny{$\Theta^*(t,\!\cdot)$}};

 %projections    
    \draw[lightgray,shift={(0.6,1.25)},rotate=20] (0,0) -- (-1,0) node[left] {\tiny{$P_{\Theta^*}(t)u\!\!\!$}};
    \draw[lightgray,shift={(0.58,1.25)},rotate=20] (0,0) -- (0,-0.98) node[anchor=north] {\tiny{$P_{\Sigma^*}(t)u$}};

    %equilibria
    \filldraw[lightgray] (0.58,1.24) circle (0.5pt) node[anchor=west]{\color{black}\tiny{$\!u$}};
  
    %equilibria
    \filldraw[lightgray] (-0.33,0.91) circle (0.5pt) ;
     
     %equilibria
    \filldraw[lightgray] (0.0,0.0) circle (0.5pt) ;
     %equilibria
    \filldraw[lightgray] (0.91,0.33) circle (0.5pt);

    \end{tikzpicture}
\endminipage\hfill 
\minipage{0.49\textwidth}\centering
    \begin{tikzpicture}[scale=1.25]
    %projections    
    \draw[gray,shift={(0.6,1.25)},rotate=20] (0,0) -- (-1,0);% node[left] {\tiny{$P_{\Sigma^*}(t)u$}};
    \draw[gray,shift={(0.58,1.25)},rotate=20] (0,0) -- (0,-0.98);% node[below] {\tiny{$P_{\Theta^*}(t)u$}};    
    
    %coordinates
    \draw[lightbrown2] (-1.5,0) -- (1.5,0) node[right] {\tiny{$Q(t)$}};
    \draw[lightbrown2] (0,-1.5) -- (0,1.5) node[above] {\tiny{$I\!-\!Q(t)$}};    

    %graphs    
    \draw[rotate=20,gray] (-1.5,0) -- (1.5,0) node[right] {\tiny{$\Sigma^*(t,\!\cdot)$}};
    \draw[rotate=20,gray] (0,-1.5) -- (0,1.5) node[left] {\tiny{$\Theta^*(t,\!\cdot)$}};    

    %finding v_u
    \draw[lightbrown2,dashed] (0.91,0.91) -- (0.91,0) node[below] {\tiny{$Q(t)v_u$}};
    \draw[lightbrown2,dashed] (0.91,0.91) -- (0.0,0.91) node[left] {\tiny{$\!(I\!-\!Q(t))v_u\!\!$}};

  %equilibria
    \filldraw[gray] (0.58,1.24) circle (0.5pt) node[anchor=west]{\color{black}\tiny{$\!u$}};
 
     %equilibria
    \filldraw[gray] (-0.33,0.91) circle (0.5pt) ;
     
     %equilibria
    \filldraw[lightgray] (0.0,0.0) circle (0.5pt) ;
    
     %equilibria
    \filldraw[lightbrown2] (0.0,0.91) circle (0.5pt) ;

     %equilibria
    \filldraw[lightbrown2] (0.91,0.0) circle (0.5pt);      
    
     %equilibria
    \filldraw[gray] (0.91,0.33) circle (0.5pt);   
      
    %v_u
    \filldraw[lightbrown2] (0.91,0.91) circle (0.5pt) node[anchor=west]{\tiny{\color{black}$\!v_u$}};
    
    \end{tikzpicture}
\endminipage
\caption{Given a point $u\in X$, we find a unique point $v_u\in X$ such that $P_{\Sigma^*}(t)u=Q(t)v_u+\Sigma^*(t,Q(t)v_u)$ and $P_{\Theta^*}(t)u=(I-Q(t))v_u+\Theta^*(t,(I-Q(t))v_u)$. } \label{FIGfp}
\end{figure}
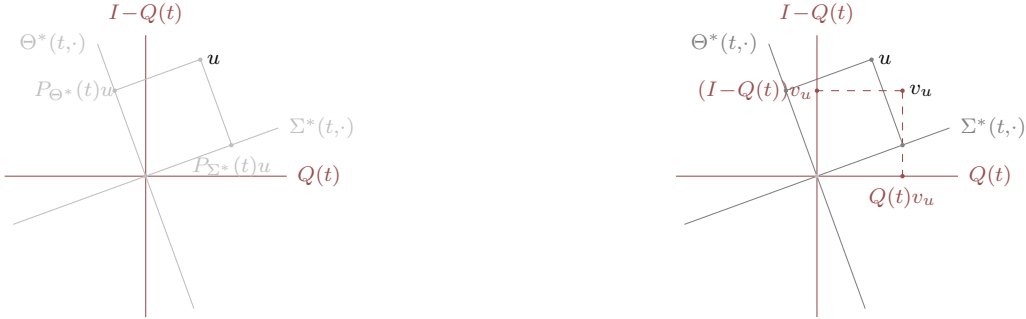

In order to show that $\mathcal{I}_u(t)$ has a unique fixed point, note that $\mathcal{I}_u(t)$ is a contraction on $X$, since
\begin{equation}\label{RoughContr}
\begin{aligned}
    \|\mathcal{I}_u(t)v-\mathcal{I}_u(t)\tilde{v}\|&=\|{\Sigma^*}(t,\tilde{v})-{\Sigma^*}(t,v)+{\Theta^*}(t,\tilde{v})-{\Theta^*}(t,v)\|,\\
    &\leqslant 2\kappa \|v-\tilde{v}\|,    
\end{aligned}
\end{equation}
for any $(t,v),(t,\tilde{v})\in \R\times X$, as ${\Sigma^*(t,\cdot)},{\Theta^*(t,\cdot)}$ are Lipschitz with constant $\kappa=\kappa_\ell>0$, for all $t\in \R$.
Thus, $\mathcal{I}_u(t)$ is a contraction, for each $t\in \mathbb{R}$ and $\kappa\in [\kappa_-,\min\{1/2,\min\{\kappa_+,\kappa_*\}\} )$, where $\kappa_-$ is given by \eqref{Delta^*}, since the hypothesis \eqref{sizeofperturb} implies that $\kappa_-<1/2$. Without loss of generality, we may choose $\kappa_\ell:=\kappa_-$.

\par Note that, for each $u\in X$, since $v_u$ is the unique element of $X$ satisfying $v_u=u-{\Sigma^*}(t,v_u)-{\Theta^*}(t,v_u)$, the map
$u\mapsto v_u$ is a bounded linear operator such that
\begin{equation}\label{eq-boundedness-for-v_u}
    \|v_u\|_X\leqslant \frac{\|u\|_X}{1-2\kappa},
\end{equation}
due to the triangle inequality.

\par For each $t\in \mathbb{R}$, define $Q_\ell (t)\in \mathcal{L}(X)$. The linear projection onto $Im(P_{\Sigma^*}(t))$ along $Im(P_{\Theta^*}(t))$, which can be written as $Q_\ell (t)u:=P_{\Sigma^*}(t)v_u$ due to the first part of the proof. Its complementary projection is given by $(I-Q_\ell (t))u=P_{\Theta^*}(t)v_u$, for each $(t,u)\in \mathbb{R}\times X$.

\par From Corollary \ref{LCum}, we have that
$\{Im(Q_\ell (t)): t\in \mathbb{R}\}$ is invariant and $\{Im(I-Q_\ell (t)): t\in \mathbb{R}\}$ is positively invariant. 
%i.e., $T(t,\tau)Im (Q_\ell (\tau))=Im(Q_\ell (t))$ and 
%$T(t,\tau)Im(I-Q_\ell (\tau))\subset Im (I-Q_\ell (t))$, 
Thus $T(t,\tau)Q_\ell (\tau)=Q_\ell (t)T(t,\tau)$, for every $t\geqslant \tau$. 
Equations \eqref{LCum-Inequalities} and \eqref{eq-boundedness-for-v_u} imply the bounds in \eqref{RoughDichoEq}.
This proves that $\{T(t,\tau):t\geqslant \tau\}$ has exponential dichotomy with constant $M_\ell:=M(1+\kappa)/(1-2\kappa)$ and exponent $\gamma_\ell := \gamma-\ell M(1+\kappa)>0$.

%
\begin{comment}

\textcolor{red}{PL: this was abrupt to me. suddenly you get the bounds below, which I have no clue where you obtained them... maybe cite a few equations you used, and how you get this... For me this is just a copy and paste of equation \eqref{RoughDichoEq}. If we are just gonna copy and paste, why not cite the equation and not repeating it?}

Define the projection $Q_\ell (t):X\to X$ given by
$$
R(P(t))=W^u(0,0)(t) \ \hbox{ e } \ R(I-P(t))=W^s(0,0)(t).
$$
\textcolor{red}{PL: I don't see why the second equality is true for this NONLINEAR projection. the sets $W^u$ and $W^s$ are NOT complementary!  }
%
In case that $-\rho=\gamma>0$, the problem \eqref{NLeq:glin} \textcolor{red}{what? this is NOT a theorem about the nonlinear problem! } has exponential dichotomy with projections $\{Q(t):t\in \R\}$, constant $M(1+\kappa)$ and exponent $\gamma-\ell M(1+\kappa)$, whenever the conditions of Theorems \ref{NLTum} and \ref{NLTsm} are satisfied, i.e., whenever \eqref{NLeq:inequalities} hold true with the additional assumption that $\kappa<1/2$.
Then, 
%
\begin{equation*}
\begin{split}
&\|T(t,\tau)Q_\ell (\tau)\|_{\mathcal{L}(X)} \leqslant \frac{M(1+\kappa)}{1-2\kappa} e^{(\gamma- M\ell(1+\kappa))(t-\tau)},\quad t\leqslant \tau\\
&\|T(t,\tau)(I-Q_\ell)\|_{\mathcal{L}(X)} \leqslant \frac{M(1+\kappa)}{1-2\kappa}e^{-(\gamma- M\ell(1+\kappa))(t-\tau)},\quad t\geqslant \tau.
\end{split}
\end{equation*}

\end{comment}

Lastly, we prove the bound in equation \eqref{NLcont_proj}, i.e., the continuous dependence of the projections $\{Q(t): t\in \mathbb{R}\}$ and $\{Q_\ell (t): t\in \mathbb{R}\}$, 
correponding to the exponential dichotomies of the respective evolution processes $\{L(t,\tau):t\geqslant \tau\}$ and $\{T(t,\tau):t\geqslant \tau\}$.%, with respect to the perturbation $B(t)$ with Lipschitz constant $\ell>0$ in equation \eqref{Leq:glin}. 

Consider $u\in X$, which can be uniquely decomposed as $u=v_u+{\Sigma^*}(t,v_u)+{\Theta^*}(t,v_u)$. Hence, $Q(t)u=Q(t)v_u+{\Theta^*}(t,v_u)$, since $Q(t){\Sigma^*}(t,v_u)=0$, and $Q_\ell(t)u=Q(t)v_u+{\Sigma^*}(t,v_u)$, by definition of $Q_\ell (t)$. Therefore,
\begin{equation}
Q(t)u-Q_\ell (t)u= {\Theta^*}(t,v_u)-{\Sigma^*}(t,v_u).
\end{equation}
Since ${\Sigma^*},{\Theta^*}$ are Lipschitz with constant $\kappa_\ell>0$, and due to equation \eqref{eq-boundedness-for-v_u}, we obtain \eqref{NLcont_proj}.
%\begin{equation}
%    \|Q(t)u-Q_\ell (t)u\|\leqslant \frac{2\kappa}{1-2\kappa}\|u\|.    
%\end{equation}
%Thus, equation \eqref{NLcont_proj} is obtained for some $\kappa\in (0,1/2)$ satisfying
%\begin{equation}
%\frac{\ell M^2 (1+\kappa )}{\gamma-M\ell (1+\kappa)}= 2\kappa.
%\end{equation}
%\textcolor{red}{PL: shouldn't we have an inequality $\leqslant$ actually? Similarly to the inertial manifold, the equality is the 'worst constant', which can be solved again by a quadratic equation... Also, I don't understand why this is needed at all? Can't we simply have equation (51)? Why is (45) better? Also, in eq (45) there is a 2 in the denominator, which I have no clue where it came from... Also, shouldnt we want $\frac{2\kappa}{1-2\kappa}$ in eq 52? }
\cqd

%\section{Parabolic PDEs}\label{sec:parabolicPDEs}

\subsection{Hyperbolic Solutions for a Nonautonomous PDE}\label{sec:parabolicPDEs}

We establish hyperbolicity for global solutions a non-autonomous parabolic PDE by reducing the infinite dimensional problem to a two-dimensional non-autonomous ODE with hyperbolic solutions. 

\par Consider the scalar parabolic non-autonomous equation
\begin{equation}\label{eq-applications-scalar-parabolic}
%\left\{ 
\begin{array}{l l} 
u_t=(a_\nu(x) u_x)_x+f(t,u), & \ x\in (0,1), \ t>\tau\\
u_x(t,0) = u_x(t,1)=0, & \  t>\tau \\
u(\tau,x)=u_0(x),&  \ x\in (0,1),
\end{array} 
%\right.
\end{equation}
where the reaction is given by $f(t,u):=u-\beta(t)u^3$, where $\beta:\mathbb{R}\to \mathbb{R}$ is an uniformly positive and bounded globally Lipschitz function. The diffusivity %, i.e. $\beta(\mathbb{R})\subset [\beta_0,\beta_1]$ for some $0<\beta_0<\beta_1<+\infty$.
$a_\nu \in C^2(\mathbb{R})$ %([0,1],\mathbb{R})$ 
satisfies $a_\nu(x)>0$, for all $x\in [0,1]$ and $\nu > 0$,  it is small in a neighborhood of a point $x_*\in (0,1)$ and it is large outside this neighborhood for sufficiently small $\nu>0$, according to \cite{Carvalho-Pereira-94}.
More precisely, 
\begin{equation}\label{eq-applications-function-a}
\begin{array}{l l} 
a_\nu(x)\geqslant \frac{1}{\nu}, & \ x\in (0,x_*-\nu \beta_\nu)\cup (x_*+\nu \beta_\nu,1),\\
a_\nu(x)\geqslant \nu \alpha_0, & \ x\in(x_*-\nu \beta_\nu,x_*+\nu \beta_\nu),\\
a_\nu(x)\leqslant \nu \alpha_\nu, & \ x\in(x_*-\nu \beta_0,x_*+\nu \beta_0).
\end{array} 
\end{equation}
for some $x_*\in(0,1)$, $\alpha_\nu,\beta_\nu\in (0,\infty)$ satisfying $\alpha_\nu\searrow \alpha_0>0,\beta_\nu\searrow \beta_0>0$, as $\nu\to 0$.

The spectral properties of the operator $A_\nu u:=-(a_\nu(x) u_x)_x$, defined by $A_\nu:D(A_\nu)\subset L^2(0,1)\to L^2(0,1)$, where $D(A_\nu):=\{u\in H^2(0,1): u_x(0)=u_x(1)=0\}$ can be found in \cite[Lemma 1.1]{Carvalho-Pereira-94}. In particular, it was shown that the spectrum consists of a sequence of eigenvalues, $\lambda_1<\lambda_2<\lambda_3<\ldots$, with corresponding normalized eigenfunctions, $\phi_1,\phi_2,\phi_3,\ldots$, such that
\begin{equation}\label{ex:nonautgap}
	\lambda_1=0, \quad  \lambda_2\xrightarrow{\nu \to 0} \frac{\alpha_0}{2\beta_0}\frac{1}{(1-x_*)}, \quad \lambda_3\xrightarrow{\nu \to 0} +\infty	
\end{equation}
and 
\begin{equation}
\phi_1(x)\equiv 1, \, \phi_2(x)=
\left\{ 
\begin{array}{l l}
-k_1+O(\nu^{\epsilon}), & x\in [0,x_*-\nu \beta_0]\\
%-k_1+O(\nu^{q/2}), & \ x\in [x_*-\nu \beta_\nu,x_*-\nu \beta_0]\\
O(1),& x\in [x_*-\nu \beta_0,x_*+\nu \beta_0]\\
%k_1^{-1}+O(\nu^{q/2}), & \ x\in [x_*+\nu \beta_0,x_*+\nu \beta_\nu]\\
k_1^{-1}+O(\nu^{\epsilon}), & x\in [x_*+\nu \beta_0,1],
\end{array} 
\right.
\xrightarrow{\nu\to 0} 
\left\{ 
\begin{array}{l l}
-k_1, & \ x\in [0,x_*)\\
k_1^{-1}, & \ x\in (x_*,1],
\end{array} 
\right.
\end{equation}
for some $\epsilon \in (0,1)$ and $k_1:=\sqrt{(1-x_*)/x_*}.$

Therefore, \eqref{eq-applications-scalar-parabolic} is associated with a skew-product semiflow in the phase-space $X^\alpha \times \mathcal{B}$, where $X^\alpha$ is the fractional power spaces associated with $A$, for some $\alpha\in (0,1 )$, and $\mathcal{B}$ is the base space consisting of the closure of time translations of 
$\beta(\cdot)$, with respect to the metric of the uniform convergence in bounded intervals. Moreover, this skew-product semiflow is dissipative and hence has a global attractor $\mathbb{A}_\nu$, with associated uniform attractor $\mathcal{A}_\nu$ for all $\nu>0$, see \cite[Chapter 6]{BCL}.
Next, we present a result regarding the existence of non-trivial hyperbolic global solutions of \eqref{eq-applications-scalar-parabolic}.% , without perturbing a hyperbolic equilibrium of an autonomous equation.
%{\color{black}Our main result in this section is that the asymptotic profiles in the global attractors (i.e., $\alpha,\omega$-limit sets) are global hyperbolic non-autonomous solutions.}

\begin{corollary}\label{HypNonautEq}
For sufficiently small $\nu>0$, there are four non-trivial hyperbolic global solutions of  \eqref{eq-applications-scalar-parabolic}.
\end{corollary}
\proof
%Define $W= \text{span}\{\phi_1,\phi_2\}$, and $W_\alpha^\perp$ its orthogonal complement in $X^\alpha$. 
We rewrite solutions $u(t,x)$ of equation \eqref{eq-applications-scalar-parabolic} according to the eigenprojections in the first two eigendirections and a remainder term, as 
\begin{equation}
u(t,x)=u_1(t)\phi_1+u_2(t)\phi_2+w(t,x),
\end{equation}
%where $u_1,u_2$ are given by
%\begin{equation*}
%    \dot{u}_1(t)=\int_{0}^{1}f(t,u(t,y))dy \hbox{ and }
%    \dot{u}_2(t)=-\lambda_2u_2(t)+\int_{0}^{1}f(t,u(t,y))\phi_2(y)dy.
%\end{equation*}
which evolve according to the following system of differential equations
\begin{equation}\label{eq-applications-scalar-parabolic-complementary_equation}
%\left\{ 
\begin{split}
\dot{u}_1=& h_1(t,u_1,u_2,w), \,\, \qquad\qquad\qquad t>\tau, \\
\dot{u}_2=&-\lambda_2u_2+h_2(t,u_1,u_2,w),  \qquad t>\tau,  \\
w_t=&(a_\nu(x)w_x)_x+h_3(t,u_1,u_2,w), \  \ t>\tau,
\end{split} 
%\right.
\end{equation}
where $w$ satisfies Neumann boundary conditions, $w_x(t,0) = w_x(t,1)=0$ for all $t>\tau$, with initial condition $w(\tau,x)=u_0(x)-u_1(\tau,x)\phi_1(x)-u_2(\tau,x)\phi_2(x)$, and the vector field is given by
\begin{equation}
\begin{split}
h_1(t,u_1,u_2,w) & :=\int_{0}^{1}f(t,u_1(t)\phi_1(y)+u_2(t)\phi_2(y)+w(t,y))dy,  \\
h_2(t,u_1,u_2,w) & :=\int_{0}^{1}f(t,u_1(t)\phi_1(y)+u_2(t)\phi_2(y)+w(t,y))\phi_2(y)dy,  \\
h_3(t,u_1,u_2,w) & :=f(t,u_1\phi_1+u_2\phi_2+w)-h_1(t,u_1,u_2,w)-h_2(t,u_1,u_2,w)\phi_2.
\end{split} 
\end{equation}

As a consequence of Theorem \ref{th-existence-local-inertial-manifolds-skew-products}, we obtain the non-autonomous version of the inertial manifold reduction in \cite[Theorem 1.2]{Carvalho-Pereira-94}, due to the arbitrarily large gap $\gamma-\rho:=\lambda_3-\lambda_2$ in \eqref{ex:nonautgap} for sufficiently small $\nu>0$. Note that we obtain a two-dimensional local inertial manifold $\{\mathcal{M}_{\text{loc}}(t): t\in \mathbb{R}\}$ according to \eqref{loc:graph} for each sufficiently small $\nu>0$, which is a graph over $Im(Q(t))=\text{span}\{\phi_1,\phi_2\}$. 
Moreover, the graph describing the local inertial manifold satisfies $\Sigma_\nu\to 0$, as $\nu \to 0$, in $C^1$-norm, in a similar manner as in \cite{Carvalho-Pereira-94}, due to the properties of the inertial manifold in Theorem \ref{NLTum}. Thus it is sufficient to analyze the corresponding two-dimensional vector field projected to $\text{span}\{\phi_1,\phi_2\}$ for sufficiently small $\nu>0$, as the remaining of the evolution of the solutions can be obtained in terms of the graph.
Furthermore, the solution in the inertial manifold is 
given by 
$u(t,x)=u_1(t)+u_2(t)\phi_2(x)+\Sigma_\nu(t,u_1(t)+u_2(t)\phi_2)$, where 
$(u_1(t),u_2(t))$ is the solution of 
\begin{equation}\label{eq-equation-on-the-inertial-manifold}
\begin{split} 
\dot{u}_1(t)&=h_1(t,u_1,u_2,\Sigma_\nu(t,u_1,u_2)),\\
\dot{u}_2(t)&=-\lambda_2u_2+ h_2(t,u_1,u_2,\Sigma_\nu(t,u_1,u_2)).
\end{split} 
\end{equation}

%\begin{theorem}\label{th-applications-reduce-to-2-dim}
%Then, for any neighborhood $U$ of $\mathcal{A}$ such that %$\overline{U}\subset V$, there exists $\nu_0>0$ and a $C^k$ function 
%	$\Sigma_\nu:\mathbb{R}\times V\to X$, 
%	for each $\nu \leqslant \nu_0$, such that 
%	$\sup_{t\in \mathbb{R}}\Sigma(t,u)\to 0,$ for $u\in \overline{U}$, and the sets
%	\begin{equation*}
%	    \mathcal{M}_{\text{loc},\nu}(t)=\{v+\Sigma_{\text{loc},\nu}(t,v):v\in U\}, t\in \mathbb{R},
%	\end{equation*}
%defines an exponentially attracting inertial manifold $\widehat{\mathcal{M}}_\nu=\{\mathcal{M}_\nu(t):t\in \mathbb{R}\}$.
%\par Furthermore, the solution in $\widehat{\mathcal{M}}_\nu$ is  given by  $u(t,x)=u_1(t)+u_2(t)\phi_2(x)+\Sigma_\nu(t,u_1(t)+u_2(t)\phi_2)$, where  $(u_1(t),u_2(t))$ is the solution of 
%\begin{equation}\label{eq-equation-on-the-inertial-manifold}
%\begin{split} 
%\dot{u}_1(t)&=h_1(t,u_1,u_2,\Sigma_\nu(t,u_1,u_2))\\
%\dot{u}_2(t)&=-\lambda_2u_2+ h_2(t,u_1,u_2,\Sigma_\nu(t,u_1,u_2)).
%\end{split} 
%\end{equation}
%\end{theorem}
%\proof
%    \textcolor{lightgreen}{Must be a consequence of Theorem \ref{th-existence-local-inertial-manifolds-skew-products} in the parabolic case. Write something about it.}
%\cqd
%
\par Since $\Sigma_\nu\to 0$, as $\nu \to 0$, the limiting ordinary differential equation, as $\nu\to 0$, is given by
\begin{equation}\label{eq-limiting-problem-variables-v}
\begin{split}
\dot{u}_1(t) & =f_1(t,u_1,u_2),\\
\dot{u}_2(t) & =- \frac{\alpha_0}{2\beta_0}\frac{1}{(1-x_*)} u_2+ f_2(t,u_1,u_2),
\end{split} 
\end{equation}
where 
\begin{equation}
\begin{split} 
f_1(t,u_1,u_2)&:=x_*f(t,u_1-k_1u_2)+(1-x_*)f(t,u_1+k_1^{-1}u_2),\\
f_2(t,u_1,u_2)&:=-x_*f(t,u_1-k_1u_2)+(1-x_*)k_1^{-1}f(t,u_1+k_1^{-1}u_2).
\end{split} 
\end{equation}
Changing the variables, $z_1(t):=u_1(t)-k_1u_2(t) \hbox{ and } z_2(t):=u_1(t)+k_1^{-1}u_2(t)$, the new variables $z_1,z_2$ evolve according to the following ODE
%\begin{equation}\label{eq-change-variables-z-to-v}
%    \begin{split}
%        u_1(t)=&x_*z_1(t)+ (1-x_*)z_2(t),\\
%        u_2(t)=&-k_1x_*z_1(t) + k_2 (1-x_*) z_2(t), 
%    \end{split}
%\end{equation}
%\begin{equation*}
%\begin{split}
%z_1(t)=\frac{1}{x_*} \int_0^{x_*} u(t,x)dx, \\
%z_2(t)=\frac{1}{1-x_*} \int_{x_*}^1 u(t,x)dx,
%\end{split} 
%\end{equation*}
%and Problem \eqref{eq-limiting-problem-variables-v} becomes
\begin{equation}\label{eq-limiting-problem-variables-z}
\begin{split} 
\dot{z}_1=&\frac{\alpha_0}{2\beta_0 x_*} (z_2-z_1)+f(t,z_1)\\
\dot{z}_2=&\frac{-\alpha_0}{2\beta_0 (1-x_*)} (z_2-z_1)+f(t,z_2)
\end{split} 
\end{equation}

%Now, we obtain hyperbolic solutions for Problem \eqref{eq-applications-scalar-parabolic} for every $\nu>0$ suitable small, for a special nonlinearity $f$.

\par For $x_*=1/2$ and $\alpha_0/\beta_0 \in (\frac13,\frac12)$, the 2-dimensional ODE \eqref{eq-limiting-problem-variables-z} is considered in \cite{Carvalho-Langa-Obaya-Rocha}.
Taking advantage of the invariant subsets 
$E_1=\{(z_1,z_2): z_1=z_2\}$ and $E_2=\{(z_1,z_2): z_1=-z_2\}$, % with respective ODEs in each invariant set 
%\begin{subequations}
%\begin{align}
%&\dot{e}_1=e_1-\beta(t)e_1^3, \qquad \qquad  t> \tau, \label{eq-E_2} \\
%&\dot{e}_2=(1-2k)e_2-\beta(t) e_2^3, \ \ t> \tau \label{eq-E_1},
%\end{align}
%\end{subequations}
the authors in \cite[Theorem 3.4]{Carvalho-Langa-Obaya-Rocha} prove that there are exactly four non-autonomous global solutions which remain away from zero hyperbolic solutions $\xi_{1,\pm}$ and $\xi_{2,\pm}$ of \eqref{eq-limiting-problem-variables-z}, corresponding to the asymptotic profiles for the longtime behaviour. %Indeed, in \cite[Theorem 3.3]{Carvalho-Langa-Obaya-Rocha}, they prove that there are global solutions of \eqref{eq-E_1} and \eqref{eq-E_2} which both remain away from zero, given by $\xi_0$ and $\xi_k$, respectively.
%Moreover, in, the authors prove the following result on the existence of hyperbolic equilibria for \eqref{eq-limiting-problem-variables-z}.
%\begin{theorem}
%    Let $k\in (1/3,1/2)$. 
%    Then \eqref{eq-limiting-problem-variables-z} has exactly four global solutions which remains way from zero and are given by
%    \begin{equation*}
%        \begin{split}
%            &\xi^*_{1,+}=(\xi_0(t),\xi_0 (t)),  \ \ t\in \mathbb{R}\\
%            &\xi^*_{1,-}=(-\xi_0(t),-\xi_0 (t)),  \ \ t\in \mathbb{R}\\
%            &\xi^*_{2,+}=(\xi_k(t),-\xi_k (t)),  \ \ t\in \mathbb{R}\\
%            &\xi^*_{2,-}=(-\xi_k(t),\xi_k (t)),  \ \ t\in \mathbb{R}.
%        \end{split}
%    \end{equation*}
%    Moreover, $\xi^*_0=(0,0)$, $\xi^*_{1,\pm}$ and $\xi^*_{2,\pm}$ are hyperbolic equilibria for $\{S(t,s): t\geqslant s\}$, i.e., the corresponding linearizations 
    %of $\{S(t,s): t\geqslant s\}$
%    around each of these solutions admit an exponential dichotomy.
%\end{theorem}

%{\color{lightgreen} ADD FIFURE??}

Since we now know there are four non-trivial hyperbolic global non-autonomous solutions of equation \eqref{eq-limiting-problem-variables-z} in the variables $(z_1,z_2)$, we can reverse gear and obtain corresponding hyperbolic solutions for \eqref{eq-limiting-problem-variables-v} in the variables $(u_1,u_2)$.
Moreover, since $\Sigma_\nu\to 0$ as $\nu \to 0$, in $C^1$-norm uniformly in $t\in \mathbb{R}$, we apply the persistence of hyperbolic solutions in \cite[Lemma 8.3]{Carvalho-Langa-Robinson-13} for each non-tirvial hyperbolic solution of the limiting case $\nu=0$ in \eqref{eq-limiting-problem-variables-v}, which ensures the existence of a hyperbolic solution of \eqref{eq-equation-on-the-inertial-manifold} for any sufficiently small $\nu>0$. Therefore, we obtain the four hyperbolic global solutions of \eqref{eq-applications-scalar-parabolic} for any sufficiently small $\nu>0$.

%{\color{black}PL: Moreover, we can ascend the solutions from the limiting case $\nu=0$ equation \eqref{eq-limiting-problem-variables-v} to sufficiently small $\nu>0$ in \eqref{eq-equation-on-the-inertial-manifold}, since $\Sigma_\nu\to 0$, as $\nu \to 0$, in $C^1$-norm. Lastly, one can obtain the four hyperbolic global solutions of \eqref{eq-applications-scalar-parabolic} for sufficiently small $\nu>0$ by means of the graph $\Sigma_\nu$.}
%
%, if $v^*$ is a hyperbolic equilibrium for \eqref{eq-limiting-problem-variables-v}, then there exist $\nu_0>0$ such that, for each $\nu\in (0,\nu_0)$ there exists $\xi^*_\nu$ hyperbolic equilibrium for \eqref{eq-equation-on-the-inertial-manifold}.
%Note that, if $z^*=(z^*_1,z^*_2)$ is a hyperbolic equilibrium for  $\{S(t,s): t\geqslant s\}$. Then, from \eqref{eq-change-variables-z-to-v}, the global solution $v^*=(v^*_1,v^*_2)$ given by
%\begin{equation}
%\begin{split}
%    &u_1^*(t)=\frac{1}{2}(z^*_1(t)+z^*_2(t)), \ \ u_2^*(t)=\frac{1}{2}(z^*_2(t)-z^*_1(t)),
%    \end{split}
%\end{equation}
%is a hyperbolic equilibrium for \eqref{eq-limiting-problem-variables-v}. 
\cqd
\begin{remark}
The family of  pullback attractors 
$\{\mathcal{A}_\nu(t)_{ t\in \mathbb{R}}: \nu \geqslant 0\}$ of \eqref{eq-applications-scalar-parabolic} and  \eqref{eq-limiting-problem-variables-v} are continuous at $\nu =0$, i.e, $d_H (\mathcal{A}_\nu(t),\mathcal{A}_0(t))\to 0$ as $\nu \to 0$, see \cite[Theorem 7.1]{Carvalho-Langa-07}. %, i.e.,
%\begin{equation*}
%    \lim_{\nu \to 0} \sup_{t\in \mathbb{R}} 
%    d_H( \mathcal{A}_\nu(t),\mathcal{A}_0(t))=0,
%\end{equation*}
%where $d_H$is the Hausdorff distance: $d_H(A,B)=dist_H(A,B)+dist_H(B,A)$, $A,B\subset X$.
Moreover, for $x_*=1/2$, $\alpha_0/\beta_0\in (1/3,1/2)$ and $\beta:\R\to [1,2]$, the skew product associated to \eqref{eq-applications-scalar-parabolic} is dynamically gradient for suitably small $\nu>0$, since this is a small perturbation of a dynamically gradient system, see \cite{Carvalho-Langa-07,Carvalho-Langa-Obaya-Rocha}. 
\end{remark}

\end{document}